\documentclass[11pt,twoside]{article}

\usepackage{latexsym}
\usepackage{amsmath}
\usepackage{amsthm}
\usepackage{amssymb}
\usepackage{vmargin}
\usepackage{amscd}
\usepackage{stmaryrd}
\usepackage{euscript}
\usepackage{mathrsfs}
\usepackage{amscd}
\usepackage[all]{xy}
\usepackage{xr}
\DeclareMathAlphabet{\mathpzc}{OT1}{pzc}{m}{it}

\setmargins{32mm}{20mm}{14.6cm}{22cm}{1cm}{1cm}{1cm}{1cm}

\newcommand\da{\!\downarrow\!}
\newcommand\ra{\rightarrow}
\newcommand\la{\leftarrow}
\newcommand\lra{\longrightarrow}
\newcommand\lla{\longleftarrow}
\newcommand\id{\mathrm{id}}

\newcommand\ten{\otimes}
\newcommand\vareps{\varepsilon}
\newcommand\eps{\epsilon}

\newcommand\CC{\mathrm{C}}

\newcommand\Ru{\mathrm{R_u}}

\renewcommand\H{\mathrm{H}}
\newcommand\z{\mathrm{Z}}

\newcommand\N{\mathbb{N}}
\newcommand\Z{\mathbb{Z}}
\newcommand\Q{\mathbb{Q}}
\newcommand\R{\mathbb{R}}
\newcommand\Cx{\mathbb{C}}

\newcommand\vv{\mathbb{V}}
\newcommand\ww{\mathbb{W}}

\newcommand\bG{\mathbb{G}}
\newcommand\bH{\mathbb{H}}
\newcommand\bI{\mathbb{I}}

\newcommand\bL{\mathbb{L}}

\newcommand\bO{\mathbb{O}}

\newcommand\bS{\mathbb{S}}

\newcommand\C{\mathcal{C}}

\newcommand\cE{\mathcal{E}}

\newcommand\cG{\mathcal{G}}

\newcommand\cL{\mathcal{L}}
\newcommand\cM{\mathcal{M}}
\newcommand\cN{\mathcal{N}}
\newcommand\cO{\mathcal{O}}

\newcommand\cU{\mathcal{U}}

\newcommand\cW{\mathcal{W}}

\newcommand\sA{\mathscr{A}}

\newcommand\sF{\mathscr{F}}

\newcommand\m{\mathfrak{m}}
\newcommand\n{\mathfrak{n}}

\newcommand\g{\mathfrak{g}}

\newcommand\fh{\mathfrak{h}}

\newcommand\fk{\mathfrak{k}}

\newcommand\fr{\mathfrak{r}}

\newcommand\Ho{\mathrm{Ho}}

\newcommand\Alg{\mathrm{Alg}}
\newcommand\Mod{\mathrm{Mod}}

\newcommand\Hom{\mathrm{Hom}}
\newcommand\HHom{\underline{\mathrm{Hom}}}

\newcommand\End{\mathrm{End}}

\newcommand\Der{\mathrm{Der}}

\newcommand\Aut{\mathrm{Aut}}
\newcommand\Out{\mathrm{Out}}
\newcommand\ROut{\mathrm{ROut}}

\newcommand\Iso{\mathrm{Iso}}

\newcommand\im{\mathrm{Im\,}}

\newcommand\Ob{\mathrm{Ob}\,}

\newcommand\Top{\mathrm{Top}}
\newcommand\Gp{\mathrm{Gp}}
\newcommand\agp{\mathrm{AGp}}
\newcommand\agpd{\mathrm{AGpd}}

\newcommand\mal{\mathrm{Mal}}

\newcommand\Spec{\mathrm{Spec}\,}

\newcommand\Set{\mathrm{Set}}

\newcommand\Aff{\mathrm{Aff}}

\newcommand\Sch{\mathrm{Sch}}
\newcommand\sch{\mathrm{sch}}

\newcommand\Sing{\mathrm{Sing}}
\newcommand\FD{\mathrm{FD}}

\newcommand\ad{\mathrm{ad}}

\newcommand\Lim{\varprojlim}
\DeclareMathOperator*{\holim}{holim}
\newcommand\into{\hookrightarrow}
\newcommand\onto{\twoheadrightarrow}
\newcommand\abuts{\implies}
\newcommand\xra{\xrightarrow}

\newcommand\pr{\mathrm{pr}}

\newcommand\alg{\mathrm{alg}}

\newcommand\bt{\bullet}
\newcommand\by{\times}

\newcommand\mc{\mathrm{MC}}

\newcommand\Gg{\mathrm{Gg}}

\newcommand\Vect{\mathrm{Vect}}
\newcommand\Rep{\mathrm{Rep}}

\newcommand\Symm{\mathrm{Symm}}

\newcommand\Tot{\mathrm{Tot}\,}
\newcommand\diag{\mathrm{diag}\,}

\newcommand\ind{\mathrm{ind}}
\newcommand\pro{\mathrm{pro}}

\newcommand\pd{\partial}
\newcommand\dc{d^{\mathrm{c}}}
\newcommand\half{\frac{1}{2}}

\newcommand\gr{\mathrm{gr}}
\newcommand\ab{\mathrm{ab}}

\newcommand\gp{\mathrm{Gp}}
\newcommand\gpd{\mathrm{Gpd}}

\renewcommand\alg{\mathrm{alg}}
\newcommand\red{\mathrm{red}}
\newcommand\Fr{\mathrm{Fr}}
\newcommand\Lie{\mathrm{Lie}}

\newcommand\sk{\mathrm{sk}}

\newcommand\opp{\mathrm{opp}}

\newcommand\co{\colon\thinspace}

\newtheorem{theorem}{Theorem}[section]
\newtheorem{proposition}[theorem]{Proposition}
\newtheorem{corollary}[theorem]{Corollary}

\newtheorem{lemma}[theorem]{Lemma}
\newtheorem{theorem*}{Theorem}
\newtheorem{proposition*}[theorem*]{Proposition}
\newtheorem{corollary*}[theorem*]{Corollary}
\newtheorem{lemma*}[theorem*]{Lemma}

\theoremstyle{definition}
\newtheorem{definition}[theorem]{Definition}

\newtheorem{definition*}[theorem*]{Definition}

\theoremstyle{remark}
\newtheorem{example}[theorem]{Example}
\newtheorem{examples}[theorem]{Examples}
\newtheorem{remark}[theorem]{Remark}
\newtheorem{remarks}[theorem]{Remarks}

\newtheorem{exercises}[theorem]{Exercises}

\newtheorem{example*}[theorem*]{Example}
\newtheorem{examples*}[theorem*]{Examples}
\newtheorem{remark*}[theorem*]{Remark}
\newtheorem{remarks*}[theorem*]{Remarks}
\newtheorem{exercise*}[theorem*]{Exercise}

\begin{document}
\title{Pro-algebraic homotopy types}
\author{J.P.Pridham\thanks{Trinity College, Cambridge, CB2 1TQ, U.K. \newline
The author is supported by Trinity College, Cambridge.\newline
}}
\maketitle

\begin{abstract}
The purpose of this paper is to generalise Sullivan's rational homotopy theory to non-nilpotent spaces, providing an alternative approach to defining To\"en's schematic homotopy types over any field $k$ of characteristic zero. New features include an explicit description of homotopy groups using the Maurer-Cartan equations,  convergent spectral sequences comparing schematic homotopy groups with cohomology of the universal semisimple local system, and a generalisation of the Baues-Lemaire conjecture. For compact K\"ahler manifolds, the schematic homotopy groups can be described explicitly in terms of this cohomology ring, giving them canonical weight decompositions.  There are also notions of minimal models, unpointed homotopy types and algebraic automorphism groups. For a space with algebraically good fundamental group and higher homotopy groups of finite rank, the schematic homotopy groups are shown to be $\pi_n(X)\ten_{\Z}k$.
\end{abstract}

%55p15, 55p62, 14l17, 32j27
%
\tableofcontents

\section*{Introduction}
\addcontentsline{toc}{section}{Introduction}

Given a group $G$ and a topological space $X$, the  non-abelian cohomology set $\H^1(X,G)$ is the set of isomorphism classes of $G$-torsors on $X$, or equivalently the set $[X,BG]$ of homotopy classes of maps from $X$ to the classifying space $BG$. There is a notion of classification space for every simplicial group. This has the property that for simplicial abelian groups $A_{\bt}$,
$$
[X,BA_{\bt}]=\bH^1(X,A^{\bt}),
$$
where  $A^{\bt}$ is $A$ regarded as a negatively graded cochain complex (i.e. $A^{-n}:=A_n$).

One of several equivalent definitions of higher non-abelian cohomology is to consider
$$
\H^1(X,G_{\bt}):=[X,BG_{\bt}],
$$
where $G_{\bt}$ ranges over all (non-abelian) simplicial groups.
In \cite{loopgp}, Kan constructed a loop group functor $G(X)$ from reduced simplicial sets (or equivalently pointed, connected topological spaces) to simplicial groups. This is left adjoint to the classification space functor, so studying non-abelian cohomology of a space is equivalent to studying its loop group. The geometric realisation of the loop group is weakly equivalent to the loop space of $X$, so $\pi_iG(X)=\pi_{i+1}X$.

In this paper, we only wish to study non-abelian cohomology over a field $k$ of characteristic $0$. In the abelian case, this would mean restricting attention to   (finite-dimensional) simplicial $k$-vector spaces. Since these are just  simplicial abelian unipotent algebraic $k$-groups,  for the non-abelian case we generalise to  simplicial complexes $G_{\bt}$ of algebraic groups over $k$. As we wish homotopy groups to be $k$-linearised, we only consider those $G_{\bt}$ for which $G_n \to \pi_0(G)$ is a unipotent extension; this is a crucial difference from the approach of \cite{schematic}. We refer to these as algebraic simplicial groups. 

This enables us to develop an analogue of Sullivan's and Quillen's rational homotopy theories for non-nilpotent spaces. This differs from the homotopy types of \cite{GHT} in that we only consider finite-dimensional local systems. 
Whereas the schematic homotopy types of \cite{chaff} can be thought of as a natural generalisation of Sullivan's theory (\cite{Sullivan}), 
the approach we use here is closer in spirit to the rational homotopy type defined by Quillen in \cite{QRat}, studying simplicial groups and their torsors. This leads to a new proof of the Baues-Lemaire conjecture (Remark \ref{baues}), comparing Quillen and Sullivan rational homotopy types.

Although pro-algebraic and schematic homotopy types are equivalent in a weak sense (Corollary \ref{eqtoen}), this new perspective leads to the discovery of much extra structure, especially for pro-algebraic homotopy groups. There are convergent Adams and reverse Adams spectral sequences (Propositions \ref{spectralh} and \ref{spectralpi}) relating cohomology and homotopy groups. In many cases, the pro-algebraic homotopy groups are just $\pi_n(X)\ten_{\Z}k$ (Theorems \ref{classicalpi} and \ref{classicalpimal}). We also show that pro-algebraic homotopy groups can be described explicitly in terms of differential forms (Corollary \ref{bigequiv}), taking the form conjectured by Deligne (Remark \ref{deligne}). For compact K\"ahler manifolds, homotopy groups can thus be described explicitly in terms of cohomology (Proposition \ref{formalpins}).  

Other important new features are a characterisation of unpointed homotopy types (\S \ref{unptd}), and  good notions of minimal models (Propositions \ref{sminimal} and \ref{dgminimal}). The latter enable us (Theorem \ref{auto}) to regard the automorphism group of a pro-algebraic homotopy type as a pro-algebraic group itself. It then makes sense to talk about a weight decomposition on the homotopy type as being a homomorphism from $\bG_m$ to the automorphism group. This has several consequences in \S \ref{hodge}, including canonical weight decompositions for compact K\"ahler manifolds, and lays a foundation for the results of \cite{hodge}.

In Section \ref{ptd}, we define the pro-algebraic homotopy type of a reduced simplicial set $X$ to be the completion $G(X)^{\alg}$  of $G(X)$ with respect to algebraic simplicial groups. In other words, it is the pro-algebraic simplicial group given by
$$
\Hom(G(X)^{\alg},G)=\Hom(G(X),G(k)),
$$
for all algebraic simplicial groups $G$.  This approach is similar to that in \cite{QRat}, where Quillen took a  pro-nilpotent (although not pro-finite-dimensional) completion of the loop group. 

We prove that the category of pro-algebraic simplicial groups has a closed model structure, and that $\pi_0(G^{\alg})=\pi_0(G)^{\alg}$. There is also a homology theory for pro-algebraic simplicial groups, with convergent Adams and reverse Adams spectral sequences between homotopy and homology. This enables us to show that if $\pi_1X$ is algebraically good, and the higher homotopy groups have finite rank, then the higher pro-algebraic homotopy groups are just $\pi_nX\ten_{\Z}k$. For schematic homotopy groups, this was previously only known for compact K\"ahler manifolds.

In Section \ref{unptd}, we extend these notions to unpointed homotopy types. This is done by replacing groups by groupoids. We begin by defining a pro-algebraic groupoid to consist of a discrete set of objects, together with affine schemes of morphisms between them, endowed with the usual  composition, identity and inversion maps. Thus a pro-algebraic group is a pro-algebraic groupoid on one object. All the important properties of pro-algebraic groups, such as Levi decompositions and Tannakian duality, carry over to pro-algebraic groupoids. The results of Section \ref{ptd} are then adapted by substituting Dwyer and Kan's path groupoid $G(X)$ for the  loop group. A closed model structure is defined on the category of pro-algebraic simplicial groupoids, and the unpointed pro-algebraic homotopy type is defined as the pro-algebraic completion $G(X)^{\alg}$. 

In Section \ref{mc}, we define the simplicial Maurer-Cartan space, which classifies higher torsors, and the action of the gauge group  on it, which detects equivalences. An alternative characterisation of the pro-algebraic homotopy type is derived by first taking a Levi decomposition $G(X)^{\alg}=R\ltimes U$, with $R$ a reductive pro-algebraic groupoid, and $U$ an $R$-representation in  simplicial pro-unipotent groups. Then the classification space $BU$ is an $R$-representation in simplicial affine schemes. This gives us an equivalence between the homotopy categories of simplicial pro-unipotent extensions of $R$, and of connected  $R$-representations in simplicial affine schemes; we can think of this as a form of Koszul duality (Remark \ref{koszul}). This fulfils the hope expressed in \cite{schematic}v1 Remark 3.19 of comparing equivariant cosimplicial algebras and simplicial groups  directly, without needing simplicial presheaves to mediate between them.

By studying the Maurer-Cartan space and its gauge action, we then show that there is a weak equivalence
$$
\Spec \CC^{\bt}(X,\bO(R)) \sim BU,
$$
where $\bO(R)$ is the universal semisimple local system on $X$, and $\CC^{\bt}$ denotes simplicial cochains. This allows us to prove an equivalence between pro-algebraic homotopy types and To\"en's schematic homotopy types.

In Section \ref{forms}, we show how to replace simplicial objects by differential graded objects. Simplicial pro-unipotent groups correspond to chain Lie algebras, and simplicial affine schemes correspond contravariantly to cochain algebras. If $X$ is a manifold and $k=\R$ or $\Cx$, this allows us to replace $\CC^{\bt}(X,\bO(R))$ by the cochain algebra
$$
A^{\bt}(X,\bO(R) ):=\Gamma(X, \bO(R)\ten\sA^{\bt}),
$$
where $\sA^{\bt}$ is the complex of differential forms on $X$. There is a notion of minimal models for chain Lie algebras, which is often more convenient than that of  Sullivan's minimal models for cochain algebras, being generated by homology rather than homotopy groups. 

Under this correspondence, the simplicial Maurer-Cartan space becomes the familiar classical space, given by the equation
$$
d\omega + \half[\omega,\omega]=0.
$$
The gauge action also corresponds to the classical gauge action
$$
g(\omega)=   g\cdot \omega \cdot g^{-1} -dg\cdot g^{-1}.
$$
These allow us to recover the pro-algebraic  homotopy groups and Whitehead products explicitly from differential forms. Indeed, this confirms that schematic homotopy groups are dual to generators of the minimal model of $A^{\bt}(X,\bO(R))$, as originally conjectured by  Deligne.

In Section \ref{hodge}, we show that  for finite simplicial complexes, the  automorphism group of $G(X)^{\alg}$ over its reductive quotient $\pi_fX^{\red}$ is itself pro-algebraic. 
We then  investigate formal topological spaces, i.e. those for which the cochain algebra $A^{\bt}(X,\bO(R) )$ is weakly equivalent to its cohomology ring. Compact K\"ahler manifolds are formal.   The constructions of \S \ref{forms} then enable us to  describe pro-algebraic homotopy groups of a formal space in terms of $\H^*(X,\bO(R))$. The formality weak equivalence gives  a canonical weight decomposition on the pro-algebraic homotopy type of a formal space, in the form of a morphism $\bG_m \to \Aut(G(X)^{\alg})_{\pi_fX^{\red}}$.
We use this in Corollary \ref{ratformal} to show that, for a finite simplicial complex, formality over $\R$ or $\Cx$ implies formality over all fields (of characteristic zero). 

Although analogues of several of our constructions have appeared in \cite{schematic}v1, all  the main results of this paper are new, even for schematic homotopy types, and simplicial presheaves play no r\^ole in our proofs.
Remark \ref{modelbad} provides a 
discussion comparing $s\agp$ and the cosimplicial Hopf algebras of [ibid.]v2.
Our results rely heavily on unipotence and Koszul
duality, so most are unlikely to be demonstrable by adapting current machinery from
schematic homotopy theory. However, as explained in Remark \ref{toenadams}, our Adams spectral sequence is essentially the weight spectral sequence of \cite{KTP}; the latter was an early demonstration of how powerful a tool relative unipotent completion can be.

I would like to thank Bertrand To\"en for alerting me to errors in a previous version. 

\section{Pointed pro-algebraic homotopy types}\label{ptd}

We begin by recalling some standard definitions from \cite{sht}.

\begin{definition}
Let $\bS_0$ be the category of reduced simplicial sets, i.e. simplicial sets with one vertex, and $s\Gp$ the category of simplicial groups. Let $\Top_0$ denote the category of pointed, connected compactly generated Hausdorff topological spaces.
\end{definition}

Note that there is a functor from $\Top_0$ to $\bS_0$ which sends $(X,x)$ to the simplicial set
$$
\Sing(X,x)_n:=\{ f \in \Hom_{\Top}(|\Delta^n|, X) : f(v)=x \quad \forall v \in \Delta^n_0\}.
$$
this is a right Quillen equivalence, the corresponding left equivalence being geometric realisation. For the rest of this section, we will therefore restrict our attention to reduced simplicial sets.

 The idea underlying this section is that for any group $G$, and any topological space $X$, $\H^1(X,G)=\pi_0 \HHom(X,BG)$. Given a simplicial group $G_{\bullet}$, we can extend this definition to:
$$
\H^1(X,G_{\bullet}):=\pi_0 \HHom(X,BG_{\bullet}).
$$

As in \cite{sht} Ch.V, the 
the classification space $BG_{\bullet}$ can be chosen canonically by the functor $\bar{W}\co s\Gp \to \bS_0$. This has a left adjoint $G\co \bS_0 \to s\Gp$, Kan's loop group functor (\cite{loopgp}), and these give a Quillen equivalence of model categories. In particular, $\pi_i(G(X))=\pi_{i+1}(X)$.

\subsection{Review of pro-algebraic groups}

\begin{definition}
Given a category $\C$, recall that the category $\pro\C$ of pro-objects in $\C$, which we will also denote by $\widehat{\C}$, has objects consisting of filtered inverse systems $\{A_{\alpha}\in \C\}$, with 
$$
\Hom_{\widehat{\C}}(\{A_{\alpha}\}, \{B_{\beta}\})= \lim_{\substack{\lla \\ \beta}} \lim_{\substack{\lra \\ \alpha}} \Hom_{\C}(A_{\alpha},B_{\beta}).
$$
\end{definition}

\begin{lemma}
 The category of  pro-algebraic groups over a field $k$ is equivalent to the category of  affine group schemes over $k$. 
\end{lemma}
\begin{proof}
\cite{tannaka} Corollary 2.7.
\end{proof}

\begin{definition}
If $G$ is a pro-algebraic group, let $O(G)$ denote the Hopf algebra of global sections of its structure sheaf.
\end{definition}

\begin{lemma}\label{grep}
The category $\Rep(G)$ of  $k$-representations of $G$ is equivalent to the category of $O(G)$-comodules. Every $G$-representation $V$ can be expressed as a directed union of finite-dimensional subrepresentations of $V$. Thus dualisation gives a contravariant equivalence of categories between $\Rep(G)$ and $\widehat{\FD\Rep}(G)$, where $\FD\Rep(G)$ is the category of finite-dimensional $G$-representations.
\end{lemma}
\begin{proof}
\cite{tannaka}   Proposition II.2.2 and Corollary 2.4.
\end{proof}

\begin{definition}
Given an abstract group $\Gamma$, the pro-algebraic completion $\Gamma^{\alg}$  of $\Gamma$ over $k$ is the pro-algebraic $k$-group representing the functor
$$
G \mapsto \Hom_{\gp}(\Gamma, G(k)),
$$ 
where $G$ ranges over all pro-algebraic groups. 
\end{definition}

\begin{remark}
Note that the category $\FD\Rep(\Gamma)$ of finite-dimensional $\Gamma$-representations is thus equivalent to $\FD\Rep(\Gamma^{\alg})$. In fact, Tannakian duality (\cite{tannaka} Theorem 2.11) shows that any proalgebraic group $G$ can be recovered from the tensor category $\FD\Rep(G)$. Thus $\Gamma^{\alg}$ can be recovered from $\FD\Rep(\Gamma)$. The key idea behind this paper is that algebraisation of groups is in many ways easier than the algebraisations of sets considered in \cite{chaff}.
\end{remark}

\subsection{Simplicial pro-unipotent groups}\label{unip} 
From now on, let  $k$ be a  field  of characteristic $0$. Given a reductive pro-algebraic group $G$ over $k$, let $\cN(G)$ be the category of finite-dimensional nilpotent Lie algebras equipped with $G$-actions (as in \cite{higgs}\S 2). Write $\hat{\cN}(G)$ for the category of pro-objects of $\cN(G)$, and $s\hat{\cN}(G)$ for the category of simplicial objects in $\hat{\cN}(G)$. Since any unipotent algebraic group is isomorphic to its Lie algebra, via exponentiation and the Campbell-Baker-Hausdorff formula, there is an equivalence between unipotent algebraic groups and finite-dimensional nilpotent Lie algebras. As $\cN(G)$ is an Artinian category, every pro-object is isomorphic to a strict pro-object (\cite{descent}). 

\begin{definition}\label{lowercentral}
Given $\g_{\bt} \in s\hat{\cN}(G)$ define the lower central series of $\g_{\bt}$ inductively by
$$
\Gamma_1\g_{\bt}=\g_{\bt}\quad \Gamma_{n+1}\g_{\bt}=[\g_{\bt},\Gamma_n\g_{\bt}].
$$
Denote the abelianisation of $\g$ by $\g^{\ab}:=\g/[\g,\g]$.
\end{definition}

\begin{definition}
Given $\g_{\bt} \in s\hat{\cN}(G)$, define the normalised complex $N(\g)_*$ by
$$
N(\g)_n= \bigcap_{i=1}^n \ker(\pd_i\co \g_n \to \g_{n-1}),
$$
and observe that $\pd_0\co N(\g)_n \to N(\g)_{n-1}$.  We then define 
$$
\pi_n(\g):= \H_n(N(\g)_*, \pd_0).
$$
\end{definition}

\begin{definition}\label{cfw}
 A morphism $f\co \g_{\bt} \to \fh_{\bt}$ in $s\hat{\cN}(G)$ is said to be:
\begin{enumerate}
\item a weak equivalence if $\pi_*(\g_{\bt}) \to \pi_*(\fh_{\bt})$ is an isomorphism of pro-finite-dimensional $G$-representations;
\item a fibration if  $N_n(f)\co  N(\g)_n \to N(\fh)_n$ is surjective (in $\widehat{\FD\Vect}$) for all $n>0$;
\item a cofibration if it has LLP with respect to all trivial fibrations  (a fibration or cofibration is called \emph{trivial} if it is also a weak equivalence).
\end{enumerate}
\end{definition}

\begin{definition}
Define $s\cN(G)$ to be the full subcategory of $s\hat{\cN}(G)$ whose objects $\g$ are finite-dimensional in the sense that $\g_n \in \cN(G)$ and  $N(\g)$ is concentrated in finitely many degrees. Note that the objects of $s\cN(G)$ are cofinite in $s\hat{\cN}(G)$ in the sense of \cite{Hovey} Definition 2.1.4.   

Define a small extension in $s\cN(G)$ to be a surjective map $\g \to \fh$ with kernel $I$, such that $[\g,I]=0$. 
\end{definition}

\begin{proposition}\label{cmsn}
  The classes of morphisms given in Definition \ref{cfw} define a closed model category structure on $s\hat{\cN}(G)$.
\end{proposition}
\begin{proof}
Apply \cite{Bou} Theorem 3.3 to $\hat{\cN}(G)^{\opp}$, with the class $\cG$ of injective models pro-representing the functors $\g \mapsto \Hom_{\pro-\FD\Vect}(V, \g)$, for each  pro-finite-dimensional vector space $V$. Note that this characterises cofibrations as maps $\g \to \fh$ for which the simplicial latching maps $(\g_n)*_{L_n \g} (L_n\fh)\to \fh_n$ are smooth in the sense of \cite{higgs}. This implies that cofibrant objects are levelwise smooth.

 Alternatively, we may use \cite{Hovey} Theorem 2.1.19 to show that this is a fibrantly cogenerated model category.  Acyclic small extensions in $s\cN(G)$ are the generating acyclic fibrations, while small extensions in $s\cN(G)$ together with arbitrary maps in $\cN(G)$ give the generating fibrations.
\end{proof}

The following can be thought of as an Adams-type spectral sequence, describing homotopy groups in terms of homology groups:
\begin{proposition}\label{spectralh1}
  Given $\g \in s\hat{\cN}(G)$ cofibrant, there is a convergent  spectral sequence (in $\widehat{\FD\Vect}$)
$$
E^1_{pq}= (\Lie_{-p} (\pi_*(\g^{\ab})))_{p+q} \abuts \pi_{p+q}(\g),
$$
where  $\Lie_nV = \Gamma_n \Lie(V)/\Gamma_{n+1}\Lie(V)$, for $\Lie$ the free graded Lie algebra functor. 
\end{proposition}
\begin{proof}
The spectral sequence comes from the  lower central series filtration of Definition \ref{lowercentral}, giving
$$
E^1_{pq}= \pi_{p+q}(\gr^{\Gamma}_{-p}\g ) \abuts \pi_{p+q}(\g).
$$
Now, $\gr^{\Gamma}_{-p}=\Lie_{-p}(\g^{\ab})$, as $\g$ is levelwise free. Similarly to \cite{QRat} Proposition I.4.5,  $\pi_*\Lie_{-p}(\g^{\ab})\cong \Lie_{-p}(\pi_*(\g^{\ab}))$, giving the spectral sequence.

We now need to show convergence.  There is a functor
$$
\varinjlim: \ind(\FD\Vect) \to \Vect
$$
from the  category of ind-finite-dimensional vector spaces to the category of all vector spaces, given by mapping a direct system $\{V_{\alpha}\}$ to $\varinjlim V_{\alpha}$. This is essentially surjective, since  any vector space is the  direct limit of its finite-dimensional subspaces. It is also full and faithful:
$$
\Hom_{\Vect}(\lim_{\substack{\lra\\ \alpha}} V_{\alpha}, \lim_{\substack{\lra\\ \beta}} W_{\beta})= \lim_{\substack{\lla\\ \alpha}}\Hom_{\Vect}(V_{\alpha}, \lim_{\substack{\lra\\ \beta}} W_{\beta})= \lim_{\substack{\lla\\ \alpha}}\lim_{\substack{\lra\\ \beta}}\Hom_{\Vect}(V_{\alpha},  W_{\beta}),
$$
since $V_{\alpha}$ is finite-dimensional. 

By taking duals, we see that $\Vect$ is equivalent to the opposite category of $\widehat{\FD\Vect}$. Convergence now follows by taking duals (replacing pro-finite-dimensional vector spaces by vector spaces), and observing that the dual filtration  is bounded below and exhaustive. The Classical Convergence Theorem (\cite{W} Theorem 5.5.1) then applies, since  the category of vector spaces satisfies AB5.  
\end{proof}

\begin{remark}\label{FDVspectral}
In general, this argument shows that for a increasing filtration $F_*$ on a complex in $\widehat{\FD\Vect}$, the associated spectral sequence will converge if $F$ is bounded above  and Hausdorff ($\bigcap F_p=0$). In the Proposition above, $F_p=\Gamma_{-p}$.
\end{remark}

\begin{proposition}\label{htpycoho}
If $\g,\fh$ are cofibrant objects in $s\hat{\cN}(G)$, then a morphism $f\co \g \to \fh$ is a weak equivalence if and only if $f^{\ab}\co  \g/[\g,\g] \to \fh/[\fh,\fh]$ is a quasi-isomorphism of simplicial vector spaces.
\end{proposition}
\begin{proof}
 For any G-representation $V$, let $V[-n]$ be the chain complex consisting of $V$ in degree $n$, so that $N^{-1}V[-n]$ is a simplicial (abelian) Lie algebra. If $f$ is a weak equivalence, then  the homotopy classes of maps from $\g$ and $\fh$ to  $N^{-1}V[-n]$ must be isomorphic, but these are precisely 
$$
\H^n(\Hom(\g/[\g,\g],V))  
$$
for $\g$, and similarly for $\fh$. Since $V$ was arbitrary, $f^{\ab}$ must be a quasi-isomorphism.

Conversely, since $\g$ is cofibrant, we have an Adams spectral sequence. Now, $f$ gives an isomorphism on the $E^1$ term of this spectral sequence, so we have an isomorphism at the limit. 
\end{proof}

The following proposition is an analogue of \cite{Q} Proposition 10.1.
\begin{proposition}[Minimal models]\label{sminimal}
Every weak equivalence class in $s\hat{\cN}(G)$ has a cofibrant element $\m$, unique up to non-unique isomorphism, with $d=0$ on  the  normalisation $N(\m/[\m,\m])$ of the  abelianisation.
\end{proposition}
\begin{proof}
Given such an equivalence class, first choose a cofibrant representative $\g_{\bt}$, with abelianisation  $V_{\bt}$. Since $\Rep(G)$ is a semisimple category, $G$ being reductive, we may decompose the simplicial  complex
$$
V_{\bt}=U_{\bt}\oplus W_{\bt},
$$ 
with $\H_*(U_{\bt})=0$ and $d=0$ on $N(W_{\bt})$, making use of the Dold-Kan correspondence.

We then have a commutative diagram 
$$
\begin{CD}
0 @>>> \g_{\bt}\\
@ViVV @VVpV\\
U_{\bt}@>>> V_{\bt}
\end{CD}
$$
in the model category of simplicial pro-finite-dimensional vector spaces, with $i$ a trivial cofibration, and $p$ a fibration. This allows us to lift $U_{\bt}$ to  a subcomplex of $\g_{\bt}$, which we denote by $\tilde{U}_{\bt}$. 

Now, we can also choose a lift $\tilde{W}_*$ of $W_{\bt}$ to $\g_{\bt}$, closed under degeneracy operations, but not necessarily face operations. Since $\g_{\bt}$ is cofibrant, it is generated by $\tilde{U}_*\oplus\tilde{V}_*$. 

Define the minimal model as the quotient
$$
\m_{\bt}:=\g_{\bt}/\langle \tilde{U}_{\bt} \rangle
$$
of $\g_{\bt}$ by the ideal generated by $\tilde{U}_{\bt}$. This Lie algebra is freely generated by the image of $\tilde{W}_*$, so is cofibrant. Since $\m_{\bt}^{\ab}\cong W_{\bt}$, the minimality criterion is satisfied.

Finally, if $\n_{\bt}$ is another minimal model in the same weak equivalence class, then there must exist a weak equivalence
$$
f\co \m_{\bt} \to \n_{\bt}, 
$$ 
since $\m_{\bt}$ is cofibrant. But
$$
f^{\ab}\co \m_{\bt}^{\ab} \to \n_{\bt}^{\ab}
$$
must then be an isomorphism, since $N(f^{\ab})$ is, $d$ being zero on both complexes. Since $\m_{\bt},\n_{\bt}$ are both free, the morphism $\gr f$ between the associated graded Lie algebras is an isomorphism, so $f$ must also be an isomorphism (or use \cite{higgs} Propositions 2.22 and 2.28).
\end{proof}

\begin{definition}
A morphism $f\co \g \to \fh$ is said to be smooth if it has the LLP with respect to all surjections. 
\end{definition}

\begin{lemma}\label{smoothcrit}
A morphism $f\co \g \to \fh$ is smooth if and only if it is cofibrant and the homology groups
$$
\H_i(\fh/([\fh,\fh]+\g))
$$
of the relative cotangent complex vanish for all $i>0$.
\end{lemma}
\begin{proof}
As for \cite{Q} Theorem 5.4.
\end{proof}

\subsection{Pro-algebraic simplicial groups}

\begin{definition}
Given a simplicial object $G_{\bullet}$ in the category of pro-algebraic groups, define $\pi_0(G_{\bullet})$ to be the coequaliser 
$$
\xymatrix@1{G_1 \ar@<1ex>[r]^{\pd_1} \ar@<-1ex>[r]_{\pd_0}& G_0 \ar[r] &\pi_0(G) }
$$ 
in the category of pro-algebraic groups.
\end{definition}

\begin{definition}
Define a pro-algebraic simplicial  group  to consist of a simplicial complex $G_{\bullet}$ of pro-algebraic groups, such that the maps $G_n \to \pi_0(G)$ are pro-unipotent extensions of pro-algebraic groups, i.e. $\ker(G_n \to \pi_0(G))$ is pro-unipotent.  We denote the category of pro-algebraic simplicial groups by $s\agp$. 
\end{definition}

\begin{remark}
If we instead took $k$ to be field of finite characteristic, it seems that some of the  results of this section  should still  hold. Levi decompositions are vital for \S\S \ref{mc}--\ref{hodge}  of this paper, so these sections cannot generalise to finite characteristic, although it is possible that an equivalence between the pro-algebraic and schematic homotopy types can be established by other means.
\end{remark}

\begin{definition}
Given a pro-algebraic simplicial group $G_{\bt}$, define the normalised complex $N(G)_*$ by
$$
N(G)_n= \bigcap_{i=1}^n \ker(\pd_i\co G_n \to G_{n-1}),
$$
and observe that $\pd_0\co N(G)_n \to N(G)_{n-1}$. These are pro-unipotent algebraic groups for $n>0$, and we define 
$$
\pi_n(G):= N(G)_n/\pd_0N(G)_{n+1},
$$
where the quotient is taken as pro-algebraic groups. Observe that this is consistent with the definition of $\pi_0$ given above, and that for $n>0$ the groups $\pi_n(G)$ are abelian pro-unipotent $k$-algebraic groups, i.e. pro-finite-dimensional vector spaces over $k$.
\end{definition}

\begin{definition}
Define a morphism $f\co G_{\bullet} \to H_{\bullet}$ in $s\agp$ to be: 
\begin{enumerate}
\item a weak equivalence if the maps $\pi_n(f)\co \pi_n(G_{\bullet}) \to \pi_n(H_{\bullet})$ are isomorphisms for all $n$; 

\item a fibration if $N_n(f)\co  N(G)_n \to N(H)_n$ is surjective  for all $n>0$;

\item a cofibration if  it has the left lifting property (LLP) with respect to all trivial fibrations.
\end{enumerate}
\end{definition}

\begin{remark}\label{modelbad}
In \cite{schematic}v2, a model category structure is being defined on the category 
$s\mathrm{GAff}$
of all simplicial objects in the
category of pro-algebraic groups, replacing the erroneous structure of [ibid.]v1. 
A map $f:G \to H$ in that category is said to be a weak equivalence  if it  induces an isomorphism 
on homotopy group schemes $\pi_nG$. The 
map $f$ is a fibration if for all trivial cofibrations $K \into L$ in $\bS$, the morphism
$$
\Hom_{\bS}(L,G)\to \Hom_{\bS}(K,G)\by_{\Hom_{\bS}(K,H)}\Hom_{\bS}(L,H)
$$ 
of group schemes is faithfully flat. The simplicial structure is given by defining $G^K$ as in Proposition \ref{sagpsimplicial}. 
It follows immediately that the inclusion functor $\iota:s\agp \to s\mathrm{GAff}$  is simplicial right Quillen, so gives a 
functor $\iota:\Ho(s\agp) \to \Ho(s\mathrm{GAff})$ on homotopy categories. 

Morphisms in $s\mathrm{GAff}$ satisfying the criteria of Corollary \ref{detectweak} are referred to in [ibid.] as $P$-equivalences, so every weak equivalence in $s\agp$ is a $P$-equivalence. The quotient category $\Ho^P(s\mathrm{GAff})$ is defined by localising at $P$-equivalences, and is shown to be a model for pointed schematic homotopy types. However, [ibid.] does not realise this category as the homotopy category of a model category, but conjectures that it should arise as a localisation of $s\mathrm{GAff}$. We now sketch such a construction.

Observe that $s\agp$ is fibrantly cogenerated. The cogenerating fibrations consist of arbitrary maps $G \to H$ of (constant) algebraic groups, together with maps $R \ltimes U' \to R \ltimes U$, for $R$ a reductive algebraic group and $U' \to U$ a small extension. The cogenerating trivial fibrations consist of maps $R \ltimes U' \to R \ltimes U$, for $R$ a reductive algebraic group and $U' \to U$ an acyclic small extension. The functor $\iota:s\agp \to s\mathrm{GAff}$ has a left adjoint $G \mapsto \hat{G}$, given by taking the pro-unipotent completion of $G \to \pi_0G$. This adjoint pair satisfies the conditions dual to \cite{Hirschhorn} Theorem 11.3.2, so gives a new model structure $s\mathrm{GAff}^P$ on $s\mathrm{GAff}$, for which $f$ is a cofibration or weak equivalence whenever $\hat{f}$ is.

Now, observe that $\iota:s\agp \to s\mathrm{GAff}^P$ is a simplicial right Quillen equivalence, that $\Ho(s\mathrm{GAff}^P)=\Ho^P(s\mathrm{GAff})$, and that the identity functor $ s\mathrm{GAff}^P \to s\mathrm{GAff}$ is simplicial right Quillen.

For computational purposes, the main disadvantage of $\Ho^P(s\mathrm{GAff})$ over $\Ho(s\agp)$ is that
the pro-algebraic groups $\pi_n(G_{\bullet})$ are not in general invariant under $P$-equivalence, so schematic homotopy groups have to be  defined more subtly.
\end{remark}

\begin{lemma}[Levi decomposition]
For $G_{\bt} \in s\agp$, there is a section of $G_{\bt} \to \pi_0(G)^{\red}$, unique up to conjugation by $\Ru(G_0)$, giving a decomposition
$$
G_{\bt}= \Ru(G)_{\bt} \rtimes \pi_0(G)^{\red}.
$$ 
\end{lemma}
\begin{proof}
Since $G_0 \to \pi_0(G)$ is a pro-unipotent extension, we have $G_0^{\red}=\pi_0(G)^{\red}$, so the Levi decomposition for $G_0$ (\cite{Levi}) gives a section of $G_0 \to \pi_0(G)^{\red}$. Composing with the degeneracy morphisms, this gives the required section to $G_{\bt}$. That we then have a Levi decomposition follows from the condition that $G_n \to \pi_0(G)$ is a pro-unipotent extension for all $n$.
\end{proof}

From now on, we denote $\pi_0(G)^{\red}$ simply by $G^{\red}$.

\begin{lemma}\label{OG}
A morphism $f\co G_{\bullet} \to H_{\bullet}$ in $s\agp$ is a weak equivalence if and only if the structural morphism $f^{\sharp}\co O(H)^{\bt} \to O(G)^{\bt}$ is a quasi-isomorphism of cosimplicial complexes, i.e. induces an isomorphism on cohomology. 
\end{lemma}
\begin{proof}
Since $\pi_0(G) =\Spec \H^0(O(G)^{\bt})$, it is equivalent for $\pi_0(f)$ and $\H^0(f^{\sharp})$ to be isomorphisms. We may therefore choose the same Levi decomposition for both $G_0$ and $H_0$, thus getting a Levi decomposition which is preserved by $f$. Since any pro-unipotent group is isomorphic to its Lie algebra, via exponentiation, the Levi decomposition gives isomorphisms
$$
O(G)^{\bt}\cong O(G^{\red})[ \fr_u(G_{\bt})^{\vee}],
$$
and similarly for $H$. Here,  $\fr_u(G_n)$ is the (pro-finite-dimensional) Lie algebra of $\Ru(G)$, and  $\fr_u(G_n)^{\vee}$ is its dual, which is an ind-finite-dimensional vector space (equivalently, just a vector space).

Exponentiation gives an isomorphism between $\pi_n(G_{\bt})$ and $\pi_n(\fr_u(G_{\bt}))$ for $n>0$. We therefore have isomorphisms on all $\pi_n(f)$ if and only if we have isomorphisms on the cohomology of the polynomial rings above, i.e. on all $\H^n(f^{\sharp})$.
\end{proof}
\begin{remark}
This lemma is one of the main reasons for our definition of pro-algebraic simplicial groups, since the corresponding statement is not true for arbitrary simplicial objects in the category of pro-algebraic groups.
\end{remark}

\begin{definition}
Recall from \cite{higgs} Definition 2.16 that $\cL \in \hat{\cN}(H)$ is said to be a hull for a functor $F:\cN(H) \to \Set$ if there is a natural transformation $\Hom(\cL, -) \to F$ such that
for every surjection $M \twoheadrightarrow N$ in $\cN(H)$, the canonical map  $\Hom(\cL,M) \to F(M)\by_{F(N)}\Hom(\cL,N)$ is surjective, and 
$\Hom(\cL, M) \to F(M)$ is an isomorphism for all $M$ abelian.
\end{definition}  

\begin{definition}
Given a morphism $G \xra{f} H$ of pro-algebraic groups, with $H$ reductive, define $G^{\rtimes,f}$ to be the object of $\hat{\cN}(H)$ which is the hull of  the functor 
$$
U \mapsto \Hom(G, U \rtimes H)_f/U,
$$
with $U$ acting by conjugation. 

 Note that from the Levi decomposition, it follows that for $H=G^{\red}$, and $f$ the canonical projection, the map $\Ru(G) \to G^{\rtimes,f}$ is an isomorphism. 
\end{definition}

\begin{definition}
Given a morphism $G \xra{f} H$ of pro-algebraic groups, with $G$ reductive, let $\Ru(H,f)$ be the maximal pro-unipotent subgroup of $H$ which is normalised by $G$. Note that in general $\Ru(H) \le \Ru(H,f)$, and that  if $G=H^{\red}$ and $f$ is a choice of Levi decomposition, then $\Ru(H,f)=\Ru(H)$.
\end{definition}

\begin{lemma}\label{cf} A morphism $G \xra{f} H$ of pro-algebraic simplicial groups is
\begin{enumerate}
\item a fibration if and only if $\Ru(G) \to \Ru(H,f\circ i)$ is a fibration in $s\hat{\cN}(G^{\red})$, for some Levi decomposition $i\co G^{\red} \to G$;
\item a cofibration if and only if $G^{\rtimes,r\circ f} \to \Ru(H)$ is a cofibration in $s\hat{\cN}(H^{\red})$, for $r\co H \to H^{\red}$ the canonical projection, and $\Ru(H)\in s\hat{\cN}(H^{\red})$ via some choice of Levi decomposition.
\end{enumerate}
\end{lemma}
\begin{proof}
\begin{enumerate}
\item For $n>0$, $N_n(G)=N_n(\Ru(G))$, and similarly for $H$. Now,  $f(\Ru(G))$ is pro-unipotent, so must be contained in $\Ru(H,f)$, giving the morphism required. To say that this morphism in $s\hat{\cN}(G^{\red})$ is a fibration is equivalent to saying that for all $n>0$, the $N_n (\Ru(G)) \to N_n(\Ru(H,f))$ are surjective, but this is precisely the condition for $f$ to be a fibration in $s\agp$.

\item Assume that we have a commutative diagram
$$
\begin{CD}
G @>n>> G'\\
@VfVV @VVf'V\\
H @>m>> H',
\end{CD}
$$
with $f'$ a trivial fibration, and $G^{\rtimes,f} \to \Ru(H)$  a cofibration in $s\hat{\cN}(H^{\red})$. For $f'$ to be a trivial fibration is equivalent to saying that $f'$ is surjective, and that $\ker(f')$ is pro-unipotent and acyclic. Therefore $(f')^{-1} \im(m) \to \im(m)$ is a trivial fibration, so we may assume that $m$ is surjective.

Now, ${H'}^{\red} \cong {G'}^{\red}$, and after fixing a choice of Levi decomposition for $H$, $m(H^{\red})$ is a reductive subgroup of $H'$, so is contained in some maximal reductive subgroup (corresponding to a Levi decomposition for $H'$). We may choose the Levi decomposition for $G'$ compatibly,  since every trivial fibration is a pro-unipotent extension. The diagram is now
$$
\begin{CD}
G @>n>> (G')^{\red} \ltimes \Ru(G')\\
@VfVV @VVf'V\\
H^{\red} \ltimes \Ru(H) @>m>> (G')^{\red} \ltimes \Ru(H'),
\end{CD}
$$
with $m$ and $f'$ preserving the Levi decompositions. Since $H^{\red} \onto (G')^{\red}$, this is equivalent to giving the commutative diagram
$$
\begin{CD}
G^{\rtimes,r\circ f} @>n>>  \Ru(G')\\
@VfVV @VVf'V\\
 \Ru(H) @>m>>  \Ru(H')
\end{CD}
$$
in $s\hat{\cN}(H^{\red})$. Observe that $\Ru(G') \xra{f'} \Ru(H')$ is a trivial fibration, so the hypothesis that $G^{\rtimes,r\circ f} \to \Ru(H)$ be a cofibration ensures existence of the lift.

Conversely, assume that $G \xra{f} H$ is a cofibration. Then for every commutative diagram
$$
\begin{CD}
G^{\rtimes,r\circ f} @>n>>  U\\
@VfVV @VVf'V\\
 \Ru(H) @>m>>  V
\end{CD}
$$
in $s\hat{\cN}(H^{\red})$,  with $f'$ a trivial fibration, we    have a commutative diagram in $s\agp$
$$
\begin{CD}
G @>n>> H^{\red} \ltimes U\\
@VfVV @VVf'V\\
H^{\red} \ltimes \Ru(H) @>m>> H^{\red} \ltimes V,
\end{CD}
$$
with $f'$ a trivial fibration. We therefore have a lift $H^{\red} \ltimes \Ru(H) \xra{\tilde{m}} H^{\red} \ltimes U$, such that $f'\tilde{m}=m$ and $\tilde{m}f=n$. Thus $ \tilde{m}(H^{\red})\ltimes U$ is an alternative Levi decomposition for $H^{\red} \ltimes U$, and $\tilde{m}$ therefore provides a lift from $\Ru(H)$ to $U$ in $\cN(H^{\red})$, as required.
\end{enumerate}
\end{proof}

We will now follow the conventions and notation of \cite{sht} for closed model categories.
\begin{theorem}\label{cmsagp}
With the classes of morphisms given above, $s\agp$ is a closed model category.
\end{theorem}
\begin{proof}
\begin{enumerate}
\item[CM1] The category of algebraic groups is clearly closed under all  limits and colimits. Given a colimit diagram in $s\agp$, let $G_{\bt}$ be the colimit as a simplicial object over algebraic groups. Taking the pro-unipotent completion with respect to the projection onto $\pi_0(G)$ gives the colimit as a pro-algebraic simplicial group. 

A product of pro-algebraic simplicial groups is clearly a pro-algebraic simplicial group. It only remains to verify that an equaliser of pro-algebraic simplicial groups is  a pro-algebraic simplicial group. Let
$$
\xymatrix@1{K \ar[r] &G \ar@<1ex>[r] \ar@<-1ex>[r]  &   H }
$$
be an equaliser diagram in the category of  simplicial objects over algebraic groups, with $G$ and $H$ pro-algebraic simplicial groups. Then $K_n$ is a subgroup of $G_n$, and the kernel of $K_n \to \pi_0(G)$ is pro-unipotent, so the kernel of $K_n \to \pi_0(K)$ must also be, and thus $K$ is a pro-algebraic simplicial group.

\item[CM2] If $f$ and $g$ are composable morphisms, with any two of $f,g, f\circ g$ weak equivalences, then so is the third, by the long exact sequence of homotopy for simplicial groups.

\item[CM3] Retracts of fibrations, cofibrations and weak equivalences are (respectively) fibrations, cofibrations and weak equivalences, by the corresponding properties for $s\hat{\cN}(G)$ and for simplicial groups. 

\item[CM4] By definition, every cofibration has the left lifting property with respect to all trivial fibrations. We need to show that every trivial cofibration has the left lifting property with respect to all  fibrations. We use Lemma \ref{cf}, so that this amounts to having a commutative diagram
$$
\begin{CD}
K \ltimes U @>n>> G \\
@VfVV @VVf'V\\
K \ltimes V @>m>> H,
\end{CD}
$$
with $K$ reductive, $f\co U \to V$ a trivial cofibration in $s\hat{\cN}(K)$, and $f'$ a fibration. This is equivalent to the diagram
$$
\begin{CD}
 U @>n>> \Ru(G, n\circ i) \\
@VfVV @VVf'V\\
 V @>m>> \Ru(H, m\circ i)
\end{CD}
$$
in $s\hat{\cN}(K)$, noting that $f'$ is also a fibration here. Therefore the lift exists, since $s\hat{\cN}(K)$ is a model category. 

\item[CM5] Given an arbitrary morphism $f\co G \to H$ in $s\agp$, choose a Levi decomposition $G^{\red} \to G$, and let $r\co H \to H^{\red}$ be the canonical quotient. These give morphisms $\Ru(G) \to \Ru(H,f\circ i)$, and $G^{\rtimes, r\circ f} \to \Ru(H)$ in $\cN(G^{\red})$ and $\cN(H^{\red})$, respectively. Since these are both model categories, we may choose factorisations
\begin{eqnarray*}
\Ru(G) \xra{a} &U& \xra{b} \Ru(H,f\circ i)\\
G^{\rtimes, r\circ f} \xra{c} &V& \xra{d} \Ru(H),
\end{eqnarray*}
with $a,c$ cofibrations, $b,d$ fibrations, and $a,d$ weak equivalences. Finally, observe that $G^{\red} \rtimes \Ru(H,f\circ i) \to H$ is a fibration, and that $G \to H^{\red}\rtimes G^{\rtimes, r\circ f}$ is a cofibration, so that the factorisations
\begin{eqnarray*}
G \xra{a} &G^{\red} \rtimes U& \to H\\
G \to &H^{\red} \rtimes V& \xra{d} H,
\end{eqnarray*}
are of the form cofibration-fibration, with $a,d$ trivial, as required.
\end{enumerate}
\end{proof}

In fact, $s\agp$ has a simplicial model structure, which we now describe.
\begin{definition}
Given $G \in s\agp$ and $K \in \bS$, consider the simplicial complex $H_{\bt}$ of pro-algebraic groups which in level $n$ is
$$
H_n:= \overbrace{G_n*G_n*\ldots*G_n}^{K_n},
$$
where $*$ denotes coproduct in the category of pro-algebraic groups. Define $G\ten K \in s\agp$ to be the pro-unipotent completion of $H$ with respect to the morphism  $H_{\bt}\to \pi_0(H)$. 
\end{definition}

\begin{proposition}\label{sagpsimplicial}
The category $s\agp$ is a simplicial category under the above construction.  
\end{proposition}
\begin{proof}
By \cite{sht} Lemma II.2.4, we must verify the following conditions:
\begin{enumerate}
\item For fixed $K \in \bS$, $\ten K\co  s\agp \to s\agp$ has a left adjoint $G \mapsto G^K$. We can describe this left adjoint explicitly by $(G^K)_n(A)=\Hom_{\bS}(\Delta \by K, G(A))$, for all $k$-algebras $A$, so that $G^K(A)=(G(A))^K$ as a simplicial set. 

\item For fixed $G \in s\agp$, the functor $G\ten \co \bS \to s\agp$ commutes with arbitrary colimits, and $G\ten \bt \cong G$. This is immediate.

\item There is an isomorphism $G\ten (K\by L) \cong (G\ten K)\ten L$ natural in $G \in s\agp$, and $K,L \in \bS$. This is also immediate.
\end{enumerate}
\end{proof}

\begin{theorem}
With the structures described above, $s\agp$ is a simplicial model category.
\end{theorem}
\begin{proof}
It only remains to establish the simplicial model axiom SM7. By the dual version of \cite{sht} Corollary II.3.12, this is equivalent to showing that for all fibrations $p\co G \to H$ in $s\agp$, the map
$$
H^{\Delta^n}\to G^{\Delta^n} \by_{G^{\pd\Delta^n}}H^{\pd\Delta^n} 
$$
is a fibration for $n\ge 0$, which is trivial if $p$ is, and that
$$
  H^{\Delta^1}\to G^{\Delta^1} \by_{G^{\{e\}}}H^{\{e\}}
$$
is the trivial fibration for $e=0$ or $1$. This follows from the corresponding results for the simplicial sets $G(A),H(A)$.
\end{proof}

\subsection{Algebraisation for connected spaces}

There is a forgetful functor $(k)\co s\agp \to s\gp$, given by sending $G_{\bt}$ to $G_{\bt}(k)$. This functor clearly commutes with all limits, so has a left adjoint $G_{\bt} \mapsto (G_{\bt})^{\alg}$. We can describe $(G_{\bt})^{\alg}$ explicitly. First let $(\pi_0(G))^{\alg}$ be the pro-algebraic completion of the abstract group $\pi_0(G)$, then let $(G^{\alg})_n$ be the relative Malcev completion  (in the sense of \cite{malcev}) of the morphism
$$
G_n \to (\pi_0(G))^{\alg}.
$$
In other words, $G_n \to (G^{\alg})_n \xra{f} (\pi_0(G))^{\alg}$ is the universal diagram with $f$ a pro-unipotent extension.

\begin{proposition}\label{algq}
The functors $(k)$ and ${}^{\alg}$ are a pair of Quillen functors.
\end{proposition} 
\begin{proof}
It suffices to observe that $(k)$ preserves fibrations and trivial fibrations. 
\end{proof}

\begin{definition}
Given a reduced simplicial set (or equivalently a pointed, connected topological space), define the pro-algebraic homotopy type $X^{\alg}$ of $X$ over $k$ to be the equivalence class of 
$$
G(X)^{\alg}
$$ 
in the closed model category of pro-algebraic simplicial groups. 
\end{definition}

\begin{remark}
In \cite{schematic}, the homotopy type associated to $X$ is denoted $GX^{\alg} \in s\mathrm{GAff}$, given by $( GX^{\alg})_n=(GX_n)^{\alg}$. Thus our homotopy type $G(X)^{\alg}$ is the relative unipotent 
completion of $GX^{\alg}$ over $\pi_0G$.
\end{remark}

It follows from Proposition \ref{algq} that the algebraisation functor from simplicial groupoids to pro-algebraic simplicial groupoids induces a functor on the homotopy categories.  
\begin{definition}
Given $X \in \bS_0$,  we  define the pro-algebraic fundamental group $\varpi_1(X):=\pi_0(G(X)^{\alg})$, where $G\co \bS_0\to s\gp$ is Kan's loop group functor. We then define the higher homotopy groups $\varpi_n(X)$  by 
$$
\varpi_n(X):=\pi_{n-1}(G(X)^{\alg}).
$$
Note that these will depend only upon the homotopy type of $X$. Given a  pointed, connected topological space $(X,x)$, define the pro-algebraic homotopy groups of $X$ by $\varpi_n(X,x):=\varpi_n(\Sing(X,x))$, i.e. the pro-algebraic homotopy groups of the simplicial set 
$$
\Sing(X,x)_n=\{ f \in \Hom_{\Top}(|\Delta^n|,X)\,|\, f(v)=x\quad \forall  v \in (\Delta^n)_0\}
$$ 
of $(X,x)$.
\end{definition} 

\begin{theorem}
Given $X \in \bS_0$, $\varpi_1(X)=\pi_1(X)^{\alg}$.
\end{theorem}
\begin{proof}
$$
\varpi_1(X)=\pi_0(G(X)^{\alg})=\pi_0(G(X))^{\alg},
$$
as required.
\end{proof}

\subsection{Homology of pro-algebraic simplicial groups}

\begin{definition}
Given a pro-algebraic simplicial group $G$ over $k$, we define $s\widehat{\FD\Rep}(G)$ to be the category of  pro-finite-dimensional simplicial $G$-representations. Explicitly, an object is a simplicial pro-finite-dimensional $k$-vector space $V$, together with  morphisms  $G\by V \xra{\rho} V$ of simplicial schemes satisfying the associativity and identity axioms. In other words, $\theta(\rho(g)\cdot v) = \rho(\theta(g))\cdot \theta(v)$, for $g \in G_n, v \in V_n$, $\theta=\pd_i,\sigma_i$.

There is a forgetful functor $s\widehat{\FD\Rep}(G)\to s\widehat{\FD\Vect}$, with left adjoint $V \mapsto V\hat{\ten}_k O(G)^{\vee}$. Note that $O(G)$ is an ind-finite-dimensional $G$-representation, so its dual $O(G)^{\vee} \in  s\widehat{\FD\Rep}(G)$. 
\end{definition}

\begin{definition}
Define a morphism $f\co V \to W$ in $s\widehat{\FD\Rep}(G)$ to be: 
\begin{enumerate}
\item a weak equivalence if the morphism of underlying pro-vector spaces in $s\widehat{\FD\Vect}$ is a weak equivalence;

\item a fibration if  the morphism of underlying pro-vector spaces in $s\widehat{\FD\Vect}$ is a fibration;

\item a cofibration if  it has the left lifting property (LLP) with respect to all trivial fibrations.
\end{enumerate}
\end{definition}

\begin{proposition}
With the classes of morphisms given above, $s\widehat{\FD\Rep}(G)$ is a closed model category.
\end{proposition}
\begin{proof}
This follows from \cite{ship}, since $s\widehat{\FD\Rep}(G)$ is the category of left modules for the monoid $O(G)^{\vee} $ in $s\widehat{\FD\Vect}$.
\end{proof}

\begin{lemma}
If $f\co U \to V$ is a morphism in $s\widehat{\FD\Rep}(G)$, such that $f_n\co U_n \to V_n$ is injective for all $n$, and there exists a set $\Sigma =\bigcup_n \Sigma_n$, with $\Sigma_n \subset V_n$, closed under degeneracy operators $\sigma_i$, such that
$$
V_n=f(U_n)\by \prod_{s \in \Sigma_n} O(G_n)^{\vee}s,
$$ 
then $f$ is a cofibration. Such a morphism is called free on generators $\Sigma$.
\end{lemma}
\begin{proof}
This is identical to the corresponding result for simplicial groups in \cite{QHA}.
\end{proof}

\begin{example}\label{owres}
The pro-finite-dimensional simplicial  $G$-representation $O(WG)^{\vee}$, given by $(O(WG)^{\vee})_n=(O(WG)^n)^{\vee}$, where
$$
O(WG)^n=O(G_n)\ten O(G_{n-1})\ten \ldots \ten O(G_0),
$$
with operations and $G$-action as  in \cite{sht} \S V.4, i.e.
\begin{eqnarray*}
\pd_i(v_n\ten v_{n-1}\ten \ldots \ten v_0)&=& \left\{ \begin{matrix} \pd_iv_n\ten\pd_{i-1}v_{n-1}\ten\ldots\ten (\pd_0v_{n-i})v_{n-i-1}\ten v_{n-i-2}\ten \ldots\ten v_0 & i<n,\\ \pd_nv_n\ten\pd_{n-1}v_{n-1}\ten \ldots\ten \pd_1v_1 & i=n, \end{matrix} \right.\\
\sigma_i(v_n\ten v_{n-1}\ten \ldots \ten v_0)&=& \sigma_iv_n\ten\sigma_{i-1}v_{n-1}\ten\ldots\ten \sigma_0v_{n-i}\ten e \ten  v_{n-i-1}\ten \ldots\ten v_0,
\end{eqnarray*}
and for $h \in G_n$, 
$$
h(v_n\ten v_{n-1}\ten \ldots \ten v_0)=(hv_n)\ten v_{n-1}\ten \ldots \ten v_0.
$$
This is cofibrant, and the morphism $ O(WG)^{\vee}\to k$ is a weak equivalence. Observe that for any pro-finite-dimensional simplicial $G$-representation $V$, $O(WG)^{\vee} \hat{\ten} V \to V$ will then be a cofibrant approximation. 
\end{example}

\begin{remark}
If $G$ is a simplicial group, then the finite-dimensional $G^{\alg}$-representations are precisely those finite-dimensional $G$-representations $V$ which admit filtrations $F^iV$ such that  the graded pieces $\gr^i_F V$ are $\pi_0(G)$-representations  for all $i$. 
\end{remark}

\begin{definition}
Given a pro-algebraic simplicial group $G$, define  $\bL(G,-)\co s\widehat{\FD\Rep}(G) \to s\widehat{\FD\Vect}$  to be the left-derived functors of the co-invariant functor  $V \mapsto V_G:=V/G$. 
\end{definition}

\begin{definition}
Given a pro-algebraic simplicial group $G$, and a pro-finite-dimensional simplicial $G$-representation $V$, define the homology groups of $G$ with coefficients in $V$ by
$$
\H_i(G,V):=\H_i(\bL(G,V)).
$$
\end{definition}

\begin{lemma}\label{owgpcoho}
For any pro-algebraic simplicial group $G$, $\H_{n+1}(G,k)$ is dual to $\H^{n+1}(O(\bar{W}G))$, where $O(\bar{W}G)$ is the invariant subspace  of $O(WG)$ (defined in Example \ref{owres}) under its canonical $G$-action.
\end{lemma}
\begin{proof}
Apply $\H_i\bL(G,-)$ to the canonical resolution $O(WG)^{\vee} \to k$. 
\end{proof}

\begin{proposition}\label{freecoho}
If $G$ is a cofibrant simplicial pro-unipotent pro-algebraic group, then
$$
\H_i(G,k)\cong\left\{ \begin{matrix} k & i=0\\ \H_{i-1}(G/[G,G]) &i>0  \end{matrix}    \right.
$$
\end{proposition}
\begin{proof}
Since $G$ is cofibrant and pro-unipotent, the Lie algebra $\g$ of $G$ is free. Let $V_*$ be a space of  free generators for $\g$, closed  under the degeneracy operators. Now, $O(G)^{\vee} $ is the universal enveloping algebra of $\g$.   Hence
$O(G)^{\vee}_n$ is the free associative (non-commutative) pro-nilpotent algebra generated by $V_n$, i.e.
$$
O(G)^{\vee}_n \cong \prod_{j\ge 0} V_n^{\hat{\ten}j}.
$$

If we now let $I$ be  augmentation ideal of the map $O(G)^{\vee} \to k$, observe that $I$ is cofibrant as a $G$-representation, since it is freely generated by $V_*$ as an $O(G)^{\vee}$-module. The mapping cone $C$ of the inclusion $I \to O(G)^{\vee}$ is then also cofibrant, and a resolution of $k$. Now $C/G$ is the mapping cone of $I/G \to O(G)^{\vee}/G$, which is just $I/I^2 \xra{0} k$. Since $\H_i(G,k)= \H_i(C/G)$ and $I/I^2 \cong G/[G,G]$, the result follows.
\end{proof}

\subsubsection{The Hochschild-Serre spectral sequence}
Let $G$ be a pro-algebraic simplicial group, and $K \lhd G$ a  pro-algebraic simplicial normal subgroup. We wish to compare the homology groups of $G, K$ and $G/K$.
\begin{theorem}\label{hs}
For $G,K$ as above, and  every pro-finite-dimensional simplicial $G$-representation $V$, there is a convergent first quadrant spectral sequence
$$
E^2_{pq}=\H_p(G/K,\H_q(K,V))\abuts \H_{p+q}(G,V)
$$
\end{theorem}
\begin{proof}
We have a commutative triangle
$$
\xymatrix{
  s\widehat{\FD\Rep}(G) \ar[rr]^{-_K} \ar[dr]_{-_G} & &  s\widehat{\FD\Rep}(G/K) \ar[dl]_{-_{G/K}}\\
           &              s\widehat{\FD\Vect}   
}
$$
of left Quillen functors of closed model categories. The Grothendieck spectral sequence thus exists, as required.
\end{proof}

\subsubsection{Reverse Adams spectral sequence}

\begin{definition}
Given a cosimplicial vector space $V$, define the symmetric power 
$$
\Symm^nV=V^{\ten n}/S_n, 
$$
where $V^{\ten n}$ is the cosimplicial complex $(V^{\ten n})^i=(V^i)^{\ten n}$.

Given a cochain complex $V$, define the symmetric power 
$$
\Symm^nV=V^{\ten n}/S_n, 
$$
where $V^{\ten n}$ is the cochain complex 
$$
(V^{\ten n})_i= \bigoplus_{i_1+i_2+\ldots + i_n=i} (V^{i_1})\ten (V^{i_2})\ten \ldots \ten (V^{i_n}),
$$
on which the symmetric group $S_n$ acts using the usual graded-commutative convention $ab=(-1)^{\deg a \deg b}ba$.
\end{definition}

\begin{theorem}\label{spectralpi}
If $G$ is  a pro-unipotent pro-algebraic simplicial group, with Lie algebra $\g$, there is a canonical convergent reverse Adams spectral sequence 
$$
E_1^{pq}=(\mathrm{Symm}^p(\pi_{*-1}(\g)^{\vee}))^{p+q} \abuts \H^{p+q}(O(\bar{W}G)),
$$
and $\exp(\pi_0\g) \cong \pi_0G$, $\pi_i\g \cong \pi_iG$ for all $i>0$.
\end{theorem}
\begin{proof}
The exponential map corresponds to an isomorphism
$$
O(\exp \g) \cong k[\g^{\vee}]
$$ 
between the structure sheaf of $G$ and the polynomial ring on the dual of $\g$, as in Lemma \ref{OG}. 

The identity map $1 \to G$ gives a cosimplicial ring homomorphism $O(\bar{W}G) \to k$. Let $I$ be the augmentation ideal of this homomorphism. Then 
$$
I/I^2 \cong (\bar{W}(\g,+))^{\vee},
$$
where $(\g,+)$ is  $\g$ with its additive  group structure (forgetting the bracket). 

The powers $I^n$ of $I$ then give a filtration of $O(\bar{W}G)$, with $I^n/I^{n+1} \cong \Symm^n((\bar{W}(\g,+))^{\vee})$. This gives us a  bounded spectral sequence
$$
E_1^{pq}=\H^{p+q}(\Symm^p(\bar{W}(\g,+))^{\vee})) \abuts \H^{p+q}(O(\bar{W}G)),
$$
which is thus convergent.

For any cosimplicial complex $V$, 
$$  
\H^*(\Symm^pV)\cong \Symm^p(\H^*V),
$$
since $k$ is of characteristic $0$. 
Thus
$$
\Symm^p(\H^*((\bar{W}(\g,+))^{\vee}))^{p+q} \abuts \H^{p+q}(O(\bar{W}G)).
$$

Now, $\H^i((\bar{W}(\g,+))^{\vee})$ is dual to $\pi_i(\bar{W}(\g,+))$, which is isomorphic to $\pi_{i-1}(\g,+)$, since $(\g,+)$ is abelian. Since $\g$ and $G$ are isomorphic  as fibrant simplicial sets via the exponential map, $\pi_{i-1}(\g,+)\cong \pi_{i-1}(G)$ are isomorphic as groups for $i>1$ and as sets for $i=1$.
\end{proof}

\begin{corollary}
A morphism $G \xra{f} K$ of pro-unipotent pro-algebraic simplicial groups is a weak equivalence if and only if 
$$
\H_i(G,k) \to \H_i(K,k)
$$
is an isomorphism for all $i$.
\end{corollary}
\begin{proof}
That isomorphisms on homotopy give isomorphisms on homology is an immediate consequence of Theorem \ref{spectralpi}. The converse follows by contradiction, considering the first $i$ for which $\pi_i(f)$ is not an isomorphism.
\end{proof}

\begin{corollary}\label{detectweak}
A morphism $G \xra{f} K$ of  pro-algebraic simplicial groups is a weak equivalence if and only if $f(\Ru(G))\le \Ru(K)$,  with the quotient map
$$
G^{\red} \to K^{\red}
$$
 an isomorphism, and for all (finite-dimensional) irreducible $K$-representations $V$, the maps
$$
\H_i(f)\co  \H_i(G,f^*V) \to \H_i(K,V)
$$
are isomorphisms for all $i>0$.
\end{corollary}
\begin{proof}
Applying Theorem \ref{hs} to the extension $G \to G^{\red}$ gives the isomorphism
$$
\H_i(G,V) \cong \H_i(\Ru(G),k)\ten_{G^{\red}}V.
$$
We now just apply the previous corollary to the map $\Ru(f)\co \Ru(G) \to \Ru(K)$.
\end{proof}

\subsubsection{Classical homotopy groups}\label{clpi}

\begin{definition}
Recall (e.g. from \cite{sht} Definition VI.3.4) that the Moore-Postnikov tower $\{X(n)\}$ of a fibrant simplicial set $X$ is given by
$$
X(n)_q:=\im( X_q \to \Hom(\sk_n\Delta^q, X)),
$$  
with the obvious simplicial structure. These form an inverse system $X\to \ldots \to X(n) \to X(n-1) \to \ldots$, with $X=\Lim X(n)$, and 
$$
\pi_q X(n)= \left\{ \begin{matrix} \pi_q X & q \le n \\ 0 & q >n. \end{matrix} \right.
$$ 
The maps $X(n)\to X(n-1)$ are fibrations. If $X$ is reduced, then so is $X(n)$, and we define $E(n)$ to be the fibre of $X \to X(n)$ over the basepoint $*$.
\end{definition}

\begin{definition}
Recall from \cite{schematic}v1 Lemma 4.11 that a discrete group $\Gamma$ is said to be algebraically good if for all finite-dimensional $\Gamma$-representations $V$, the map $\H^*(\Gamma^{\alg},V) \to \H^*(\Gamma,V)$ is an isomorphism.
\end{definition}

\begin{theorem}\label{classicalpi}
If $X$ is a reduced topological space with fundamental group $\Gamma$, such that 
\begin{enumerate}
\item $\Gamma$ is  algebraically good,  and 
\item $\pi_n(X)$ is of finite rank for all $n>1$,
\end{enumerate}
then the canonical map
$$
  \pi_n(X)\ten_{\Z} k \to \varpi_{n}(X) 
$$
is an isomorphism for all $n>1$.
\end{theorem}
\begin{proof}
This can be thought of as a generalisation of \cite{KTP} Lemma 4.3.2 to arbitrary topological spaces. Observe that $G(X)_q= G(X(n))_q$ for all $q <n$, so $\varpi_q(X)=\Lim \varpi_q(X(n))$. It therefore suffices to prove the theorem for all $X(n)$, and we proceed by induction on $n$. 

$G(X(1))$ is a cofibrant resolution of $\Gamma$, so for any finite-dimensional $\Gamma^{\alg}$-representation $V$, $\H^*(X(1), V)\cong \H^*(\Gamma, V)$.
Since $\Gamma$ is algebraically good, the map
$$
G(X(1))^{\alg} \to \Gamma^{\alg}
$$
gives isomorphisms 
$$
\H^*(\Gamma^{\alg}, V) \to \H^*(\Gamma, V) 
$$
on cohomology for all such $V$, so is a weak equivalence by Corollary \ref{detectweak}.

Now, assume that  $G(X(n-1))^{\alg}$ satisfies the inductive hypothesis, and consider the fibration $X(n) \to X(n-1)$. This is determined up to homotopy  by a k-invariant (\cite{sht} \S VI.5) $\kappa \in \H^{n+1}(X(n-1), \pi_{n}(X))$. Since $A:=\pi_{n}(X))\ten_{\Z}k$ is a finite-dimensional $\Gamma^{\alg}$-representation by hypothesis, the element 
$$
\kappa \in \H^{n+1}(X(n-1), A) \cong \H^{n+1}(X(n-1)^{\alg}, A)
$$
comes from a map
$$
G(X(n-1))^{\alg}\to (N^{-1}A[-n])\rtimes R.
$$

Let $L$ be the chain complex  $A[-n] \xra{\id} A[1-n]$, and define $\cG$ to be the pullback of this map along the surjection
$ N^{-1}L \rtimes R \to(N^{-1}A[-n])\rtimes R$. This gives an
extension 
$$
N^{-1}A[1-n]  \to \cG \to  G(X(n-1))^{\alg}.
$$

Taking $k$-valued points gives the fibration
$$
\bar{W}N^{-1}A[1-n] \to \bar{W}\cG(k) \to \bar{W}G(X(n-1))^{\alg}(k)
$$
in $\bS$, corresponding to the k-invariant $f^*\kappa \in \H^n(\bar{W}G(X(n-1))^{\alg}(k),A)$, for $f:X(n-1) \to \bar{W}G(X(n-1))^{\alg}(k)$. This in turn gives a map $X(n) \to \bar{W}\cG(k)$, compatible with the fibrations. 

From the long exact sequence of homotopy, it follows that $\cG$ has the required homotopy groups, so it will suffice to show that $F:G(X(n))^{\alg} \to \cG$ is a weak equivalence. We now apply the Hochschild-Serre spectral sequence (Theorem \ref{hs}), giving
$$
\H_p(X(n-1), \H_q( N^{-1}A[1-n],V))= \H_p(G(X(n-1))^{\alg},\H_q(N^{-1}A[1-n],V))\abuts \H_{p+q}(\cG,V).
$$
Similarly
$$
\H_p(X(n-1),\H_q(E(n),V))\abuts \H_{p+q}(X(n),V),
$$
for all finite-dimensional $\Gamma^{\alg}$-representations $V$.

Now, since $E(n)$ is simply connected, it follows from \cite{QRat} (or just the Curtis Convergence Theorem) that $G(E(n))^{\alg}\to N^{-1}A[1-n]$ is a weak equivalence. Corollary \ref{detectweak} then implies that $\H_q(N^{-1}A[1-n],V) \cong \H_q(E(n),V)$, so $F$ induces isomorphisms on homology groups, hence must be a weak equivalence, as required.
\end{proof}

\section{Unpointed pro-algebraic homotopy types}\label{unptd}

\begin{definition}
Let $\bS$ be the category of simplicial sets, and $s\gpd$ the category of simplicial groupoids on a constant set of objects (as in \cite{sht}). Let $\Top$ denote the category of compactly generated Hausdorff topological spaces.
\end{definition}

Note that there is a functor from $\Top$ to $\bS$ which sends $X$ to the simplicial set
$$
\Sing(X)_n= \Hom_{\Top}(|\Delta^n|, X).
$$
this is a right Quillen equivalence, the corresponding left equivalence being geometric realisation. From now on, we will thus restrict our attention to simplicial sets.

Whereas in Section \ref{ptd} we studied the loop group of a reduced simplicial set, here we will study the path groupoid of a simplicial  set.
As in \cite{sht} \S V.7, there is a  functor $\bar{W}\co s\gpd \to \bS$, with left adjoint $G\co \bS \to s\gpd$, Dwyer and Kan's  path groupoid functor (\cite{pathgpd}), and these give a Quillen equivalence of model categories. This section is devoted to defining and studying pro-algebraic completions of simplicial groupoids.

\subsection{Pro-algebraic groupoids}\label{gpdsn}
\begin{definition}
Define a pro-algebraic groupoid $G$ over $k$ to consist of the following data:
\begin{enumerate}
\item A discrete set $\Ob(G)$.
\item For all $x,y \in \Ob(G)$, an affine scheme $G(x,y)$ (possibly empty) over $k$.
\item A groupoid structure on $G$, consisting of an associative multiplication  morphism $m\co G(x,y)\by G(y,z) \to G(x,z)$, identities $\Spec k \to G(x,x)$ and inverses $G(x,y) \to G(y,x)$
 \end{enumerate}
We say that a pro-algebraic groupoid is reductive (resp. pro-unipotent) if the pro-algebraic groups $G(x,x)$ are so for all $x \in \Ob(G)$. An algebraic groupoid is a pro-algebraic groupoid for which the $G(x,y)$ are all of finite type.
\end{definition}
If $G$ is a pro-algebraic groupoid, let $O(G(x,y))$ denote the global sections of the structure sheaf of $G(x,y)$.

\begin{definition}
Given morphisms $f,g\co G \to H$ of pro-algebraic groupoids, define a natural isomorphism $\eta$ between $f$ and $g$ to consist of morphisms
$$
\eta_x\co  \Spec k \to H(f(x),g(x))
$$
for all $ x\in \Ob(G)$, such that the following diagram commutes, for all $x,y \in \Ob(G)$:
$$
\begin{CD}
G(x,y) @>f(x,y)>> H(f(x),f(y))\\
@Vg(x,y)VV  @VV{\cdot\eta_y}V \\
 H(g(x),g(y)) @>{\eta_x\cdot}>>  H(f(x),g(y)).
\end{CD}
$$
 
A morphism $f\co G \to H$ of pro-algebraic groupoids is said to be an equivalence if there exists a morphism $g\co H \to G$ such that $fg$ and $gf$ are both naturally isomorphic to identity morphisms. This is the same as saying that for all $y \in \Ob(H)$, there exists $x \in \Ob(G)$ such that $H(f(x),y)(k)$ is non-empty (essential surjectivity), and that for all $x_1,x_2 \in \Ob(G)$, $G(x,y) \to G(f(x_1),f(x_2) )$ is an isomorphism.
\end{definition}

\begin{remark}\label{pathological}
Note that if $k$ is not algebraically closed, then $G(x,y)$ might be non-empty and $G(x,y)(k)$ empty. An example is to take $k=\R$, set $\Ob(G)=\{x,y\}$, with $G(x,x), G(y,y)$ both the constant groups $\{1,\tau\}$, and $G(x,y)=\Spec \Cx$, on which $\tau$ acts by complex conjugation. This is an example of an algebraic groupoid which is not equivalent to a disjoint union of algebraic groups. 
\end{remark}

\begin{definition}\label{gpdrep}
Given a pro-algebraic groupoid $G$, define a finite-dimensional linear $G$-representation to be a functor $\rho\co G \to \mathrm{FDVect}_k$ respecting the algebraic structure. Explicitly, this  consists of a set $\{V_x\}_{x \in \Ob(G)}$ of finite-dimensional $k$-vector spaces, together with morphisms $\rho_{xy}\co G(x,y) \to \Hom(V_y,V_x)$ of affine schemes, respecting the multiplication and identities. 

A morphism $f\co (V,\rho)\to (W,\varrho)$ of $G$-representations consists of $f_x \in \Hom(V_x,W_x)$ such that 
$$
f_x\circ\varrho_{xy}=\rho_{xy}\circ f_y\co G(x,y) \to \Hom(V_x,W_y).
$$

There is a similar definition for linear representations of an abstract groupoid $G$, where instead we have morphisms $\rho_{xy}\co G(x,y) \to \Hom(V_y,V_x)(k)$
 of sets. 
Thus the category of local systems on a topological space $X$ is equivalent to the category of representations of the fundamental groupoid $\pi_f X$, and the category of linear representations of an abstract groupoid $G$ is equivalent to the category of local systems on the classifying space $BG$.
\end{definition}

\begin{remark}
There is a form of Tannakian duality for pro-algebraic groupoids, extending Tannakian duality (\cite{tannaka} Theorem II 2.11) for pro-algebraic groups. 
For each $x \in \Ob(G)$, we  have a fibre functor
$$
\omega_x\co \FD\Rep(G) \to \mathrm{FD}\Vect(k),
$$
where $\FD\Rep(G)$ is the tensor category of finite-dimensional $G$-representations. We can recover $G$ from these data (a rigid tensor category together with a jointly faithful collection of fibre functors) by setting
$$
G(x,y):=\Iso^{\ten}(\omega_x,\omega_y).
$$
Equivalent multi-fibred categories then give isomorphic pro-algebraic groupoids. 

In the example of Remark \ref{pathological}, $\Rep(G)$ is equivalent to the category of $C_2$-representations, with the fibre functors $\omega_x,\omega_y$ having the same underlying functor, but for the non-trivial representation $V$, the structural isomorphisms
\begin{eqnarray*}
\omega_x(V^{\vee}) &\to& \omega_x(V)^{\vee}\\
\omega_y(V^{\vee}) &\to& \omega_y(V)^{\vee}
\end{eqnarray*}
differ by a factor of $-1$.
\end{remark}

\begin{definition}
Given a pro-algebraic groupoid $G$, define the reductive quotient $G^{\red}$ of $G$ by setting $\Ob(G^{\red})=\Ob(G)$, and
$$
G^{\red}(x,y)=G(x,y)/\Ru(G(y,y))= \Ru(G(x,x))\backslash G(x,y),
$$
the equality arising since if $f\in G(x,y),\, g \in \Ru(G(y,y))$, then $f gf^{-1}\in \Ru(G(x,x))$, so both equivalence relations are the same. Multiplication and inversion descend similarly. Observe that $G^{\red}$ is then a reductive pro-algebraic groupoid.
\end{definition}

\begin{definition}
Let $\agpd$ denote the category of pro-algebraic groupoids over $k$, and observe that this category is contains all (inverse) limits. There is functor from $\agpd$ to $\gpd$, the category of abstract groupoids, given by $G \mapsto G(k)$. This functor preserves all limits, so has a left adjoint, the algebraisation functor, denoted $\Gamma \mapsto \Gamma^{\alg}$. This can be given explicitly by $\Ob(\Gamma)^{\alg}=\Ob(\Gamma)$,  and 
$$
\Gamma^{\alg}(x,y)=\Gamma(x,x)^{\alg}\by^{\Gamma(x,x)}\Gamma(x,y),
$$
where $\Gamma(x,x)^{\alg}$ is the pro-algebraic completion of the group $\Gamma(x,x)$. 

The finite-dimensional linear representations of $\Gamma$ (as in Definition \ref{gpdrep}) correspond to those of $\Gamma^{\alg}$, and these can be used to recover $\Gamma^{\alg}$, by Tannakian duality.
\end{definition}

\subsection{Groupoid cohomology and obstructions}\label{gpdcoho}

As is proved in \cite{tannaka} Proposition II2.2 for pro-algebraic groups, a (not necessarily finite-dimensional) $G$-representation is an $O(G)$-comodule.
\begin{definition}
Given a pro-algebraic groupoid $G$, and a $G$-representation $V$, then $V^{\vee}$ is a pro-finite-dimensional $G$-representation with  $V_x^{\vee}=\Spec k[V_x]$,  and  we define
\begin{eqnarray*}
\CC^n(G,V^{\vee})&:=&\prod_{x_0,\ldots, x_n \in \Ob(G)} \Hom_{\Sch}(G(x_0,x_1)\by G(x_1,x_2) \by \ldots G(x_{n_1},x_n), V^{\vee}_{x_0})\\
&=& \prod_{x_0,\ldots, x_n \in \Ob(G)} \Hom_{\Vect}(V_{x_0}, O(G(x_0,x_1))\ten \ldots \ten O(G(x_{n_1},x_n))).
\end{eqnarray*}

This has the natural structure of cochain complex (in fact, a cosimplicial complex), with coboundary
\begin{eqnarray*}
df (g_1,\ldots, g_{n+1})&=& g_1 \cdot f(g_2,\ldots, g_{n+1}) \\&&+\sum_{i=1}^{n}f(g_1, \ldots, g_jg_{j+1}, \ldots, g_{n+1}) \\&&+ (-1)^{n+1}f(g_1,\ldots, g_n),
\end{eqnarray*}
for $f\in \CC^n(G,V^{\vee})$, and we define $\H^n(G,V^{\vee})$ to be the $n$th cohomology group of this chain complex. 
\end{definition}

\begin{exercises}\label{gpdex}
\begin{enumerate}
\item If $G \xra{f} K$ is an equivalence of pro-algebraic groupoids, then the maps $f^{\sharp}\co \H^n(K,W) \to \H^n(G,f^*W)$ are isomorphisms.

\item If $G$ is a disjoint union of groups $G(x)$, then $\H^n(G,W) =\prod_{x \in \Ob(G)}\H^n(G(x), W_x)$, where the latter is group cohomology defined in the usual way.
\end{enumerate}
\end{exercises}

\begin{definition}
A morphism $G \xra{f} H$ of pro-algebraic groupoids is said to be a pro-unipotent extension if the following hold:
\begin{enumerate}
\item $\Ob(f)\co \Ob(G) \to \Ob(H)$ is an isomorphism.

\item For all $x,y \in \Ob(G)$ and all $k$-algebras $A$, $f_{xy}\co G(x,y)(A) \to H(fx,fy)(A)$ is surjective (equivalently, $f_{xy}$ has a section).

\item For all $x \in \Ob(G)$, the pro-algebraic group $U(x):=\ker (G(x,x) \to H(fx,fx))$ is pro-unipotent.
\end{enumerate}
Such an extension is called abelian if the groups $U(x)$ are all abelian.
\end{definition}

\begin{proposition}\label{obs}
 If $\theta\co K \to H$ is a morphism of pro-algebraic groupoids, and $f\co G\to H$ an abelian  pro-unipotent extension, then the groups $U(x):=\ker (G(x,x) \to H(fx,fx))$ form a pro-finite-dimensional $H$-representation. The obstruction to lifting $\theta$ to $G$ then lies in $\H^2(K, \theta^*U)$.
\end{proposition}
\begin{proof}
$U(x)$ has the natural structure of a $G$-representation, given by the adjoint action. Since the $U(x)$ are all abelian and $H(x,y)=G(x,y)/U(y,y)$, this action descends to $H$, making $\{U(x)\}$ an $H$-representation.

To see that $\H^2(K, \theta^*U)$ is an obstruction space, choose morphisms $\psi_{xy}\co K(x,y) \to G(x,y)$ of affine schemes lifting the $\theta_{xy}$. Then the obstruction to lifting is the datum 
$$
o_{xyz}(\psi)=\psi_{xy}\psi_{yz}\psi_{xz}^{-1} \in \z^2(G,\theta^*U).
$$ 
A different choice of lift would amount to maps 
$\alpha_{xy}\co K(x,y) \to U(x)$ of affine schemes, with $\psi'_{xy}=\alpha_{xy}\psi_{xy}$. We then see that $o(\psi')=o(\psi)+d\alpha$, so the obstruction lies in 
$$\H^2(K, \theta^*U),$$
as required.
\end{proof}

\begin{proposition}
If $G$ is a reductive pro-algebraic groupoid, then $\H^n(G,V)=0$ for all $n>0$ and all pro-finite-dimensional $G$-representations $V$.
\end{proposition}
\begin{proof}
For some field extension $k \subset K$, every non-empty $G(x,y)$ has a $K$-valued point. Define an equivalence relation on $\Ob(G)$ by $x \sim y$ if $G(x,y)$ is non-empty, and choose a set $S$ of representatives of these equivalence classes. Let $G_{K/k}$ denote the pro-algebraic groupoid over $K$ represented by the $O(G(x,y))\ten_k K$. Then the morphism
$$
\coprod_{s \in S} G_{K/k}(s,s) \to G_{K/k}
$$
is an equivalence, so for all $G_{K/k}$-representations $V$, 
$$
\H^n(G_{K/k},V^{\vee})\cong \prod_{s \in S} \H^n(G_{K/k}(s,s),V^{\vee}_s),
$$
which equals $0$ for $n>0$, since the $G(s,s)$ are reductive and $K$ is of characteristic $0$. 

Now, 
$$
\H^n(G_{K/k},V^{\vee})= \H^n(G,V^{\vee}),
$$
where the dual on the left is as a $K$-vector space, and the dual on the right as a $k$-vector space. 
For any $G$-representation $W$ over $k$, the representation $W\ten_k K$ is therefore flabby, and we can write $W$ as a direct summand $W\ten_k K =W \oplus U$ as a $G$-representation over $k$, so $W$ is flabby, as required.
\end{proof}

\begin{corollary}\label{redlifts}
If $K$ is a reductive pro-algebraic groupoid, and $\theta\co K \to H$ is a morphism of pro-algebraic groupoids, for $f\co G\to H$ a  pro-unipotent extension with kernel $U$, then $\theta$ lifts to $G$, the lift being unique up to conjugation by $U$.
\end{corollary}
\begin{proof}
Let $G_n(x,y)=G(x,y)/\Gamma_{n+1}U(y,y)$. Then $G_n$ is a pro-algebraic groupoid, and $G_n \to G_{n-1}$ is an abelian pro-unipotent extension. By Proposition \ref{obs}, we may lift $\theta$ to each $G_n$ compatibly, and hence to $G = \Lim G_n$, since $G_0=H$, and the second cohomology of $K$ is trivial.

It only remains to show that for any pair $\psi,\phi\co  K \to G$ of lifts, there exist elements $u_x \in U(\theta x,\theta x)$, such that $\psi_{xy}u_y=u_x\phi_{xy}$. Again, we may prove this inductively on $n$, so it suffices to consider the case when $U$ is abelian. As in Proposition \ref{obs}, the difference between $\phi$ and $\psi$ is given by some $\alpha \in \z^1(K,\theta^*U)$. Now, conjugation by $u \in \CC^0(K,\theta^*U)$ just amounts to sending $\alpha$ to $\alpha + du$. But $\H^1(K, \theta^*U)=0$, so there exists such a $u$ with $u\psi u^{-1}=\phi$. 
\end{proof}

\subsection{Levi decompositions}\label{levisn}

\begin{definition}\label{semidirect}
Given a pro-algebraic groupoid $G$, and $U=\{U_x\}_{x \in \Ob(G)}$ a collection of pro-algebraic groups parametrised by $\Ob(G)$, we say that $G$ acts on $U$ if there are morphisms $ U_x\by G(x,y) \xra{*} U_y$ of affine schemes, satisfying the following conditions:
\begin{enumerate}
\item $(uv)*g= (u*g)(v*g)$, $1*g=1$ and $(u^{-1})*g= (u*g)^{-1}$,  for $g \in G(x,y)$ and $u,v \in U_x$.

\item $u*(gh)=(u*g)*h$ and  $u*1=u$, for $g \in G(x,y), h \in G(y,z)$ and $u \in U_x$. 
\end{enumerate}

If $G$ acts on $U$, we write $G \ltimes U $ for the groupoid given by  
\begin{enumerate}
\item $\Ob(G \ltimes U):=\Ob(G)$.

\item $( G \ltimes U)(x,y):= G(x,y)\by U_y$.

\item $(g,u)(h,v):= (gh, (u*h)v)$ for $g \in G(x,y), h \in G(y,z)$ and $u \in U_y, v \in U_z$. 
\end{enumerate}
\end{definition}

\begin{definition}
Given a pro-algebraic groupoid $G$, define $\Ru(G)$ to be the collection $\Ru(G)_x=\Ru(G(x,x))$ of pro-unipotent pro-algebraic groups, for $x \in \Ob(G)$. $G$ then acts on $\Ru(G)$ by conjugation, i.e.
$$
u*g:= g^{-1}u g,
$$
for $u \in \Ru(G)_x$, $g \in G(x,y)$.
\end{definition}

\begin{proposition}\label{leviprop}
For any pro-algebraic groupoid $G$, there is a Levi decomposition $G=G^{\red} \ltimes \Ru(G)$, unique up to conjugation by $\Ru(G)$.
\end{proposition}
\begin{proof}
By Corollary \ref{redlifts}, the pro-unipotent extension $G \to G^{\red}$ admits a section $s$, unique up to conjugation by $\Ru(G)$. If we write $H$ for the subgroupoid $s(G^{\red})$ of $G$ isomorphic to $G^{\red}$, then the $G$-action on $\Ru(G)$ restricts to an $H$-action, and  it only remains to show that the canonical map
\begin{eqnarray*}
H \ltimes \Ru(G) &\xra{f}& G\\
(h,u) &\mapsto& hu
\end{eqnarray*} 
is an isomorphism, since it is clearly a groupoid homomorphism. 

Given $g \in G(x,y)$, let $h \in H(x,y)$ be the image of $g$ under the composition $G \to G^{\red} \xra{s} H$. Then $h^{-1}g \in \Ru(G)_y$, so $f$ is surjective. If $f(h,u)=f(h',u')$, then $hu=h'u'$, so $h'h^{-1} \in \Ru(G)_x$, so $h=h'$ and therefore $u=u'$, as required.
\end{proof}

\begin{definition}
Given a morphism $G \xra{f} H$ of pro-algebraic groupoids, with $H$ reductive, define $G^{\rtimes,f}$ to be the object of $\hat{\cN}(H)$ which is the hull of the functor 
$$
U \mapsto \Hom(G,  H\ltimes U)_f/U,
$$
where $U$ acts by conjugation.
 Note that from the Levi decomposition, it follows that for $H=G^{\red}$ and $f$ the canonical projection, the map $\Ru(G) \to G^{\rtimes,f}$ is an isomorphism. 
\end{definition}

\begin{definition}
Given a morphism $G \xra{f} H$ of pro-algebraic groupoids, with $G$ reductive, let 
$$
\{\Ru(H,f)_x\}_{x \in \Ob(G)}
$$ be the maximal collection of pro-unipotent subgroups of $\{H(fx,fx)\}$ which is normalised by $G$, in the sense that for all $k$-algebras $A$, all $g \in G(x,y)(A)$, and all $u \in \Ru(H,f)_y(A)$, we have 
$$
f(g)uf(g)^{-1} \in \Ru(H,f)_x.
$$ 
In particular, this means that $\Ru(H,f) \in \hat{\cN}(G)$.   Note that in general $\Ru(H) \le \Ru(H,f)$, and that  if $G=H^{\red}$ and $f$ is a choice of Levi decomposition, then $\Ru(H,f)=\Ru(H)$. Observe that 
$$
\Ru(H,f)_x= \Ru(H(fx,fx), f_x),
$$ 
by maximality.
\end{definition}

\subsection{Pro-algebraic simplicial groupoids}\label{sagpdsn}

\begin{definition}
Given a reductive pro-algebraic groupoid $G$, let $\cN(G)$ be the category of $G$-representations in  finite-dimensional nilpotent Lie algebras. Write $\hat{\cN}(G)$ for the category of pro-objects of $\cN(G)$, and $s\hat{\cN}(G)$ for the category of simplicial objects in $\hat{\cN}(G)$. This generalises the definition of \S \ref{unip}, and the definitions and results of that section will be applied for this generality without further comment.
\end{definition}

\begin{definition}
Given a simplicial object $G_{\bullet}$ in the category of pro-algebraic groupoids,   with $\Ob(G_{\bt})$ constant, define the fundamental groupoid $\pi_f(G_{\bullet})$ of $G_{\bt}$ to have objects $\Ob(G)$, and for $x,y \in \Ob(G)$, set 
$$
\pi_f(G)(x,y):= G_0(x,y) / \pd_0N_1(G(y,y)_{\bt}) =  \pd_0N_1(G(x,x)_{\bt})\backslash G_0(x,y),
$$ 
the equality arising since for $f\in G_0(x,y), g \in N_1(G(y,y)_{\bt})$, we have $(\sigma_0f) g(\sigma_0f)^{-1} \in N_1(G(x,x)_{\bt})$, so both equivalence relations are the same. Multiplication and inversion descend similarly.
\end{definition}

\begin{definition}
Define a pro-algebraic simplicial  groupoid  to consist of a simplicial complex $G_{\bullet}$ of pro-algebraic groupoids, such that $\Ob(G_{\bt})$ is constant  and for all $x \in \Ob(G)$, $G(x,x)_{\bt} \in s\agp$, i.e. the maps $G_n(x,x) \to \pi_0(G)(x,x)$ are pro-unipotent extensions of pro-algebraic groups.  We denote the category of pro-algebraic simplicial groupoids by $s\agpd$. 
\end{definition}

\begin{definition}
Define a morphism $f\co G_{\bullet} \to H_{\bullet}$ in $s\agpd$ to be: 
\begin{enumerate}
\item a weak equivalence if the map  $\pi_f(f)\co \pi_f(G_{\bullet})\to \pi_f(H_{\bullet})$ is an equivalence of algebraic groupoids, and the maps  $\pi_n(f,x)\co \pi_n(G_{\bullet}(x,x)) \to \pi_n(H_{\bullet}(fx,fx))$ are isomorphisms for all $n$ and for all $x \in \Ob(G)$. 

\item a fibration if $G(x,x) \to H(fx,fx)$ is a fibration in $s\agp$ for all $x \in \Ob(G)$, and $f$ satisfies the path-lifting condition that for all $x \in \Ob(G), y \in \Ob(H)$, and $h \in H_0(fx,y)(k)$, there exists $z \in \Ob(G)$, $g \in G_0(x,z)(k)$ with $fg=h$.   Equivalently, this says that $G(k) \to H(k)$ is a fibration in the category of simplicial groupoids.

\item a cofibration if   it has the left lifting property (LLP) with respect to all trivial fibrations.
\end{enumerate}
\end{definition}

\begin{lemma}\label{cfgpd} A morphism $G \xra{f} H$ of pro-algebraic simplicial groupoids is
a cofibration if and only if $\Ob(G) \to \Ob(H)$ is injective and 
$G^{\rtimes,r\circ f} \to \Ru(H)$ is a cofibration in $s\hat{\cN}(H^{\red})$, for $r\co H \to H^{\red}$ the canonical projection, and $\Ru(H)\in s\hat{\cN}(H^{\red})$ via some choice of Levi decomposition.
\end{lemma}
\begin{proof}
Assume that we have a commutative diagram
$$
\begin{CD}
G @>n>> G'\\
@VfVV @VVf'V\\
H @>m>> H',
\end{CD}
$$
with $f'$ a trivial fibration, $\Ob(G) \into \Ob(H)$ and $G^{\rtimes,f} \to \Ru(H)$  a cofibration in $s\hat{\cN}(H^{\red})$. 

Since $f'$ is a weak equivalence, for all $y \in \Ob(H')$ there exists $x \in \Ob(G')$ and $h \in H'_0(f'x,y)(k)$. But by the path-lifting property, there then exists $z \in \Ob(G')$ with $f'(z)=y$, so $\Ob(f')$ is surjective. This means that we may lift $\Ob(m)$ to give a map $\Ob(\tilde{m})\co \Ob(H) \to \Ob(G')$, compatible with $\Ob(n)$.

Then define groupoids $G'', H''$, with $\Ob(G'')=\Ob(H'')=\Ob(H)$, and $G''(x,y)=G'(\tilde{m}x,\tilde{m}y)$, $H''(x,y)=H'(mx,my)$. Observe that $G''\to H''$ is a trivial fibration, so we may replace $G'$, $H'$ by $G''$, $H''$, i.e. without loss of generality, $\Ob(m)$ and $\Ob(f')$ are isomorphisms. This means that $\pi_fG' \cong \pi_fH'$, so  $(H')^{\red} \cong (G')^{\red}$, which  means  that $f'$ is now a pro-unipotent extension by an acyclic simplicial group. Therefore $(f')^{-1} \im(m) \to \im(m)$ is a trivial fibration, so we may assume that $m$ is surjective.

After fixing a choice of Levi decomposition for $H$, by Corollary \ref{redlifts} we may lift $m|_{H^{\red}}$ to $G'$.  We may choose the Levi decomposition for $G'$ compatibly, with $\tilde{m}(H^{\red}) = (G')^{\red}$, since $m$ is full. The diagram is now
$$
\begin{CD}
G @>n>> (G')^{\red} \ltimes \Ru(G')\\
@VfVV @VVf'V\\
H^{\red} \ltimes \Ru(H) @>m>> (G')^{\red} \ltimes \Ru(H'),
\end{CD}
$$
with $m$ and $f'$ preserving the Levi decompositions, and giving isomorphisms on objects. Since $H^{\red} \to (G')^{\red}$ is full, this is equivalent to giving the commutative diagram
$$
\begin{CD}
G^{\rtimes,r\circ f} @>n>>  \Ru(G')\\
@VfVV @VVf'V\\
 \Ru(H) @>m>>  \Ru(H')
\end{CD}
$$
in $s\hat{\cN}(H^{\red})$. Observe that $\Ru(G') \xra{f'} \Ru(H')$ is a trivial fibration, so the hypothesis that $G^{\rtimes,r\circ f} \to \Ru(H)$ be a cofibration ensures existence of the lift.

Conversely, assume that $G \xra{f} H$ is a cofibration. Define $K$ to be the groupoid with $\Ob(K)=\Ob(G)$, and $K(x,y)=H(fx,fy)$. Then the natural map $K \to H$ is a trivial fibration, and the LLP for the diagram
$$
\begin{CD}
G @>>> K \\
@VfVV @VVV \\
H @= H
\end{CD}
$$
is equivalent to saying that $\Ob(f)$ is injective.

For every commutative diagram
$$
\begin{CD}
G^{\rtimes,r\circ f} @>n>>  U\\
@VfVV @VVf'V\\
 \Ru(H) @>m>>  V
\end{CD}
$$
in $s\hat{\cN}(H^{\red})$,  with $f'$ a trivial fibration, we    have a commutative diagram in $s\agp$
$$
\begin{CD}
G @>n>> H^{\red} \ltimes U\\
@VfVV @VVf'V\\
H^{\red} \ltimes \Ru(H) @>m>> H^{\red} \ltimes V,
\end{CD}
$$
with $f'$ a trivial fibration. We therefore have a lift $H^{\red} \ltimes \Ru(H) \xra{\tilde{m}} H^{\red} \ltimes U$, such that $f'\tilde{m}=m$ and $\tilde{m}f=n$. Thus $ \tilde{m}(H^{\red})\ltimes U$ is an alternative Levi decomposition for $H^{\red} \ltimes U$, and $\tilde{m}$ therefore provides a lift from $\Ru(H)$ to $U$ in $\cN(H^{\red})$, as required.
\end{proof}

\begin{theorem}\label{cmsagpd}
With the classes of morphisms given above, $s\agpd$ is a closed model category.
\end{theorem}
\begin{proof}
\begin{enumerate}
\item[CM1] As in Theorem \ref{cmsagp}, $s\agpd$ is closed under all limits, defined by the formula $(\Lim G_{\alpha})(A)=\Lim (G_{\alpha}(A))$ for all $k$-algebras $A$. Colimits are defined by taking the colimit $H$ in the category of simplicial objects over pro-algebraic groupoids, then taking pro-unipotent completion of $H$ over $\pi_fH$.

\item[CM2] If $f$ and $g$ are composable morphisms, with any two of $f,g, f\circ g$ weak equivalences, then so is the third, by the long exact sequence of homotopy for simplicial groupoids.

\item[CM3] Retracts of fibrations, cofibrations and weak equivalences are (respectively) fibrations, cofibrations and weak equivalences, by the corresponding properties for $s\hat{\cN}(G)$ and for simplicial groupoids.

\item[CM4] By definition, every cofibration has the left lifting property with respect to all trivial fibrations. We need to show that every trivial cofibration has the left lifting property with respect to all  fibrations. Assume we are given a commutative diagram
$$
\begin{CD}
G @>n>> G'\\
@VfVV @VVf'V\\
H @>m>> H',
\end{CD}
$$
with $f'$ a  fibration, and $f$ a trivial cofibration. 
Since $\pi_fG \to \pi_fH$ is an equivalence, it is full, giving $G^{\red} \to H^{\red}$ (also an equivalence).  

Now for each $y \in \Ob(H)$, there exists $x \in \Ob(G)$ with $H^{\red}(fx,y)(k)$ non-empty. As $\Ob(f)$ is injective, choose $\theta_y \in  H^{\red}_0(fx,y)(k)$ for each $y$, such that $\theta_y =\id_y$ if $y \in \im \Ob(f)$. Fix a Levi decomposition for $G$, lift each $\theta_y$ to $H_0$, and define a Levi decomposition for $H$ by $H^{\red}(y_1,y_2)=\theta_{y_1}^{-1}fG^{\red}(x_1,x_2)\theta_{y_2}$.

By the path-lifting property, we can then choose  $\phi_y \in  G'(nx,z)(k)$ with $f'\phi_y = m\theta_y$, such that $\phi_y=\id_{nx}$ if $y =fx$. We can then extend $\tilde{m}$ to the whole of $H^{\red}$ by sending $h \in H(y_1,y_2)$ to
$$
\phi_{y_1}^{-1} n(\theta_{y_1}h\theta_{y_2}^{-1})\phi_{y_2};
$$ 
this extends $m$ and restricts to $n$.

This amounts to having the commutative diagram
$$
\begin{CD}
G^{\rtimes,f}  @>n>> \Ru(G, \tilde{m}) \\
@VfVV @VVf'V\\
 \Ru(H) @>m>> \Ru(H, m)
\end{CD}
$$
in $s\hat{\cN}(H^{\red})$.  Note that $f'$ is also a fibration here, since for $i>0$, $N_i(\Ru(G, \tilde{m}))=N_i(G)$, and similarly for $H$. Since 
$$
G^{\rtimes,f}(y_1,y_2)=\theta_{y_1}^{-1}\Ru(G)(x_1,x_2)\theta_{y_2} \cong \Ru(G)(x_1,x_2),
$$
the map $G^{\rtimes,f} \to \Ru(H)$ is a trivial cofibration in the model category $s\hat{\cN}(K)$, so the lift exists.

\item[CM5] 
Given an arbitrary morphism $f\co G \to H$ in $s\agp$, we need to show that $f$ can be factorised as:
\begin{enumerate}
\item  $f=pi$, with $p$ a fibration, an $i$ a trivial cofibration. To do this, let 
$$
T =\coprod_{ x \in \Ob(G), y \in \Ob(H)} H_0(fx,y)(k).
$$
 We have a canonical map $a\co \Ob(G) \to T$ given by $x \mapsto \id(fx,fx)$. There is a canonical retraction   $r\co T \to \Ob(G)$, and define $\hat{G}$ to have objects $T$ and morphisms $\hat{G}(s,t)=G(rs,rt)$, with the obvious composition law.  There is also a map $b\co \hat{G} \to H$ defined on objects by $b\theta=y$, for $\theta \in H_0(fx,y)(k)$, and sending $\alpha \in \hat{G}(\theta_1,\theta_2)=G(x_1,x_2)$ to $\theta_1^{-1}f(\alpha)\theta_2$.   

This gives us a factorisation of $f$ as $G \xra{a} \hat{G} \xra{b} H$, with $a$ a trivial cofibration. Now, in $s\hat{\cN}(\hat{G}^{\red})$, choose a factorisation of $\Ru(\hat{G}) \xra{b} \Ru(H,b)$ as
$$
\Ru(\hat{G}) \xra{t} U \xra{p} \Ru(H,b),
$$
for $t$ a trivial cofibration, and $p$ a fibration. Finally, setting $i=ta$, we have factorised $f$ as
$$
G \xra{i} \hat{G}^{\red} \ltimes U \xra{p} H.
$$
Note that $p$ is a fibration, since for $n>0$, $N_n(\hat{G}^{\red} \ltimes U(x,x))=N_n(U(x,x))$, and $N_n(H(px,px))= N_n(\Ru(H,b)(x,x))$. The path-lifting condition is automatically satisfied, by the construction of $b$.

\item $f=qj$, with $q$ a trivial fibration, and $j$ a cofibration. To do this, choose a factorisation $\Ob(G) \xra{a} S \xra{b} \Ob(H)$, with $a$ injective and $b$ surjective.   Define $\tilde{H}$ to have objects $S$, and morphisms $\tilde{H}(s,t) =H(bs,bt)$, with the obvious composition law. This gives us a factorisation of $f$ as $G \xra{a} \tilde{H} \xra{b} H$, with $b$ a trivial fibration. Now, in $s\hat{\cN}(\tilde{H}^{\red})$, choose a factorisation of  $G^{a,\rtimes} \xra{a} \Ru(\tilde{H})$ as 
$$
G^{a,\rtimes} \xra{j}U \xra{t} \Ru(\tilde{H}),
$$
 for $t$ a trivial fibration, and $j$ a cofibration. Finally, setting $q=bt$, we have factorised $f$ as 
$$
G \xra{j} \tilde{H}^{\red} \ltimes U \xra{q} H,
$$ 
as required. 
\end{enumerate}
\end{enumerate}
\end{proof}

\subsection{Algebraisation}

There is a forgetful functor $(k)\co s\agpd \to s\gpd$, given by sending $G_{\bt}$ to $G_{\bt}(k)$. This functor clearly commutes with all limits, so has a left adjoint $G_{\bt} \mapsto (G_{\bt})^{\alg}$. We can describe $(G_{\bt})^{\alg}$ explicitly. First let $(\pi_f(G))^{\alg}$ be the pro-algebraic completion of the abstract groupoid $\pi_f(G)$, then let $(G^{\alg})_n$ be the relative Malcev completion  of the morphism
$$
G_n \to (\pi_f(G))^{\alg}.
$$
In other words, $G_n \to (G^{\alg})_n \xra{f} (\pi_f(G))^{\alg}$ is the universal diagram with $f$ a pro-unipotent extension.

\begin{proposition}\label{algqd}
The functors $(k)$ and ${}^{\alg}$ are a pair of Quillen functors.
\end{proposition} 
\begin{proof}
It suffices to observe that $(k)$ preserves fibrations and trivial fibrations. 
\end{proof}

\begin{definition}
Given a  simplicial set (or equivalently a  topological space), define the pro-algebraic homotopy type of $X$ over $k$ to be the equivalence class of 
$$
G(X)^{\alg}
$$ 
in the closed model category of pro-algebraic simplicial groupoids. 
\end{definition}

\begin{theorem}\label{varpi}
Given a simplicial groupoid $G$, $\pi_f(G^{\alg})$ is the pro-algebraic completion of the fundamental groupoid $\pi_f(G)$.  
\end{theorem}
\begin{proof}
  This  is immediate, since algebraisation preserves coequalisers. 
\end{proof}

\subsection{The homotopy category}

The category of simplicial groupoids is not a simplicial model category, so we cannot expect $s\agpd$ to be so. However, localising at weak equivalences still gives a homotopy category $\Ho(s\agpd)$. In order to describe this explicitly, we need to construct path  objects. 

Recall from \cite{GMO} that a path object $H^I$ in the category of simplicial groupoids is given by
$$
\Ob(H^I) =H_0,\quad (H^I)(a,b)= \HHom_{\bS}(\Delta^1, H(sa,sb)),
$$
where $s,t\co H_n \to \Ob H$ are the source and target maps in $H$. The commutative triangle corresponding to this path object is
$$
\xymatrix{
& H^I \ar[d]^{p=(p_0,p_1)}\\H \ar[ur]^w \ar[r]^{\Delta}& H\by H,
}
$$
with $w$ defined on objects by  $wx=\id_x$, and on morphisms $h \in H_n$  by the composition 
$$
\Delta^1\by \Delta^n \to \Delta^n \xra{h} H(sh,th).
$$ 
For $a \xra{h} b$ in $(H^I)_n$, we define 
$$
p_0(h)=( sa \xra{h_0} sb)  
$$ 
and $p_1$ as the composition
$$
p_1(h)= (ta \xra{(\sigma_0)^na^{-1}} sa  \xra{h_1} ta \xra{(\sigma_0)^n b } tb). 
$$

This the motivates the following:
\begin{lemma}
A path object for $H$ in $s\agpd$ is given by 
$$
\Ob(H^I) =H_0(k),\quad (H^I)(a,b)= \HHom_{\bS}(\Delta^1, H(sa,sb)),
$$
with structure morphisms as for $s\gpd$.
\end{lemma}
\begin{proof}
We first need to prove that this defines an object of $s\agpd$. $H^I$ is clearly a simplicial object in the category of pro-algebraic groupoids, with constant objects, so we need to prove that $(H^I)_n \to \pi_f(H^I)$ is a pro-unipotent extension for all $n$. To see this, observe that 
$$
H^I(a,a)= H(sa,sa)^{\Delta^1}, 
$$
for the latter  defined using the simplicial structure on the category of pro-algebraic simplicial groups. 

We now need to show that $w$ is a weak equivalence, and $p$ a fibration. Observe that $\pi_fH \xra{w} \pi_f(H^I)$ is an equivalence, and that $H(x,y) \to H^I(wx,wy)$ is an isomorphism, so $w$ is a weak equivalence. The map $p$ clearly satisfies the path-lifting condition, so it remains to show that  
$H^I(a,a) \xra{p} H(sa,sa)\by H(ta,ta)$ is a fibration of pro-algebraic simplicial groups for all $a \in H_0(k)$. This follows from the observation that conjugation by $a$ on the second factor gives  the canonical fibration $H(a,a)^{\Delta^1} \to  H(a,a)\by H(a,a)$.
\end{proof}

\begin{proposition}\label{htpymaps}
For $G,H \in s\agpd$, with $G$ cofibrant, 
$$
\Hom_{\Ho(s\agpd)}(G,H)= \Hom_{s\agpd}(G,H)/\sim,
$$
where $\sim$ is  the equivalence relation generated by $F_0\sim F_1$ for all 
$$
F \in \Hom_{s\agpd}(G,\Hom_{\bS}(\Delta^1,H)),
$$
and for $f \in \Hom_{s\agpd}(G,H)$ and $a$ making the following diagram commute
$$
\xymatrix{\Ob G \ar[r]^a \ar[dr]_{\Ob f}& H_0(k)\ar[d] \\ &\Ob H,   }
$$
$f \sim afa^{-1}$.

In particular, if $G,H$ are groups, then
$$
\Hom_{\Ho(s\agpd)}(G,H)= \Hom_{\Ho(s\agp)}(G,H)/H_0,
$$
where $H_0$ acts on the $\Hom$-set by conjugation.
\end{proposition}
\begin{proof}
This follows from the observation that $H$ is automatically cofibrant, and from the description of the path object.
\end{proof}

\section{The Maurer-Cartan space and gauge group}\label{mc}

\subsection{The Maurer-Cartan functor}

\begin{definition}
Given a simplicial set $X$ and a set $S$, define the cosimplicial set $\CC^{\bt}(X,S)$ by
$$
\CC^n(X,S):=S^{X_n}, 
$$
with operations $\pd^i:=S^{\pd_i}, \sigma^i:=S^{\sigma_i}$.
\end{definition}

\begin{definition}
Given a simplicial set $X$ and a simplicial group $G$,  define the Maurer-Cartan space by
$$
\mc(X,G):=\Hom_{\bS}(X,\bar{W}G).
$$
\end{definition}

\begin{lemma}\label{maurercartan}
The Maurer-Cartan space $\mc(X,G)$ consists of sets $\{\omega_n\}_{n\ge 0}$, with $\omega_n \in \CC^{n+1}(X,G_n)$, such that 
\begin{eqnarray*}
\pd_i\omega_n &=& \left\{\begin{matrix} \pd^{i+1}\omega_{n-1}  & i>0 \\ (\pd^1\omega_{n-1})\cdot(\pd^0\omega_{n-1})^{-1} & i=0,\end{matrix} \right.\\
\sigma_i\omega_n &=& \sigma^{i+1}\omega_{n+1},\\
\sigma^0\omega_n&=& 1.
\end{eqnarray*}
\end{lemma}
\begin{proof}
Recall (e.g. from \cite{sht} \S V.5) that the loop group $G(X)$ of a simplicial set $X$ is defined by setting $G_n$ to be the free group $\Fr(X_{n+1})/\sigma_0\Fr(X_n)$.
 The face and degeneracy operators are given by:
\begin{eqnarray*}
\pd_i[x]&=& \left\{\begin{matrix} [\pd_{i+1}x] & i>0 \\ [\pd_1x][\pd_0x]^{-1} & i=0,\end{matrix} \right.\\
\sigma_i[x]&=&[\sigma_{i+1}x].
\end{eqnarray*}
Since the functors $G$ and $\bar{W}$ are adjoint, the result follows from this description. 

Explicitly, by \cite{sht} Lemma V.5.3, $\omega \in \mc(X,G)$ corresponds to the map $X \to \bar{W}G$ defined on $X_{n}$ by
$$
x \mapsto (\omega_{n-1}x, \omega_{n-2}\pd_0x,\ldots,\omega_0\pd_0^{\phantom{0}n-1}x).
$$
\end{proof}

\begin{remark}
If $G$ is a constant group, then this is equivalent to the Maurer-Cartan functor defined in  \cite{paper1} for governing deformations of $G$-torsors. 
\end{remark}

\begin{definition}\label{gauge}
Given $X \in \bS$ and $G \in s\gp$, define the gauge group $\Gg(X,G)$ to  be the subgroup of  $\prod_n \CC^n(G_n)$ consisting of those $g$ satisfying
\begin{eqnarray*}
\pd_ig_n &=& \pd^{i}g_{n-1}  \quad \forall i>0, \\
\sigma_ig_n &=& \sigma^{i}g_{n+1} \quad \forall i.
\end{eqnarray*}
Note that $\CC^0(X,G_0)$ can be regarded as a subgroup of $\Gg(X,G)$, setting $g_n= (\pd^1)^n(\sigma_0)^ng$, for $g\in \CC^0(X,G_0)$.
\end{definition}

\begin{definition}\label{vdef}
Define $V\co  s\gp \to\bS$ by setting 
$$
VG_n:= G_n\by G_{n-1}\by \ldots \by G_0,
$$
with operations
\begin{eqnarray*}
\pd_i(g_n,g_{n-1},\ldots,g_0)&=&(\pd_ig_n,\pd_{i-1}g_{n-1},\ldots, \pd_1g_{n-i+1}, g_{n-i-1},g_{n-i-2},\ldots, g_0)\\
\sigma_i(g_n,g_{n-1},\ldots,g_0)&=&(\sigma_ig_n,\sigma_{i-1}g_{n-1},\ldots, \sigma_0g_{n-i}, g_{n-i},g_{n-i-1},\ldots, x_0).
\end{eqnarray*}

Define $H\co \bS \to s\gp$ by 
$$
H(X)_n=\Fr(X_{n+1}),
$$
with operations defined on a generator $[x] \in H(X)_n$, for $x \in X_{n+1}$ by
 \begin{eqnarray*}
\pd_i[x]&=&  [\pd_{i+1}x] \\
\sigma_i[x]&=&  [\sigma_{i+1}x].
\end{eqnarray*}
\end{definition}

\begin{lemma}\label{hvgauge}
There are natural isomorphisms
$$
\Hom_{\bS}(X,VG)\cong \Gg(X,G)\cong \Hom_{s\gp}(H(X),G).
$$
\end{lemma}
\begin{proof}
We associate to $g \in \Gg(X,G)$ the map from $X$ to $VG$ given by
$$
x \mapsto (g_{n}x, g_{n-1}\pd_0x,\ldots,g_0\pd_0^{\phantom{0}n}x)
$$
for $x \in X_n$.

This corresponds to the map $f\co [x]\mapsto \pd_0 g_{n+1}(x)$, for $x\in X_{n+1}$, $[x] \in H(X)_n$, from which $g$ can be recovered by $g(x)=f[\sigma_0 x]$.
\end{proof}

\begin{definition}\label{defdef}
The gauge group acts on the Maurer-Cartan space, with the action
$
\phi\co  \Gg(X,G)\by \mc(X,G) \to \mc(X,G)
$
 given by 
$$
(g*\omega)_n= (\pd_0g_{n+1}) \cdot \omega_n \cdot (\pd^0g_n^{-1}) .
$$ 
Define the torsor space $\pi(X,G)$ by
$$
\pi(X,G):=\mc(X,G)/\Gg(X,G).
$$
\end{definition}

\begin{proposition}
For $X \in \bS$ and $G \in s\gp$, 
$$
\Hom_{\Ho(\bS)}(X,\bar{W}G)= \pi(X,G).
$$
\end{proposition}
\begin{proof}
Since $X$ is cofibrant and $\bar{W}G$ is fibrant (automatically), we need to show that 
$$
\xymatrix@1{ \bar{W}G \ar[r]^-{(V1,\id,)}&   VG\by \bar{W}G     \ar@<0.5ex>[r]^-{\pr_2} \ar@<-0.5ex>[r]_-{\phi} & \bar{W}G.}
$$
is a path object for $\bar{W}G$ in $\bS$, where
$$
\phi(h_n,h_{n-1},\ldots,h_0,g_{n-1},g_{n-2},\ldots,g_0)= ((\pd_0h_n)g_{n-1}h_{n-1}^{-1},\ldots, (\pd_0h_1)g_0h_0^{-1}).
$$

Now, there is an isomorphism $\psi:VG \to WG$ given by $\psi_n(\underline{h})= (h_n^{-1}, \phi_n(\underline{h},1))$, for $WG$ defined as in  Example \ref{owres}. Thus  $VG \to *$ is a trivial fibration, and $(\id,V1)$  a weak equivalence. 

The proof that every simplicial group is fibrant can be adapted to show that
$$
VG \by \bar{W}G \xra{\pr_2,\phi} \bar{W}G \by \bar{W}G 
$$
is a fibration;
this completes the proof.
\end{proof}

\begin{remark}
Given a simplicial groupoid $G$, we can define a gauge groupoid $\Gg(X,G)$ by the same formulae as for the gauge group. It has objects $\CC^0(X, \Ob G)$. We may then define  the Maurer-Cartan space $\mc(X,G):=\Hom(X,\bar{W}G)$. This has the structure of a $\Gg(X,G)$-representation by similar formulae, and
$$
\Hom_{\Ho(\bS)}(X,\bar{W}G)=\mc(X,G)/\Gg(X,G).
$$
\end{remark}

\subsection{Relative Malcev homotopy types}\label{malcev}

\begin{definition}\label{zardense}
Assume we have an abstract groupoid $G$, a reductive pro-algebraic groupoid $R$, and a representation $\rho\co G \to R(k)$ which is an isomorphism on objects and Zariski-dense on morphisms (i.e. $\rho\co G(x,y) \to R(k)(\rho x, \rho y)$ is Zariski-dense for all $x,y \in \Ob G$). Define the Malcev completion $(G,\rho)^{\mal}$   of $G$ relative to $\rho$  to be the universal diagram
$$
G \to (G,\rho)^{\mal} \xra{p} R,
$$
with $p$  a pro-unipotent extension, and the composition equal to $\rho$. Explicitly, $\Ob(G,\rho)^{\mal}=\Ob G$ and 
$$
(G,\rho)^{\mal}(x,y)=(G(x,x), \rho)^{\mal}\by^{G(x,x)}G(x,y).
$$
If $G$ and  $R$ are groups, observe that this agrees with the usual definition.  

If $\varrho\co G \to R(k)$  is any  Zariski-dense representation (i.e. essentially surjective on objects and Zariski-dense on morphisms) to a reductive pro-algebraic groupoid (in most examples, we take $R$ to be a group), we can define another reductive groupoid $\tilde{R}$ by setting $\Ob \tilde{R}=\Ob G$, and $ \tilde{R}(x,y)=R(\varrho x, \varrho y)$. This gives a representation $\rho\co \pi_fX \xra{\rho} \tilde{R}$ satisfying the above hypotheses, and we define the Malcev completion of $G$ relative to $\varrho$ to be the Malcev completion of $G$ relative to $\rho$. Note that $\tilde{R} \to R$ is an equivalence of pro-algebraic groupoids.
\end{definition}

\begin{definition}
Given a  Zariski-dense morphism $\rho\co \pi_fX \to R(k)$, let the Malcev completion $G(X,\rho)^{\mal}$ of $X$ relative to $\rho$ be the pro-algebraic simplicial group $(G(X), \rho)^{\mal}$.  Observe that the Malcev completion of $X$ relative to $(\pi_fX)^{\red}$ is just $G(X)^{\alg}$. Let $\varpi_f(X^{\rho,\mal}):=\pi_fG(X,\rho)^{\mal}$ and $\varpi_n(X^{\rho,\mal}):=\pi_{n-1}G(X,\rho)^{\mal}$.
\end{definition} 

\begin{definition}
Given a reductive pro-algebraic groupoid $R$,       define $\cE(R)$ to be the full subcategory of $\agpd\da R$ consisting of those morphisms $\rho\co G\to  R$ of  proalgebraic groupoids which are pro-unipotent extensions. Similarly, define $s\cE(R)$ to consist of the pro-unipotent extensions in $s\agpd\da R$, and $\Ho(s\cE(R))$ to be the full subcategory of $\Ho(s\agpd\da R)$ on objects $s\cE(R)$.
\end{definition}

\begin{definition}
Given a pro-algebraic groupoid $R$, define the category $s\cM(R)$ to have the same objects as $s\hat{\cN}(R)$, with morphisms given by
$$
\Hom_{s\cM(R)}(\g,\fh)=\Hom_{\Ho(s\hat{\cN}(R))}(\g,\fh)/\exp(\fh_0^R),
$$
where $\fh_0^R$ is the Lie algebra $\Hom_R(k,  \fh_0)$, acting by conjugation on the set of homomorphisms.
\end{definition}

\begin{proposition}\label{meequiv}
For any reductive pro-algebraic groupoid $R$, the categories $\Ho(s\cE(R))$ and $s\cM(R)$ are equivalent.
\end{proposition}
\begin{proof}
The map $F$ from $s\cM(R)$ to $\Ho(s\agpd\da R)$ is given by 
$$
\g \mapsto R\ltimes \exp(\g).
$$
From the existence of Levi decompositions, it is immediate that $F$ is essentially surjective.
To see that this is full and faithful (and indeed well defined), observe that given a  morphism
$$
f\co R\ltimes \exp(\g) \to R\ltimes \exp(\fh)
$$
in $s\agpd\da R$, there exists a family $a(x) \in \exp(\fh)(x)$ such that $afa^{-1}$ preserves the Levi decomposition. But this is homotopic to $f$ (by Proposition \ref{htpymaps}), so $F$ is full. Now, given
$$
f,g\co \exp(\g) \to  \exp(\fh),
$$ 
we will have $Ff=Fg$ if and only if $f$ is homotopic to $aga^{-1}$ in $\hat{\cN}(R)$, for some inner automorphism $a\in \fh_0$ which preserves the Levi decomposition. This is equivalent to saying that $a \in \fh_0^R$, as required.
\end{proof}

\begin{definition}\label{relmaldef}
We can now define the Malcev homotopy type of $X$ relative to $\rho$ to be the image of   $G(X,\rho)^{\mal}$ in $\Ho(s\cE(\tilde{R}))$, or equivalently $\Ru G(X,\rho)^{\mal}$ in  $s\cM(\tilde{R})$. Since $\tilde{R} \to R$ is an equivalence of groupoids, there is an equivalence $s\hat{\cN}(R) \to s\hat{\cN}(\tilde{R})$, so we may regard the Malcev homotopy type as an object of $s\cM(R)$ (or of $\Ho(s\cE(R))$). It pro-represents the functor
$$
\g \mapsto \Hom_{\Ho(\bS\da BR(k))}(X,\bar{W}(R(k)\ltimes \exp(\g))), 
$$
for $\g \in s\cM(R)$.
\end{definition}

\subsubsection{Classical homotopy groups}\label{clpimal}

\begin{theorem}\label{malcevpi}  
The fundamental groupoid $\pi_f((X,\rho)^{\mal})$ is the relative Malcev completion of $\pi_f\rho\co \pi_fX \to R(k)$. 
\end{theorem}
\begin{proof}
This is immediate, since Malcev completion commutes with quotients. 
\end{proof}

We now extend the notion of algebraically good to the relative Malcev case:
\begin{definition}\label{relgood}
Define a groupoid $\Gamma$ to be good with respect to a Zariski-dense representation $\rho: \Gamma \to R(k)$ to a reductive pro-algebraic groupoid if the map
$$
\H^n(\Gamma^{\rho, \mal}, V) \to \H^n(\Gamma, V)
$$  
is an isomorphism for all $n$ and all finite-dimensional $\Gamma^{\rho, \mal}$-representations $\Gamma$.
\end{definition}

\begin{lemma}\label{goodh}
Assume that for all $x \in \Ob \Gamma$,  $\Gamma(x,x)$ is finitely presented, with $\H^n(\Gamma, - )$ commuting with filtered direct limits of  $\Gamma^{\rho, \mal}$-representations, and $\H^n(\Gamma, V)$ finite-dimensional for all finite-dimensional $\Gamma^{\rho, \mal}$-representations $V$. 

Then $\Gamma$ is good with respect to $\rho$ if and only if for any finite-dimensional $\Gamma^{\rho, \mal}$-representation $V$, and $\alpha \in \H^n(\Gamma, V)$, there exists an injection $f:V \to W_{\alpha}$ of finite-dimensional $\Gamma^{\rho, \mal}$-representations, with $f(\alpha)=0 \in \H^n(\Gamma, W_{\alpha})$.
\end{lemma}
\begin{proof}
As for \cite{schematic}v1 Lemma 4.15.
\end{proof}

\begin{examples}\label{goodexamples}
\begin{enumerate}
\item All finite groups are good with respect to all representations, as are all finitely generated free groups. 

\item\label{ttwo} All finitely generated nilpotent groups are good with respect to all representations.

\item\label{tthree} The fundamental group of a compact Riemann surface is algebraically good with respect to all representations.

\item\label{ffour} If $1 \to F \to \Gamma \to \Pi \to 1$ is an exact sequence of groups, with $F$ finite and $\rho: \Gamma \to R(k)$ Zariski-dense, assume that $\Pi$ is good relative to $R/\overline{\rho(F)}$, where $\overline{\phantom{F}}$ denotes Zariski closure. Then $\Gamma$ is good relative to $\rho$.
 \end{enumerate}
\end{examples}
\begin{proof}
\begin{enumerate}

\item[\ref{ttwo}.] If $\Gamma$ is a finitely generated nilpotent group, then the trivial representation $k$ is the only irreducible $R$-representation $V$ for which $\H^1(\Gamma, V)\ne 0$. Thus $R$ acts trivially on the abelianisation of $\Ru(\Gamma,\rho)^{\mal}$, and hence nilpotently on $\Ru(\Gamma,\rho)^{\mal}$. But $R$ is reductive, so we must have $(\Gamma,\rho)^{\mal}=R \by \Ru(\Gamma,\rho)^{\mal}$. The Malcev completion is then  given by $\Gamma \ten k=\Ru(\Gamma,\rho)^{\mal}$. Since $\Gamma$ is nilpotent, $\H^*(\Gamma\ten k, k)= \H^*(\Gamma, k)$, giving the  required isomorphism on cohomology.

\item[\ref{tthree}.] Although this is motivated by \cite{schematic}v1 Proposition 4.13(2) (which takes $R=\varpi_1(X)^{\red}$), that proof does not carry over. Instead, apply Corollary \ref{kformalall} to see that $\varpi_n(X^{\rho, \mal})= \H_{n-1}(G(\H^*(X, \bO(R))))$, in the notation of Corollary \ref{bigequiv}. Since $\H^2(X,\bO(R))\cong k$, with higher cohomology $0$, $ G(\H^*(X, \bO(R)))$ is the free chain Lie algebra on generators $V:=\H^1(X, \bO(R)))^{\vee}$ in degree $0$ and $k$ in degree $1$, with differential determined by $d:k \to \bigwedge^2V$ dual to the cup product. Since the cup product is non-degenerate, the higher homology groups of this Lie algebra are all $0$ (similarly to \cite{Sullivan} \S 12), so $G(X)^{\rho, \mal} \sim \pi_1(X)^{\rho, \mal}$, and their cohomology groups agree. 

\item[\ref{ffour}.] This is essentially \cite{schematic}v1 Theorem 4.16 Step 1. Observe that $\H^*(\Gamma,V)\cong \H^*(\Pi, V^{F})$, so any inclusion $V^{F} \into W$ which satisfies Lemma \ref{goodh} for $(\Pi,R/\overline{\rho(F)})$ gives $V \into (W\oplus V)/V^F$, which satisfies Lemma \ref{goodh}  for $(\Gamma,R)$. Note that we cannot adapt the arguments of \cite{schematic}v1 Theorem 4.16 to consider instead the case when  $F$ is free, since we do not know that $\H^1(F,V)$ is an extension of $R$-representations whenever $V$ is.  
\end{enumerate}
\end{proof}

\begin{theorem}\label{classicalpimal}
If $X$ is a  topological space with fundamental groupoid $\Gamma$, equipped with a Zariski-dense representation $\rho: \Gamma \to R(k)$ to a reductive pro-algebraic groupoid  for which: 
\begin{enumerate}
\item $\Gamma$ is  good with respect to $\rho$,
\item $\pi_n(X,-)$ is of finite rank for all $n>1$,  and
\item the $\Gamma$-representation  $\pi_n(X,-)\ten_{\Z} k$ is an extension of $R$-representations (i.e. a $\Gamma^{\rho, \mal}$-representation),
\end{enumerate}
then the canonical map
$$
  \pi_n(X,-)\ten_{\Z} k \to \varpi_{n}(X^{\rho,\mal}) 
$$
is an isomorphism for all $n>1$.
\end{theorem}
\begin{proof}
As for Theorem \ref{classicalpi}, replacing $\Gamma^{\alg}$ by $\Gamma^{\rho, \mal}$ throughout.
\end{proof}

\begin{remark}
Combined with the results of \S\ref{forms}, this theorem can be regarded as  generalising both \cite{KTP} Lemma 4.3.2 (which takes $R=(\pi_1X)^{\red}$ and $X$ a projective complex manifold) and \cite{Sullivan} Theorem 10.1 (which takes $R$=1). 
\end{remark}

\subsection{Cosimplicial algebras}

Fix a reductive pro-algebraic groupoid $R$ over $k$. 

\begin{definition}
 Let $c\Alg(R)$ be the category of of $R$-representations in cosimplicial $k$-algebras. Explicitly, an object $A$ of $c\Alg(R)$ is a collection $\{A(x)\}_{x \in \Ob R}$ of cosimplicial $k$-algebras, equipped with $k$-algebra homomorphisms
$$
A(y) \to A(x)\ten_k O(R(x,y)), 
$$
satisfying the obvious coassociativity and coidentity conditions.
\end{definition}

\begin{definition}
A map $f\co A \to B$ in in $c\Alg(R)$ is said to be 
\begin{enumerate}
\item a  weak equivalence if  $\H^i(f)\co \H^i(A) \to \H^i(B)$ is an isomorphism in $\Rep(R)$ for all $i$; 
\item a fibration if $f^i(x)\co A^i(x) \to B^i(x)$ is a surjection for all $x \in \Ob(R)$ and all $i$;
\item a cofibration if it has the LLP with respect to all trivial fibrations.
\end{enumerate}
\end{definition}

\begin{definition}
Define a simplicial structure on $c\Alg(R)$ by setting $(A^K)^i=(A^i)^{K_i}$, for $A \in c\Alg$, $K \in \bS$. 
\end{definition}

\begin{proposition}\label{calgmodel}
With the definitions above,  $c\Alg(R)$ is a simplicial model category.
\end{proposition}
\begin{proof}
In the case $R=1$, this is standard (see e.g.\cite{chaff} \S 2.1). The same proof applies in the general case, since $\Rep(R)$ is a semisimple category. Equivalently, we may apply \cite{Hirschhorn} Theorem 11.3.2 to the forgetful functor $c\Alg(R)\to c\Mod(R)$, obtaining a cofibrantly generated  model structure. 
\end{proof}

\begin{definition}
Let $s\Aff(R)$ denote the category  of simplicial affine schemes over $k$, i.e. the category opposite to $c\Alg(R)$. We give $s\Aff(R)$ the closed model structure opposite to that on $c\Alg(R)$, so that $\Spec B \to \Spec A$ is a fibration if and only if $A \to B$ is a cofibration, and so on.
\end{definition}

\begin{definition}
Define $\bar{W}\co s\hat{\cN}(R) \to s\Aff(R)$ to represent the functor
$$
(\bar{W}\g)(x)(A):= \bar{W}(\exp(\g(x)\ten A)),
$$
from $k$-algebras to simplicial sets, for $x \in \Ob(R)$. Similarly, define 
$$
V\g(x)(A):= V(\exp(\g(x)\ten A)).
$$
Observe that these functors are continuous, so are right adjoints, but that $\bar{W}$ is  not right Quillen, since it does not preserve fibrations.
\end{definition}

\begin{proposition}
$\bar{W}\co s\hat{\cN}(R) \to s\Aff(R)$ descends to a functor
$$
\bar{W}\co \Ho(s\hat{\cN}(R)) \to \Ho(s\Aff(R)).
$$
\end{proposition}
\begin{proof}
The reverse Adams spectral sequence of Theorem \ref{spectralpi} shows that $\bar{W}$ sends weak equivalences to weak equivalences, as required.
\end{proof}

\subsection{Maurer-Cartan for cosimplicial algebras}

\begin{definition}\label{mcc}
Given $A \in c\Alg(R)$, define the Maurer-Cartan functor
\begin{eqnarray*}
\mc(A,-)\co s\hat{\cN}(R) &\to& \Set\\
 \g &\mapsto& \Hom_{c\Alg(R)}(O(\bar{W}\g), A),
\end{eqnarray*}
where $\Spec O(\bar{W}\g)$ denotes global sections of the structure sheaf of  the simplicial affine scheme $\bar{W}\g:=\bar{W}\exp(\g)$.
\end{definition}

\begin{definition}
Given $V,W \in \Rep(R)$, define $V\ten^RW:=\Hom_{\Rep(R)}(k,V\ten W)$.
\end{definition}

\begin{proposition}
Given $A \in c\Alg(R)$ and $\g \in s\hat{\cN}(R)$,  $\mc(A,\g)$ consists of sets $\{\omega_n\}_{n\ge 0}$, with $\omega_n \in \exp(A^{n+1}\hat{\ten}^R\g_n)$, satisfying the equations of Lemma \ref{maurercartan}. 
\end{proposition}
\begin{proof}
The calculation is the same as that from \cite{sht} \S V.5 used in   Lemma \ref{maurercartan}.
\end{proof}

\begin{definition}
Define the loop group functor $G\co s\Aff(R) \to s\hat{\cN}(R)$  to be  the  left adjoint to $\bar{W}$, and define $H$ to be left adjoint to $V$.
\end{definition}

\begin{corollary}
For $A \in c\Alg(R)$, let $\tilde{A}^n=\ker(\sigma^0\co A^{n+1} \to A^n)$. Then $G(A)_n$ is the free Lie algebra $\Lie((\tilde{A}^{n+1})^{\vee})$ generated by the dual of $\tilde{A}^{n+1}$, with operations
\begin{eqnarray*}
\sigma_i&=&(\sigma^{i+1})^*\co G(A)_n \to G(A)_{n+1}\\
\pd_i&=&(\pd^{i+1})^*\co G(A)_n \to G(A)_{n-1} \quad \text{for } i>0,
\end{eqnarray*} 
 and $\pd_0$ is defined on generators by
$$
\pd_0 = (\pd^1)^*\star (-\pd^0)^*\in \Hom_R((\tilde{A}^{n+1})^{\vee}, \Lie((\tilde{A}^{n})^{\vee}))\cong \Lie((\tilde{A}^{n})^{\vee})\hat{\ten}^R\tilde{A}^{n+1},
$$  
where $\star$ denotes the Campbell-Baker-Hausdorff formula on the Lie algebra $\Lie((\tilde{A}^{n})^{\vee})\hat{\ten}^R\tilde{A}^{n+1}$. 
\end{corollary}

\begin{definition}
Given $A \in c\Alg(R)$ and $\g \in s\hat{\cN}(R)$, define the gauge group $\Gg(A,\g)\le \prod_n \exp(A^n\hat{\ten}^R\g_n)$ by the equations of Definition \ref{gauge}. This acts on $\mc(A,\g)$ by the formula of Definition \ref{defdef}, and we define the torsor space by
$$
\pi(A,\g):=\mc(A,\g)/\Gg(A,\g).
$$
\end{definition}

\begin{lemma}
There are natural isomorphisms
$$
\Hom_{s\Aff(R)}(\Spec A,V\g)\cong \Gg(A,\g)\cong \Hom_{s\hat{\cN}(R)}(H(A),\g),
$$
where $H(A)_n$ is the free Lie algebra $\Lie((A^{n+1})^{\vee})$, with operations $\sigma_i= (\sigma^{i+1})^*$, $\pd_i=(\pd^{i+1})^*$, for all $i$.
\end{lemma}
\begin{proof}
As for Lemma \ref{hvgauge}.
\end{proof}

\begin{theorem}\label{mcgaugeworks}
Given $A \in c\Alg(R)$ and $\g \in s\hat{\cN}(R)$,
$$
\pi(A,\g)=\Hom_{\Ho(s\Aff(R))}(\Spec A,\bar{W}\g).
$$
\end{theorem}
\begin{proof}
It follows from Lemma \ref{smoothcrit} that $G$ maps trivial cofibrations to smooth maps, so $\bar{W}$ maps surjections to fibrations, and hence $\bar{W}\g$ is fibrant. We need to show that 
$$
\xymatrix@1{ \bar{W}\g \ar[r]^-{(\id,V1)}&   \bar{W}\g \by V\g   \ar@<0.5ex>[r]^-{\pr_1} \ar@<-0.5ex>[r]_-{\phi} & \bar{W}\g.}
$$
is a path object for $\bar{W}\g$ in $s\Aff(R)$. 

As in \cite{W} Exercise 8.3.7, The map $(\sigma^0)^*$ makes $H(A) \xra{(\pd^1)^*} \Lie((A^0)^{\vee})$  into a left-contractible augmented simplicial Lie algebra (in the sense of \cite{W} 8.4.6), so $H(A)$ is weakly equivalent to $\Lie((A^0)^{\vee})$. This means that $H$ maps cofibrations to smooth maps, so $V$ maps surjections to trivial fibrations, and $(\id,V1)$ is thus a weak equivalence.

To see that $ \bar{W}\g \by V\g \xra{(\pr_1,\phi)}\bar{W}\g \by\bar{W}\g   $ is a fibration, take a commutative diagram
$$
\begin{CD}
\Spec B @>>> \bar{W}\g \by V\g\\
@VfVV @VV(\pr_1,\phi)V\\
\Spec A @>>> \bar{W}\g \by\bar{W}\g,
\end{CD}
$$ 
with $A \xra{f^{\sharp}} B$ a trivial fibration. We need to show that $(\pr_1,\phi)$ has the  RLP with respect to $f$. This is equivalent to taking $\omega,\omega' \in \mc(A,\g)$ and $g \in \Gg(B,\g)$ such that $g(f^*\omega)=f^*\omega'$, and seeking $h \in \Gg(A,\g)$ such that $h(\omega)=\omega'$ and $f^*h=g$. Since $V\g$ is acyclic and fibrant, we know that $\Gg(A,\g) \to \Gg(B,\g)$ is surjective, so taking a lift $\tilde{g}$ of $g$ and replacing $\omega'$ by $\tilde{g}^{-1}\omega'$, we may assume that $g=1$.

Consider the following map  in $\g \da s\hat{\cN}(R)$:
$$
\xymatrix@1{\g*G(B)  \ar@<0.5ex>[r]^-{G(f)} \ar@<-0.5ex>[r]_-{\omega'} & \g*G(A) \ar[rr]^-{(\omega*\id)\circ \Delta} &&\g*(H(A)/H(B)),  }
$$
where $\Delta\co G(A) \to G(A)*H(A)$ corresponds to the action of $\Gg$ on $\mc$. This commutes, and it will suffice for us to show that the map from the coequaliser $C$  of the pair on the left to $\g*(H(A)/H(B))$ is smooth over $\g$, since 
$$
\Hom_{\g \da s\hat{\cN}(R)}(C,\fh)=\{\eta \in \mc(A,\fh)\,:\, f_*\eta=  \omega'\},
$$
and 
$$
\Hom_{\g \da s\hat{\cN}(R)}(\g*(H(A)/H(B),\fh)=\ker( \Gg(A,\fh) \to \Gg(B,\fh)).
$$

 Each entry in the diagram is cofibrant over $\g$, and the maps on relative cotangent spaces over $\g$ are:
$$
\xymatrix@1{\tilde{B}^{\vee}[1]  \ar@<0.5ex>[r]^-{f[1]} \ar@<-0.5ex>[r]_-{0} & \tilde{A}^{\vee}[1] \ar[rr]^-{\id -(\pd^0\sigma^0)^*} && (A^{\vee}/B^{\vee})[1],  }
$$
where $\sigma_i=(\sigma^{i+1})^*$ for all $i$,   $\pd_i=(\pd^{i+1})^*$ for all $i>0$, while $\pd_0=(\pd^1-\pd^0)^*$ for the two complexes on the left, and $\pd_0=(\pd^1)^*$ on the right. Since $f^{\sharp}$ is surjective, the standard smoothness criterion (\cite{higgs} Proposition \ref{higgs-nSSC}) shows that $C$ is cofibrant, with relative cotangent space  $(\tilde{A}^{\vee}/\tilde{B}^{\vee})[1]$. Since this embeds into $(A^{\vee}/B^{\vee})[1]$ under the map $\id -(\pd^0\sigma^0)^*$, the criterion shows that $C \to \g*(H(A)/H(B))$ is cofibrant. By Lemma \ref{smoothcrit}, it suffices to observe that $\id -(\pd^0\sigma^0)^*$ gives an isomorphism on the higher homology groups of the cotangent complexes, all of which are $0$.
\end{proof}

\begin{corollary}\label{factors}
The composition $ s\hat{\cN}(R) \xra{\bar{W}} s\Aff(R) \to \Ho(s\Aff(R))$ factors through $s\cM(R)$.
\end{corollary}
\begin{proof}
By Theorem \ref{spectralpi}, we know that this composition factors through $\Ho(s\hat{\cN}(R))$. Since $\exp(A^0\hat{\ten}^R\g_0)\le \Gg(A,\g)$, Theorem \ref{mcgaugeworks} gives the factorisation through $s\cM(R)$.
\end{proof}

\subsection{Pro-representability}

In this section we give conditions for a covariant functor 
$$
F\co  s\cN(R) \to \Set
$$
to be representable in $\Ho(s\hat{\cN}(R))$.

\begin{definition}
We say that a natural transformation $F \to G$ of functors $F,G\co  s\cN(R) \to \Set$ is:
\begin{enumerate}
\item
 smooth if 
$$
F(\g) \to G(\g)\by_{G(\fh)} F(\fh)
$$
is surjective for all small extensions $\g \onto \fh$;

\item
 unramified if 
$$
F(V) \to G(V) 
$$
is an isomorphism for all  abelian $V$;

\item \'etale if it is smooth and unramified.
\end{enumerate}
\end{definition}

\begin{definition}
Given $L \in s\hat{\cN}(R)$, define $h_L$ by
\begin{eqnarray*}
h_L\co  s\cN(R) &\to& \Set\\
\g &\mapsto& \Hom_{\Ho(s\hat{\cN}(R))}(L, \g).
\end{eqnarray*}
\end{definition}

\begin{definition}
We say that $F$ has a hull $(L,\xi)$ if there exists $L \in s\hat{\cN}(R)$ 
and $\xi \in \hat{F}(L)$ such  that 
$
 h_L \xra{\xi^*} F
$ 
is \'etale.

We say that    $F$ is homotopy pro-representable (by $(L,\xi)$)  if 
$
 h_L \xra{\xi^*} F
$ 
is a natural isomorphism.
\end{definition}

\begin{lemma}\label{hull!}
If $F$ has a hull, then the hull is unique up to (non-unique) weak equivalence.
\end{lemma}
\begin{proof}
If $L,L'$ are hulls for $F$, then the smoothness properties give us maps $\alpha\co L \to L'$, $\beta\co L' \to L$. On abelianisations, the unramification conditions ensure that the maps $\Ho(\alpha^{\ab}), \Ho(\beta^{\ab})$ (on homotopy classes of the abelianisations) are mutually inverse. 
Since $L,L'$ must necessarily be cofibrant, this shows that $\alpha$ is a weak equivalence. 
\end{proof}

The following can be thought of as an analogue of Schlessinger's Theorem (\cite{Sch} Theorem 2.11).
\begin{proposition}\label{schrep}
Assume  $F$ satisfies the following conditions:
\begin{enumerate}
\item[(R1)] The natural map $F(\prod_{i \in I} \g_i) \to \prod_{i \in I} F(\g_i)$ is an isomorphism for any finite set of $\{\g_i\}_{i \in I}$ of objects of $s\cN(R)$ (including the case $I =\emptyset$).

\item[(R2)] The natural map $F(\g\by_{\fh} \fk) \to F(\g)\by_{F(\fh)} F(\fk)$ is  surjective for all small extensions $\g \to \fh$ in $s\cN(R)$.

\item[(R3)] The  map $F(\g) \to F(\fh)$ is surjective for all acyclic small extensions $\g \to \fh$ in  $s\cN(R)$.

\item[(R4)] If  $V$ is  acyclic abelian, then  $F(V)=\bt$, the one-point set.
\end{enumerate} 
Then $F$ has a hull.

If in addition
\begin{enumerate}
\item[(R5)]
$
 F(\g\by_{\fh} \g) \to F(\g)\by_{F(\fh)} F(\g)
$
is an isomorphism for all small extensions $\g \to \fh$ in the subcategory $\cN(R)$, 
\end{enumerate}
then $F$ is homotopy pro-representable.
\end{proposition}
\begin{proof}
First assume that (R1)--(R4) hold, and consider the functor $\tau F\co  s\cN(R) \to \Set$ given by $\tau F(\g):=F(\pi_0\g)$. This can be regarded as a functor on $\cN(R)$, and satisfies the conditions of \cite{higgs} Theorem \ref{higgs-nSch}, so has a hull $\tau L \in \hat{\cN}(R)$. Now  define a functor $\tilde{F}$   by 
$$
\tilde{F}(\g) = \Hom_{\cN(R)}(\tau L, \pi_0\g)\by_{\tau F(\g)}F(\g),
$$
and observe that $\tilde{F} \to F$ is \'etale, so a pro-representation for $\tilde{F}$ will be a hull for $F$.
If (R5) holds, then $\tau L$ pro-represents $\tau F$, so $\tilde{F}$ is naturally isomorphic to $F$. It therefore suffices to show that $\tilde{F}$ is pro-representable, so we replace $F$ by $\tilde{F}$.

Now define 
$$
\hat{F}\co s\hat{\cN}(R) \to \Set
$$
by sending an inverse system $\{\g_i\}_{\bI}$ to $\Lim_{\bI} F(\g_i)$.
Observe that 
$$
\eta\co \hat{F}(\g\by_{\fh} \fk) \to \hat{F}(\g)\by_{\hat{F}(\fh)} \hat{F}(\fk)
$$
 is surjective for all morphisms in $\cN(R)$, as $\pi_0F$ is pro-represented by $\pi_0L$.
Since small extensions and morphisms in $\cN(R)$ are the  generating  fibrations, $\eta$ is therefore surjective for all fibrations $\g \to \fh$ in $s\hat{\cN}(R)$. 

Given a small extension $\g \to \fh$ with kernel $I$, note that $\g\by_{\fh} \g \cong \g \by I$, so (R1) and (R2) combine to give
$$
F(\g)\by F(I) \onto  F(\g)\by_{F(\fh)}F(\g).
$$
Therefore $F(\g) \to F(\fh)$ is injective if  $F(I)=\bt$. Combined with (R3), we then see that $F(f)$ is an isomorphism for all 
acyclic small extensions $f$, and hence $\hat{F}(f)$ is an isomorphism  for all trivial fibrations $f$. The factorisation property of model categories then implies that $\hat{F}(f)$ must be an isomorphism for all weak equivalences $f$.

We have now shown that $\hat{F}$ satisfies the conditions of \cite{jarrep} Theorem 8 (after re-expressing the theorem in terms of covariant, rather than contravariant, functors). This gives the required representability, since $(s\hat{\cN}(R))^{\opp}$ also satisfies the conditions required for that theorem.
\end{proof}

\subsection{The equivalence}

\begin{definition}
Let $c\Alg(R)_0$ be the full subcategory of $c\Alg(R)$ whose objects satisfy $\H^0(A) \cong k$. Let $\Ho(c\Alg(R)_0)$ be the full subcategory of $\Ho(c\Alg(R)_0)$ with objects in $c\Alg(R)_0$. Let $s\Aff(R)_0$ be the opposite category to $c\Alg(R)_0$, and   $\Ho(s\Aff(R)_0)$ opposite to $\Ho(c\Alg(R)_0)$.
\end{definition}

From Section \ref{malcev}, we know that we can regard the relative Malcev homotopy type of $X$ over $R$ as an object of the quotient category $s\cM(R)$ of $\Ho(s\hat{\cN}(R))$.
By Corollary \ref{factors}, we also know that $\bar{W}\co s\hat{\cN}(R)\to s\Aff(R)_0$ descends to a functor $\bar{W}\co s\cM(R) \to \Ho(s\Aff(R)_0)$. The purpose of this section is to show that the map is, in fact, an equivalence. 

\begin{lemma}
For $A \in c\Alg(R)  $, the functor  $\pi(A, -)$ has a hull.
\end{lemma}
\begin{proof}
Apply Proposition \ref{schrep}.
\end{proof}

\begin{definition}\label{qdef}
Given $A \in c\Alg(R)  $ define $\bar{G}(A)\in s\hat{\cN}(R)$ to be the hull of $\pi(A, -)$.
\end{definition} 

\begin{proposition}\label{qwell}
If $A \in c\Alg(R)_0$,  then for all $\g \in s\hat{\cN}(R)$, the map
$$
\Hom_{s\cM(R)}(\bar{G}(A), \g)  \to \pi(A,\g)
$$
is an isomorphism, so $\bar{G}$ defines a functor from $\Ho(s\Aff(R)_0)$ to $s\cM(R)$.
\end{proposition}
\begin{proof}
First note that the map is well defined, since $\exp(\g_0^R)= \exp(\g_0\hat{\ten}^R\H^0(A)\le \Gg(A,\g)$. Since $\bar{G}(A)$ is a hull for both functors, it suffices to show that
$$
\theta\co  \tau\Hom_{\cM(R)}(\bar{G}(A), -) \to \tau\pi(A,-)
$$
is an isomorphism, as
we may apply the proof of Proposition \ref{schrep}, noting that for a functor $F$ with hull $L$, 
$$
F=h_{L}\by_{h_{\tau L}}\tau F.
$$ 

Since $\tau\bar{G}(A)=\pi_0 \bar{G}(A)$ is a hull for $\pi(A,-)$, $\theta$ is automatically surjective. for injectivity, take 
$\omega, \omega' \in \Hom_{\hat{\cN}(R)}(\bar{G}(A),\g)$, and define a functor $ F\co \g \da \hat{\cN}(R)\to \Set$ by
$$
\fh \mapsto \{h \in \exp(A^0\hat{\ten}^R \fh)\,:\, \ad_h\omega=\omega'\},
$$ 
and define the subfunctor $E\co \g \da \hat{\cN}(R)\to \Set$ by 
$$
\fh \mapsto \{h \in \exp(\H^0(A)\hat{\ten}^R\fh)\,:\, \ad_h\omega=\omega'\}.
$$
Using the terminology of \cite{higgs}\S 2, $E$ and $F$ both have  tangent space $\H^0(A)$. 
 Both functors have obstruction space $H^1(A)$, so $E \to F$ is an isomorphism of deformation functors, by the Standard Smoothness Criterion (\cite{higgs} Proposition \ref{higgs-nSSC}),  as required.
\end{proof}

\begin{proposition}\label{wequiv}
The functor $\bar{W}\co s\cM(R) \to \Ho(s\Aff(R))_0$  is an equivalence (with quasi-inverse $\bar{G}$).
\end{proposition}
\begin{proof}
The unit of the adjunction
$$
\xymatrix@1{ s\Aff(R)  \ar@<1ex>[r]^{G}_{\bot}  & s\hat{\cN}(R) \ar@<1ex>[l]^{\bar{W}}}
$$
is a natural transformation $\id \to \bar{W}G$ on $s\Aff(R)$, which combines with the quotient maps $G \to \bar{G}$ to give a  transformation $\id \to \bar{W}\bar{G}$ on $\Ho(s\Aff(R))_0$. The co-unit of the adjunction is a natural transformation $G\bar{W} \to \id$ on $s\hat{\cN}(R)$, which gives rise to  a transformation $\bar{G}\bar{W} \to \id$ in $s\cM(R)$.

By Lemma \ref{owgpcoho}, $\H^i(O(\bar{W}\g))$ is dual to $\H_i(\exp(\g), k)$. By   Proposition \ref{freecoho}, 
$$
\H_i(\bar{G}A,k)= \left\{ \begin{matrix} k &=\H^0(A)^{\vee} & i=0\\ \H_{i-1}( (\bar{G}A)^{\ab}) &= \H^i(A)^{\vee} & i>0, \end{matrix} \right. 
$$
since $\bar{G}A$ is cofibrant.
Thus both transformations defined above give isomorphisms in the homotopy categories.
\end{proof}

\subsection{Equivariant cochains}
\begin{definition}
Given a groupoid $\Gamma$, let $\bS(\Gamma)$ be the category of functors from $\Gamma$ to $\bS$. This has a model category structure given in \cite{sht} \S VI.4. 
\end{definition}

\begin{definition}
Observe that the structure ring $O(R)$ of $R$ is an $R\by R$-representation in algebras, with $O(R)(a,b)$ the fibre over $(a,b) \in \Ob(R \by R)$.
\end{definition}

\begin{definition}\label{bo}
Given $X \in \bS(R(k))$, let $\CC^{\bt}(X,\bO(R)) \in c\Alg(R)$, the algebra of equivariant cochains, be given by
$$
\CC^n(X,\bO(R))(a):= \Hom_{R(k)}(X_n, O(R)(a,-)),
$$
where $O(R)(a,-)$ is endowed with the canonical right $R$-action, and $\CC^n(X,\bO(R))(a)$ has a 
canonical left $R$-action.
\end{definition}

\begin{lemma}
The  functor $s\Aff(R) \to \bS(R(k))$ given by $X \mapsto X(k)$, is right adjoint to the functor  given by 
$$
X \mapsto \Spec \CC^{\bt}(X,\bO(R)).
$$
These form a Quillen pair.
\end{lemma}
\begin{proof}
The adjunction follows from the characterisation of $R$-representations as $O(R)$-comodules, so that there is a canonical isomorphism of $R$-representations
$$
\Hom_R(V,O(R)) \cong V^{\vee}. 
$$
It is immediate that $X \mapsto \Spec \CC^{\bt}(X,\bO(R))$ preserves cofibrations and trivial cofibrations, so is left Quillen. 
\end{proof}

\begin{lemma}\label{wworks}
There is a commutative diagram of functors
$$
\xymatrix@=2ex{
              & s\Aff(R) \ar[rd]^-{(k)} &
\\
s\hat{\cN}(R)\ar[ur]^-{\bar{W}} \ar[dr]^-{\exp(-)(k)} \ar[dd]_{R\ltimes\exp}&         & \bS(R(k)) \ar[dd]^{\substack{{\holim} \\ \lra} \scriptscriptstyle{R(k)}}
\\ 
& s\gp(R(k))\ar[dd]|{R(k)\ltimes}\ar[ur]^-{\bar{W}} &
\\
 s\agpd \da R \ar[rd]^-{(k)}&         & \bS\da BR(k)\\
              & s\gpd\da R(k)\ar[ur]^-{\bar{W}},&
}
$$
where the right-hand  arrow is (as  in \cite{sht} \S VI.4)
$$
(\holim_{\substack{\lra \\ R(k)}}X )_n = \coprod_{\substack{a(0)\to a(1) \to \ldots \to a(n)\\ \text{in } R(k)}} X(a(0))_n.
$$
\end{lemma}
\begin{proof}
For a $\Gamma$-representation $U$ in simplicial groups, an element of 
$$
(\holim_{\substack{\lra \\ R(k)}}\bar{W}U )_n
$$
consists of a string $a(0)\xra{r_{1}} a(1) \xra{r_{2}} \ldots \xra{r_n} a(n)$, and 
$$
(u_{n-1},u_{n-2}, \ldots, u_0) \in U(a_0)_{n-1}\by U(a_0)_{n-2}\by \ldots \by U(a_0)_{0}.
$$
Mapping this  to 
$$(  r_nr_{n-1}\cdots r_1 u_{n-1} r_1^{-1}r_2^{-1}\cdots r_{n-1},\ldots,  r_2r_1 u_{1}r_1^{-1}, r_{1}u_{0} ) \in \bar{W}(\Gamma\ltimes U)
$$
gives an isomorphism  
$$
\bar{W}(\Gamma \ltimes U) \cong \holim_{\substack{\lra \\ \Gamma}}\bar{W}U.
$$

The other squares commute, trivially.
\end{proof}

\begin{definition}
Let  $X$ be a simplicial set, $R$ a reductive pro-algebraic groupoid, and  $\rho\co\pi_fX \to R(k)$ any  representation, then define $\CC^{\bt}(X,\bO(R)) \in c\Alg(R)$, the algebra of equivariant cochains, by
$$
\CC^{\bt}(X,\bO(R)):=\CC^{\bt}(\widetilde{X},\bO(R)),
$$
where $\widetilde{X}$ is the covering system of $\rho$, left adjoint to  $\mathrm{ho}\!\varinjlim_{R(k)}$,   as defined in \cite{sht} \S VI.4.
\end{definition}

\begin{theorem}\label{eqhtpy}
Let  $X$ be a simplicial set, $R$ a reductive pro-algebraic groupoid, and  $\rho\co \pi_fX \to R(k)$ a  representation which is essentially surjective on objects and Zariski-dense on morphisms. Then the Malcev homotopy type of $X$ relative to $\rho$ is given by 
$$
\Spec \CC^{\bt}(X, \bO(R)) \in s\Aff(R)_0,
$$
or equivalently
$$
R\ltimes \bar{G}(\CC^{\bt}(X, \bO(R))) \in \Ho(s\cE(R)),
$$
where $\bar{G}$ is the functor given in Definition \ref{qdef}.
\end{theorem}
\begin{proof}
By definition, the Malcev homotopy type in $s\cM(R)$ represents the functor
$$
\g \mapsto \Hom_{\Ho(s\gpd\da R(k))}(G(X), (R\ltimes\exp(\g))(k)). 
$$ 
Now, we have isomorphisms
\begin{eqnarray*}
\Hom_{\Ho(s\gpd\da R(k))}(G(X), (R\ltimes\exp(\g))(k)) &\cong& \Hom_{\Ho(\bS\da BR(k))}(X,\bar{W}((R\ltimes\exp(\g))(k)))\\
&\cong& \Hom_{\bS(R(k))}(\widetilde{X}, \bar{W}\exp(\g)(k))\\
&\cong& \Hom_{s\Aff(R)}(\Spec \CC^{\bt}(X,\bO(R)), \bar{W}\g),
\end{eqnarray*}
coming from the Quillen pairs, since $X$ is automatically cofibrant, and $\g$ automatically fibrant, so $\bL G(X)=G(X)$, $\bar{W}\g=\R\bar{W}\g$, and so on. 

Now, the hypotheses on $\rho$ ensure that $\H^0(X,\bO(R))=k$, so
\begin{eqnarray*}
\Hom(\Ru G(X,\rho)^{\mal},\g)&\cong& \Hom_{s\Aff(R)}(\Spec \CC^{\bt}(X,\bO(R)), \bar{W}\g)\\
&\cong & \Hom_{s\cM(R)}(\bar{G}(\CC^{\bt}(X,\bO(R))), \g),
\end{eqnarray*}
by Proposition \ref{wequiv}. Thus the relative Malcev homotopy type is
$$
\Ru G(X,\rho)^{\mal}\cong \bar{G}(\CC^{\bt}(X,\bO(R))) \in s\cM(R),
$$
or equivalently $\CC^{\bt}(X,\bO(R))\in \Ho(c\Alg(R)_0)$.
\end{proof}

\begin{remark}
This means that taking $R=(\varpi_fX)^{\red}$ allows us to recover the pro-algebraic homotopy type from the $(\varpi_fX)^{\red}$-representation $\CC^{\bt}(X,\bO((\varpi_fX)^{\red}))$ in cosimplicial algebras.
\end{remark}

\begin{corollary}\label{eqtoen}
Pro-algebraic homotopy types are equivalent to the schematic homotopy types of \cite{chaff}, in the sense that the full subcategory  on objects $X^{\sch}$ of the homotopy category $\Ho(s\mathrm{Pr})$ of simplicial presheaves is equivalent to the full subcategory of $\Ho(s\agpd)$ on objects $G(X)^{\alg}$. Under this correspondence, pro-algebraic homotopy groups and schematic homotopy groups are isomorphic.
\end{corollary}
\begin{proof}
The functor $\bar{W}\co s\agpd \to s\mathrm{Pr}$ has a left adjoint, which we denote by $G$. Since $(\bar{W}G)(k)=\bar{W}(G(k))$, $G(X^{\sch})=G(X)^{\alg}$. 

In \cite{schematic}v1 Corollary 3.22,  it is shown that $\CC^{\bt}(X,\bO(\varpi_f(X,x)^{\red}))$ determines the schematic homotopy type of a pointed connected topological space $(X,x)$, which is given by the simplicial affine scheme
$$
(X,x)^{\sch}= \holim_{\substack{\lra \\ \varpi_f(X,x)^{\red}}}\Spec \CC^{\bt}(X,\bO(\varpi_f(X,x)^{\red})).
$$
As a simplicial presheaf, this is weakly equivalent to
$$
\holim_{\substack{\lra \\ \varpi_f(X)^{\red}}}\Spec \CC^{\bt}(X,\bO(\varpi_f(X)^{\red})), 
$$
which by Theorem \ref{eqhtpy} is equivalent to
$$
\bar{W}G(X)^{\alg}.
$$
Note that this is not the same as the pointed pro-algebraic homotopy type, since the Levi decomposition is only unique up to homotopy, not up to pointed homotopy.   

Since both the pro-algebraic and schematic homotopy types preserve arbitrary disjoint unions, it follows that $X^{\sch} \to \bar{W}G(X)^{\alg}$ is a weak equivalence in $s\mathrm{Pr}$, for any simplicial set $X$. This proves that 
$$
\xymatrix@1{ s\agpd \ar@<0.5ex>[r]^-{\bar{W}} &s\mathrm{Pr}\ar@<0.5ex>[l]^-{G}}
$$
is full and faithful on the objects $X^{\sch},G(X)^{\alg}$.
\end{proof}

\section{Relation with differential forms}\label{forms}
Fix a reductive pro-algebraic groupoid $R$.

\subsection{Chain Lie algebras}
This section is a summary of results we will use from \cite{QRat} \S I.4, modified slightly to take into account our Artinian hypotheses on $\cN(R)$, and replacing vector spaces by $R$-representations.

\begin{definition}
Define $dg\cN(R)$ to be the category of $R$-representations in finite-dimensional nilpotent non-negatively graded chain Lie algebras. Let $dg\hat{\cN}(R)$ be the category of pro-objects in the Artinian category $dg\cN(R)$.

Here, a chain Lie algebra is  a  chain complex $\g=\bigoplus_{i \in \N_0} \g_i$ over $k$, equipped with a bilinear Lie bracket $[,]\co \g_i \by \g_j\ra \g_{i+j}$,  satisfying:

\begin{enumerate}
\item $[a,b]+(-1)^{\bar{a}\bar{b}}[b,a]=0$,

\item $(-1)^{\bar{c}\bar{a}}[a,[b,c]]+ (-1)^{\bar{a}\bar{b}}[b,[c,a]]+ (-1)^{\bar{b}\bar{c}}[c,[a,b]]=0$,

\item $d[a,b] = [da,b] +(-1)^{\bar{a}}[a,db]$,
\end{enumerate}
where $\bar{a}$ denotes the degree of $a$, mod $2$, for $a$ homogeneous.

Define a small extension in $dg\cN(R)$ to be a surjective map $\g \to \fh$ with kernel $I$, such that $[\g,I]=0$. Note that the objects of $dg\cN(R)$ are cofinite in $dg\hat{\cN}(R)$ in the sense of \cite{Hovey} Definition 2.1.4.   
\end{definition}

\begin{lemma}\label{cmdgn}
There is a closed model structure on $dg\hat{\cN}(R)$ in which a morphism $f\co \g \to \fh$ is:
\begin{enumerate}
\item a weak equivalence if $\H_i(f)\co \H_i(\g) \to \H_i(\fh)$ is an isomorphism in $\widehat{\FD\Rep}(R)$ for all $i$;

\item a fibration if $f_i\co \g_i \to \fh_i$ is a surjection for all $i>0$;

\item a cofibration if it has LLP with respect to all trivial fibrations.
\end{enumerate}
\end{lemma}
\begin{proof}
As for \cite{QRat} Theorem II.5.1. Alternatively, we may use \cite{Hovey} Theorem 2.1.19 to show that this is a fibrantly cogenerated model category.  Acyclic small extensions in $dg\cN(R)$ are the generating acyclic fibrations, while small extensions in $dg\cN(R)$ together with arbitrary maps in $\cN(R)$ give the generating fibrations.
\end{proof}

\begin{definition}
As in \cite{QRat} Theorem I.4.4, we say that a morphism $f\co \g \to \fh$ in $dg\hat{\cN}(R)$ is free if there exists a (pro-finite-dimensional)  sub-$R$-representation $V \subset \fh$ such that  $\fh$ is the free pro-nilpotent graded Lie algebra over $\g$ on generators $V$.  
\end{definition}

\begin{lemma}\label{dgfree}
A map  $f\co \g_{\bt} \to \fh_{\bt}$ is a cofibration if and only if it is free.
\end{lemma}
\begin{proof}
As for \cite{QRat} Proposition II.5.5.
\end{proof}

\begin{proposition}
  These classes of morphisms define a closed model category structure on $dg\hat{\cN}(R)$.
\end{proposition}
\begin{proof}
This is essentially the same as \cite{QRat} \S II.5, noting that  a morphism is a fibration of the underlying sets if and only if it is a fibration of the underlying vector spaces, which is equivalent to the condition given.
\end{proof}

\begin{proposition}
If $\g,\fh$ are cofibrant objects in $dg\hat{\cN}(R)$, then a morphism $f\co \g \to \fh$ is a weak equivalence if and only if $f^{\ab}\co  \g/[\g,\g] \to \fh/[\fh,\fh]$ is a quasi-isomorphism of simplicial vector spaces.
\end{proposition}
\begin{proof}
As for Proposition \ref{htpycoho}.
\end{proof}

\begin{proposition}[Minimal models]\label{dgminimal}
Every weak equivalence class in $dg\hat{\cN}(R)$ has a cofibrant element $\m$, unique up to non-unique isomorphism, with $d=0$ on  the abelianisation $\m/[\m,\m]$.
\end{proposition}
\begin{proof}
This is much the same as Proposition \ref{sminimal} (without needing to use the Dold-Kan correspondence). 
\end{proof}

\begin{definition}
Let $k[t,dt]$ be the  $\Z$-graded chain algebra over $k$, freely generated by $t$ in degree $0$, with $dt$ in degree $-1$, so $(dt)^2=0$. Consider the completion $R$ of this algebra with respect to the ideal $J=(t(t-1), dt)$. Note that the quotients $R/(J^n+tJ^{n-1})$ are quasi-isomorphic to $k$.
 Given $\g \in dg\hat{\cN}(R)$, define
$$
\g^I:= \tau_{\ge 0}(\g\hat{\ten}R),
$$ 
where $\tau$ denotes good truncation. There is a natural map $\g \to \g^I$ in $dg\hat{\cN}(R)$, and two maps $\g^I \to \g$ given by sending $t$ to $0$ and $1$ respectively. 
\end{definition}

\begin{lemma}\label{dgpath}
With these structures, $\g^I$ is a path object for $\g$ in $dg\hat{\cN}(R)$.
\end{lemma}

\begin{definition}
Let $dg\cM(R)$ be the category with the same objects as $dg\hat{\cN}(R)$, and morphisms given by
$$
\Hom_{dg\cM(R)}(\g,\fh)= \Hom_{\Ho(dg\hat{\cN}(R))}(\g,\fh)/\exp(\fh^R_0).
$$
\end{definition}

\begin{definition} As in \cite{QRat} I.4.2--4.3, we can define a normalisation functor $N\co s\cN(R) \to dg\cN(R)$, with the Lie bracket the normalised complex $N(\g)$ given by composing with the Eilenberg-Zilber shuffle product:
$$
\llbracket a,b \rrbracket:=\sum_{\begin{smallmatrix} p+q=n\\ (\mu, \nu) \in \mathrm{Sh}(p,q)\end{smallmatrix}}(-1)^{(\mu,\nu)} [\sigma_{\nu_q}\ldots \sigma_{\nu_1}a,\sigma_{\mu_p}\ldots \sigma_{\mu_1}b],
$$
where $\mathrm{Sh}(p,q)$ denotes the set of $(p,q)$ shuffle permutations
 and $(-1)^{(\mu,\nu)}$ is the sign of the permutation $(\mu,\nu)$.
Since the normalisation functor preserves finite limits, it extends to a functor $N\co s\hat{\cN}(R) \to dg\hat{\cN}(R)$.
\end{definition}

\begin{proposition}\label{nequiv}
The functor $N\co s\hat{\cN}(R) \to dg\hat{\cN}(R)$ is a right Quillen equivalence. Thus $\Ho(s\hat{\cN}(R))\simeq \Ho(dg\hat{\cN}(R))$, and $s\cM(R) \simeq dg\cM(R)$.
\end{proposition}
\begin{proof}
Similar to \cite{QRat} Theorem II.5.4.
\end{proof}

\begin{definition}
As in \cite{W} Lemma 8.3.7, we can define a chain complex  $\bar{N}(\g)_n=\g_n/\sum \sigma_i \g_{n-1}$, which is naturally isomorphic to $N(\g)$. 
\end{definition}

\begin{lemma}\label{barng}
The subspace $\sum \sigma_i \g \le \g$ is an ideal for the graded Lie bracket $\llbracket,\rrbracket$. Thus the bracket descends to $\bar{N}(\g)$, and  $N(\g) \to \bar{N}(\g)$ is an isomorphism of chain Lie algebras. 
\end{lemma}
\begin{proof}
The calculation that  $\sum \sigma_i \g$ is an ideal is standard (and the same as that given in \cite{paper1} Lemma 4.16).
\end{proof}

\begin{definition}
We say that a natural transformation $F \to G$ of functors $F,G\co  dg\cN(R) \to \Set$ is:
\begin{enumerate}
\item
 smooth if 
$$
F(\g) \to G(\g)\by_{G(\fh)} F(\fh)
$$
is surjective for all small extensions $\g \onto \fh$;

\item
 unramified if 
$$
F(V) \to G(V) 
$$
is an isomorphism for all  abelian $V$;

\item \'etale if it is smooth and unramified.
\end{enumerate}
\end{definition}

\begin{definition}
Given $L \in dg\hat{\cN}(R)$, define $h_L$ by
\begin{eqnarray*}
h_L\co  dg\cN(R) &\to& \Set\\
\g &\mapsto& \Hom_{\Ho(dg\hat{\cN}(R))}(L, \g).
\end{eqnarray*}
\end{definition}

\begin{definition}
We say that $F$ has a hull $(L,\xi)$ if there exists $L \in dg\hat{\cN}(R)$ 
and $\xi \in \hat{F}(L)$ such  that 
$
 h_L \xra{\xi^*} F
$ 
is \'etale.

We say that    $F$ is homotopy pro-representable (by $(L,\xi)$)  if 
$
 h_L \xra{\xi^*} F
$ 
is a natural isomorphism.
\end{definition}

\begin{lemma}\label{dghull!}
If $F$ has a hull, then the hull is unique up to (non-unique) weak equivalence.
\end{lemma}
\begin{proof}
As for Lemma \ref{hull!}
\end{proof}

\begin{proposition}\label{dgschrep}
Assume  $F$ satisfies the following conditions:
\begin{enumerate}
\item[(R1)] The natural map $F(\prod_{i \in I} \g_i) \to \prod_{i \in I} F(\g_i)$ is an isomorphism for any set $\{\g_i\}_{i \in I}$ of objects of $dg\hat{\cN}(R)$ (including the case $I =\emptyset$).

\item[(R2)] The natural map $F(\g\by_{\fh} \fk) \to F(\g)\by_{F(\fh)} F(\fk)$ is  surjective for all small extensions $\g \to \fh$ in $dg\cN(R)$.

\item[(R3)] The  map $F(\g) \to F(\fh)$ is surjective for all acyclic small extensions $\g \to \fh$ in  $dg\cN(R)$.

\item[(R4)] If  $V$ is  acyclic abelian, then  $F(V)=\bt$, the one-point set.
\end{enumerate} 
Then $F$ has a hull.

If in addition
\begin{enumerate}
\item[(R5)]
$
 F(\g\by_{\fh} \g) \to F(\g)\by_{F(\fh)} F(\g)
$
is an isomorphism for all small extensions $\g \to \fh$ in the subcategory $\cN(R)$, 
\end{enumerate}
then $F$ is homotopy pro-representable.
\end{proposition}
\begin{proof}
As for Proposition \ref{schrep}.
\end{proof}

\subsection{Cochain algebras}

This section summarises   standard results on cochain algebras, stated here in the generality we need (if $R$ is a group, these can all be found in \cite{schematic}v1 \S 4.1).

\begin{definition}
Define $DG\Alg(R)$ to be the category of $R$-representations in  non-negatively graded cochain algebras. 
Here, a cochain  algebra is  a  cochain complex $A=\bigoplus_{i \in \N_0} A^i$ over $k$, equipped with an associative  product $A^i \by A^j\ra A^{i+j}$,  satisfying:

\begin{enumerate}
\item $ab=(-1)^{\bar{a}\bar{b}}ba$,

\item $d(ab) = (da)b +(-1)^{\bar{a}}a(db)$,
\end{enumerate}
where $\bar{a}$ denotes the degree of $a$, mod $2$, for $a$ homogeneous, and a multiplicative identity $1 \in A^0$.
\end{definition}

\begin{lemma}\label{dgalgmodel}
There is a closed model structure on $DG\Alg(R)$ in which a morphism $f\co A \to B$ is:
\begin{enumerate}
\item a weak equivalence if $\H^i(f)\co \H^i(A) \to \H^i(B)$ is an isomorphism in $\Rep(R)$ for all $i$;

\item a fibration if $f^i\co A^i \to B^i$ is a surjection for all $i$;

\item a cofibration if it has LLP with respect to all trivial fibrations.
\end{enumerate}
\end{lemma}
\begin{proof}
This is standard (see e.g. \cite{schematic}v1 Proposition 4.1).
\end{proof}

\begin{definition}\label{dgaff}
Define $dg\Aff(R)$ to be the category opposite to $DG\Alg(R)$, equipped with the opposite model structure. 
\end{definition}

\begin{definition}
Let $DG\Alg(R)_0$ be the full subcategory of $DG\Alg(R)$ whose objects $A$ satisfy $\H^0(A)=k$. Let $\Ho(DG\Alg(R))_0$ be the full subcategory of $\Ho(DG\Alg(R))$ on the objects of  $DG\Alg(R)_0$. Let $dg\Aff(R)_0$ and $\Ho(dg\Aff(R))_0$ be the opposite categories to $DG\Alg(R)_0$ and $\Ho(DG\Alg(R))_0$, respectively.
\end{definition}

\begin{definition}\label{DKD}
Recall that the Dold-Kan correspondence gives a denormalisation functor $D$ from cochain complexes to cosimplicial complexes by setting 
$$
D^n(V)=\bigoplus_{\begin{smallmatrix} m+s=n \\ 1 \le j_1 < \ldots < j_s \le n \end{smallmatrix}} \pd^{j_s}\ldots\pd^{j_1}V^m,
$$
where we define the $\pd^j$ and $\sigma^i$ using the simplicial identities, subject to the conditions that  for all $v \in V^n$,\, $dv = \sum_{i=0}^{n+1}(-1)^i \pd^i v$ and  $\sigma^i v =0$.
This functor is quasi-inverse to the normalisation functor
$$
N^n(V):= \{v \in V^n \,:\, \sigma^iv =0 \quad \forall i\},\quad d = \sum_{i=0}^{n+1}(-1)^i \pd^i.
$$
\end{definition}

\begin{definition}
Given a bicosimplicial complex $V$, the Eilenberg-MacLane shuffle product $\nabla\co N(\diag V) \to \Tot (NV) $ is a quasi-isomorphism of cochain complexes given by summing:
$$
\nabla^{pq}=\sum_{ (\mu, \nu) \in \mathrm{Sh}(p,q)}(-1)^{(\mu,\nu)} \sigma^{\nu_1}_h\ldots \sigma^{\nu_q}_h\sigma_v^{\mu_1}\ldots \sigma_v^{\mu_p}\co  N(V^{p+q,p+q}) \to (NV)^{pq}
$$
This is associative and commutative.

The Dold-Kan correspondence then gives, for any bi-cochain complex $W$, a quasi-isomorphism
$$
\nabla:=D\nabla N\co \diag(DW) \to D(\Tot W)
$$
of cosimplicial complexes.
\end{definition}

\begin{definition}
The denormalisation functor gives rise to a functor $D\co DG\Alg(R) \to c\Alg(R)$, with the product  $\nabla\co N^n(\diag(DA\ten DA)) \to A^n$  given  by
$$
\nabla(a\ten b)=\sum_{\begin{smallmatrix} p+q=n \\ (\mu, \nu) \in \mathrm{Sh}(p,q)\end{smallmatrix}}(-1)^{(\mu,\nu)} (\sigma^{\nu_1}\ldots \sigma^{\nu_q}a)(\sigma_v^{\mu_1}\ldots \sigma_v^{\mu_p}b).
$$
This extends to a product $D(A)\ten D(A) \to D(A)$ by $\nabla(\pd^ia\ten \pd^ib)=\pd^i\nabla(a\ten b)$.
\end{definition}

\begin{proposition}\label{affequiv}
The functor $D\co DG\Alg(R) \to c\Alg(R)$ is a right Quillen equivalence. Thus $\Ho(c\Alg(R))\simeq \Ho(DG\Alg(R))$, and $\Ho(c\Alg(R))_0\simeq \Ho(DG\Alg(R))_0$.
\end{proposition}
\begin{proof}
This is standard (as in \cite{schematic}v1 Proposition 4.1). The left adjoint $D^*$ is given by Thom-Sullivan cochains.
\end{proof}

\subsection{The Maurer-Cartan functor for cochains}

\begin{definition}\label{mcDG}
Given a cochain algebra $A \in DG\Alg(R)$, and a chain Lie algebra $\g \in dg\hat{\cN}(R)$, define the Maurer-Cartan space by 
$$
\mc(A,\g):=\{\omega \in \prod_n A^{n+1}\hat{\ten}^R \g_n \,|\,d\omega+\half[\omega,\omega]=0\}.
$$
\end{definition}
\begin{remark}
The Maurer-Cartan functor is essentially the same as the functor of twisting functions from the chain coalgebra $A^{\vee}$ to $\g$, as defined in \cite{QRat} \S B.5. It is, of course, also the same as the classical Maurer-Cartan functor on the DGLA $A\hat{\ten}^R\g$.  
\end{remark}

\begin{definition}
Define functors $\xymatrix@1{ dg\Aff(R) \ar@<1ex>[r]^G & dg\hat{\cN}(R) \ar@<1ex>[l]^{\bar{W}}}_{\bot}$ as follows.
For $\g \in dg\hat{\cN}(R)$, the Lie bracket gives a linear map $\bigwedge^2\g \to \g$. Write $\Delta$ for the dual $\Delta\co \g^{\vee} \to \bigwedge^2\g^{\vee}$.
This is equivalent to a  map $\Delta\co \g^{\vee}[-1] \to \Symm^2(\g^{\vee}[-1])$, and we define 
$$
O(\bar{W}\g):= \Symm(\g^{\vee}[-1])
$$
to be the graded polynomial ring on generators $\g^{\vee}[-1]$, with a derivation defined on generators by $D:=d +\Delta$. The Jacobi identities ensure that $D^2=0$. Note that
$$
\Hom_{dg\Aff(R)}(\Spec A, \bar{W}\g) \cong \mc(A,\g).
$$  

We define $G$ by writing $\sigma A^{\vee}[1]$ for the brutal truncation (in non-negative degrees) of $A^{\vee}[1]$, and setting
$$
G(A)= \Lie(\sigma A^{\vee}[1]),
$$
the free graded Lie algebra, with differential similarly defined on generators by $D:=d +\Delta$, $\Delta$ here being the coproduct on $A^{\vee}$. Note that 
$$
\Hom_{dg\hat{\cN}(R)}(G(A), \g) \cong \mc(A,\g),
$$  
so $(G,\bar{W})$ form an adjoint pair. Note also that $G(A)$ is cofibrant for all $A$.
\end{definition}

\begin{remark}\label{koszul}
The adjoint pair $G\dashv \bar{W}$ corresponds to the pair $\cL\dashv \C$ given in \cite{QRat}. Chain Lie algebras are dual to strong homotopy commutative algebras (SHAs), so a  generating space for a cofibrant object of $dg\hat{\cN}$ is dual to an $E_{\infty}$-algebra.  $G$ is then the the canonical functor from DGAs to SHAs, by regarding a  DGA as an $E_{\infty}$-algebra.
This construction is precisely the functor given in $\cite{Kon}$ for comparing DGLAs and SHLAs (i.e. duals of graded pro-Artinian algebras) in algebro-geometric deformation theory, making it a form of Koszul duality. The only difference between  that setting and infinitesimal topology is that for the former the algebras are Artinian, whereas in the latter  the Lie algebras are. 
\end{remark}

\begin{definition}\label{dgdef}\label{dgdefgauge}
Given $A \in DG\Alg(R)$ and $\g \in dg\hat{\cN}(R)$, we define the gauge group by
$$
\Gg(A,\g):= \exp(\prod_n A^n\hat{\ten}^R\g_n).
$$ 
Define a gauge action of $\Gg(A,\g)$ on $\mc(A,\g)$ by 
$$
g(\omega):= g\cdot \omega \cdot g^{-1} -(dg)\cdot g^{-1}.
$$
Here, $a\cdot b$ denotes multiplication in the universal enveloping algebra $\cU(A\hat{\ten}^R\g)$ of the DGLA $A\hat{\ten}^R\g$. That $g(\omega) \in \mc(A,\g)$ is a standard calculation (see \cite{Kon} or \cite{Man}).

We define the
torsor space by
$$
\pi(A,\g)=\mc(A,\g)/\Gg(A,\g).
$$
\end{definition}

\begin{lemma}\label{defworks}
Given $A \in DG\Alg(R)_0$ and $\g \in dg\hat{\cN}(R)$, 
$$
\Hom_{\Ho(dg\Aff(R))}(\Spec A, \bar{W}\g)= \pi(A,\g).
$$
\end{lemma}
\begin{proof}
The gauge functor is represented in $dg\Aff(R)$ by $V\g:=\Spec( k[ \g^{\vee}t \oplus \g^{\vee}dt])$, for $t$ of degree $0$. Observe that $\bar{W}\g$ is fibrant and that $V\g \to\Spec k$ is a trivial fibration (LLP can be established using pro-nilpotence of the augmentation ideals of $O(V\g) \to k$ and $O(\bar{W}\g)\to k$). The gauge action gives us a map $\phi\co  V\g \by \bar{W}\g \to \bar{W}\g$.

We now make $V\g \by\bar{W}\g$ into   a path object  for $\bar{W}\g$ in  $dg\Aff(R)$  via the maps
$$
\xymatrix@1{ \bar{W}\g \ar[r]^-{(\id,1)} &\bar{W}\g \by H  \ar@<0.5ex>[r]^-{\pr_1} \ar@<-0.5ex>[r]_-{\phi} & \bar{W}\g,}
$$
again using pro-nilpotence of the augmentation ideals of $\bar{W}\g,V\g$ to see that $(\pr_1,\phi)$ is a fibration.
\end{proof}

\begin{definition}
Given $A \in DG\Alg(R)$, we may  define $\bar{G}(A)\in dg\hat{\cN}(R)$ to be the  hull of $\pi(A, -)$.
\end{definition}

\begin{remark}
For an explicit description of $\bar{G}(A)$, take a decomposition $A^1=dA^0\oplus U$ as $R$-representations (possible since $R$ is reductive). Then
$$
\bar{G}(A)= G(A)/\langle U^{\perp} \rangle, ,
$$
noting that $U^{\perp} \cong (dA^0)^{\vee}$.
\end{remark}

\begin{proposition}\label{dgqwell}
If $A \in dg\Alg(R)_0$,  then for all $\g \in s\hat{\cN}(R)$, the map
$$
\Hom_{dg\cM(R)}(\bar{G}(A), \g)  \to \pi(A,\g)
$$
is an isomorphism, so $\bar{G}$ defines a functor from $\Ho(dg\Aff(R)_0)$ to $dg\cM(R)$.
\end{proposition}
\begin{proof}
As for Proposition \ref{qwell}.
\end{proof}

The following result can be thought of as a converse to Theorem \ref{spectralpi}.
\begin{proposition}\label{spectralh}
Given $A \in DG\Alg(R)$, there is a convergent Adams spectral sequence (in $\widehat{\FD\Vect}$)
$$
E^1_{pq}= (\Lie_{-p} (\tilde{\H}^{*+1}(A)^{\vee}))_{p+q} \abuts \H_{p+q}(\bar{G}(A)),
$$
where $\tilde{\H}$ denotes reduced cohomology.

Moreover, the differential
$$
d^1_{-1,q}\co \tilde{\H}^q(A)^{\vee} \to ({\bigwedge}^2 \tilde{\H}^{*+1}(A)^{\vee})_{q-2}= ((\Symm^2\tilde{\H}^*(A))^q)^{\vee} 
$$
is dual to the cup product on $\tilde{\H}^*(A)$.
\end{proposition}
\begin{proof}
Define a filtration $F$ on $\bar{G}(A)$ using the lower central series, so $F_{-p}\bar{G}(A)=\Gamma_p\bar{G}(A)$. The proof is now as for Proposition \ref{spectralh1}, with the description of $d^1_{-1,q}$ an immediate consequence of the definition of $\bar{G}$. 
\end{proof}

\begin{remark}\label{toenadams}
When $A$ is the homotopy type of a reduced topological space, Corollary \ref{eqtoen} implies this spectral sequence is equivalent to the weight spectral sequence of \cite{KTP} \S 3, which is only shown to converge for complex projective varieties (assuming that all the necessary results from \cite{schematic}v1 have analogues  in [ibid.]v2). The weight spectral sequence was defined in the category of  commutative unipotent affine group schemes; by Tannakian duality, this is equivalent to $\widehat{\FD\Vect}$. Thus it seems that  Remark \ref{FDVspectral} might have been  the only new ingredient needed to prove convergence in general. However, whereas the  weight spectral sequence is an $E^1$ spectral sequence,  the proof of Proposition \ref{spectralh} gives an $E^0$ spectral sequence --- this additional structure  is significant for many computations. 
\end{remark}

\subsection{Chain versus simplicial}

\begin{theorem}\label{qs}
Given a cochain algebra $A\in DG\Alg(R)$, and $\g \in s\hat{\cN}(R)$, there is a canonical isomorphism
$$
N\co \mc(DA, \g) \to \mc(A, N(\g)),
$$
functorial in $A$ and $\g$. This gives us a  diagram
$$
\xymatrix{
dg\Aff(R) \ar@<1ex>[r]^{\Spec D} \ar@<-1ex>[d]_G& \ar@<1ex>[l]^{\Spec D^*} s\Aff(R)\ar@<1ex>[d]^G\\
dg\hat{\cN}(R)\ar@<-1ex>[r]_{N^*}\ar@<-1ex>[u]_{\bar{W}} &s\hat{\cN}(R)\ar@<1ex>[u]^{\bar{W}}\ar@<-1ex>[l]_{N}
}
$$
in which both the outer and inner squares commute, where $D^*,N^*$ are the left adjoints to $D,N$.
\end{theorem}

\begin{proof}
We begin by constructing the map $\mc(DA, \g) \to \mc(A, N(\g))$. Let $G=\exp(\g)$, so that $O(G)$ is dual to the universal enveloping algebra $\cU(\g)$. Let $\vareps\co O(G)\to k$ be the co-unit, with augmentation ideal $I$. The multiplication on $O(G)$ gives a coproduct $\rho$ on $I^{\vee}$, with 
$$
\g=\ker( \rho\co  I^{\vee} \to I^{\vee}\hat{\ten} I^{\vee} ),
$$
so the Dold-Kan equivalence gives
$$
\bar{N}\g=\ker(\bar{N}\rho\co  \bar{N}I^{\vee} \to \bar{N}(I^{\vee}\hat{\ten} I^{\vee} )),
$$
for $\bar{N}$ as in Lemma \ref{barng}.

Given $\omega \in \mc(DA,\g)$, we can then regard $\omega_n$ as a ring homomorphism
$$
\omega_n\co O(G_n)\to D^{n+1}(A),
$$
 giving
$$
N(\omega_n)\co N^nI\to  A^{n+1},
$$
using the Maurer-Cartan identities for $\sigma^i$.

Now, consider the composition
$$
N^n(I\ten I) \to N^nI \xra{N(\omega_n)} A^{n+1}.
$$
The Maurer-Cartan identities for the $\sigma^i$ show that
$$
\sigma^{\mu_1}\ldots\sigma^{\mu_p}\omega_n(x)=\omega_{n-p}(\sigma^{\mu_1-1}\ldots\sigma^{\mu_p-1}x),$$
where we set $\sigma^{-1}:=\vareps=0$ on $I$. Since any $(p,q)$-shuffle must have either $\mu_1=0$ or $\nu_1=0$, this means that $N(\omega_n)(xy)=0$. 
Thus 
$$
N(\omega_n) \in \ker \bar{N}(\rho_n)=A^{n+1}\hat{\ten}^R\bar{N}_n(\g).
$$

\begin{lemma}
$N(\omega) \in \mc(A,\bar{N}(\g))$.
\end{lemma}
\begin{proof}[Proof of lemma]
For $x \in N^nI$, we have
\begin{eqnarray*}
dN(\omega)(x)&=& \sum_{i=0}^{n+2}(-1)^i\pd^i\omega_n(x) +\sum_{i=0}^{n+1} (-1)^i\omega_{n+1}(\pd^ix)\\
&=& \pd^0\omega_n(x)-\pd^1\omega_n(x) + (\pd^1\omega_n\ten\pd^0\omega_n )(x\ten 1 -1\ten x -\psi(x))\\
&=& -(\pd^1\omega_n\nabla\pd^0\omega_n )\psi(x) 
\end{eqnarray*}
where $\psi\co I \to I\ten I$ is defined by saying that
\begin{eqnarray*}
I &\to& O(G)\ten O(G)\\
x &\mapsto& x\ten 1  -1\ten x -\psi(x)
\end{eqnarray*}
corresponds to 
\begin{eqnarray*}
G\by G &\to&G\\
(g,h) &\mapsto& gh^{-1}.
\end{eqnarray*}

Now, for $x\ten y \in N^{n-1}(I\ten I)$, consider
$$
\pd^1\omega_{n-1}(x)\ten \pd^0\omega_{n-1}(y) \in D^{n+1}(A)\ten D^{n+1}(A).
$$
We may write this as 
$$
\pd^1\omega_{n-1}(x)\ten \pd^0\omega_{n-1}(y)= ((\pd^1-\pd^0)\omega_{n-1}(x))\ten \pd^0\omega_{n-1}(y)+\pd^0(\omega_{n-1}(x)\ten \omega_{n-1}(y)),
$$
noting that 
\begin{eqnarray*}
((\pd^1-\pd^0)\omega_{n-1}(x))\ten \pd^0\omega_{n-1}(y)&\in& N^{n+1}(DA\ten DA)\\
\omega_{n-1}(x)\ten \omega_{n-1}(y) &\in& N^{n}(DA\ten DA),
\end{eqnarray*}
thus allowing us to calculate $\nabla$ using shuffles. 

Now, $\omega_{n-1}(x)\nabla \omega_{n-1}(y)=0$, leaving us to consider terms of the form
$$
(\sigma^{\nu_1}\ldots \sigma^{\nu_q}(\pd^1-\pd^0)\omega_{n-1}(x))( \sigma^{\mu_1}\ldots\sigma^{mu_p}\pd^0\omega_{n-1}(y)),
$$
for  shuffle permutations $(\mu,\nu)$. Since $\sigma^0\omega_{n-1}x=0$, this expression will be $0$ if $\nu_1=0$, so we must have $\mu_1=0$. Similarly, it will then be $0$ if $\mu_2=1$, so we must have $\nu_1=0$. Setting $p':=p-1, q':=q-1,\nu'_i:=\nu_{i+1}-2,\mu'_i:=\mu_{i+1}-2$, we get
\begin{eqnarray*}
(\pd^1\omega_{n-1}x)\nabla(\pd^0\omega_{n-1}y)&=& \!\!\!\!\sum_{\begin{smallmatrix} p'+q'=n-1\\ (\mu', \nu') \in \mathrm{Sh}(p,q)\end{smallmatrix}}\!\!\!\!\!\!\!\!(-1)^{(\mu',\nu')} \omega_{p'}(\sigma^{\nu'_1}\ldots \sigma^{\nu'_{q'}}x) \cdot \omega_{q'}(\sigma^{\mu'_1}\ldots \sigma^{\mu'_{p'}}y),\\
&=& \sum_{p'+q'=n-1}(N(\omega_{p'})\ten N(\omega_{q'}))(x\nabla^{p'q'} y). 
\end{eqnarray*}

Therefore $dN(\omega)_n $ is equal to the composition
$$
N^{n-1}I \xra{-\psi} N^{n-1}(I\ten I) \xra{\nabla} \bigoplus_{p'+q'=n-1} N^{p'}I\ten N^{q'}I \xra{N(\omega_{p'})\ten N(\omega_{q'})} A^{n+1}.
$$

Now, since $N(\omega)$ is $0$ on $I\cdot I$, we may replace $I$ by $I/(I\cdot I)=\g^{\vee}$. Then $\psi$ is just dual to $\half[,]$, by the Campbell-Baker-Hausdorff formula. But $\nabla\psi$ is therefore dual to $\half\llbracket,\rrbracket$, so we have 
$$
dN(\omega) +\half \llbracket \omega,\omega\rrbracket=0,
$$
giving $N(\omega) \in \mc(A,N(\g))$, as required.  
\end{proof}

It is a straightforward exercise in the Dold-Kan correspondence to show that if $\g$ is abelian, then $\omega \mapsto N(\omega)$ is an isomorphism of Maurer-Cartan spaces.

Now, for fixed $A$, $\mc(A, N(\g))$ is represented on $s\hat{\cN}(R)$ by $N^*G(A)$, which is cofibrant, since $N^*$ is left Quillen and $G(A)$ cofibrant. $\mc(A, N(\g))$ is represented by $GD(A)$, which is similarly cofibrant. Thus $\theta$ gives rise to a natural map $\theta^*\co N^*G(A) \to GD(A)$ which is an isomorphism on abelianisations (considering $\g$ abelian). By considering the graded Lie algebras associated to the lower central series, we see that $\theta^*$ is an isomorphism, as required.
\end{proof}

\begin{corollary}\label{bigequiv}\label{qqs}
We have the following commutative diagram of equivalences of categories:
$$
\xymatrix{
\Ho(dg\Aff(R))_0 \ar@<1ex>[r]^{\Spec D} \ar@<1ex>[d]^{\bar{G}}& Ho(s\Aff(R))_0 \ar@<1ex>[l]^{\R(\Spec D^*)} \ar@<1ex>[d]^{\bar{G}} \\
dg\cM(R) \ar@<1ex>[u]^{\bar{W}} \ar@<1ex>[r]^{\bL N^*} & s\cM(R), \ar@<1ex>[u]^{\bar{W}} \ar@<1ex>[l]^{N}
}
$$
with $\bL N^*\bar{G}=N^*\bar{G}$, since everything in the image of $\bar{G}$ is cofibrant.
\end{corollary}
\begin{proof}
By Theorem \ref{qs}, we know that $\R(\Spec D^*)\bar{W}=\bar{W}N$. Propositions \ref{qwell} and  \ref{dgqwell} then show that $\bL N^*\bar{G}=\bar{G}\Spec D$. Finally, Propositions \ref{wequiv}, \ref{affequiv} and \ref{nequiv} show that three sides of the square are pairs of equivalences, so the fourth must also be.\end{proof}

\begin{remark}\label{baues}
Taking $R=1$, this allows us to make a direct comparison between Sullivan's and Quillen's rational homotopy types. Sullivan's  is $\bL D^*\CC^{\bt}(X,k)$, while  Quillen's homotopy type is a chain coalgebra $C$ whose dual algebra is given by  $\Spec C^{\vee}= \bar{W}N\bar{G}\CC^{\bt}(X,k)$. This proves a generalisation of the Baues-Lemaire conjecture (\cite{baueslemaire} Conjecture 3.5), which was first proved by Majewski (\cite{majewski}).
\end{remark}

\begin{remark}\label{deligne}
There is a notion of minimal models in $DG\Alg(R)_0$ extending that given for cochain algebras in \cite{Sullivan}. Explicitly, $M$ is minimal if it is cofibrant, and there are subrepresentations $V^p \le A^p$ for $p\ge 1$, freely generating $A$ as a graded algebra, with $dV \subset (V\cdot V)$.
 The homotopy groups can easily be recovered from $M$ by using the fact that $M$ is cofibrant, as
\begin{eqnarray*}
\pi_{n-1}(\bar{G}(M))&\cong&\Hom_{dg\cM(R)}(\bar{G}(k\oplus O(R)^{\vee}[-n]\eps), \bar{G}(M))\\
&\cong&\Hom_{\Ho(DG\Alg(R))}(M, k\oplus O(R)^{\vee}[-n]\eps)\\
&\cong& (V^n)^{\vee},
\end{eqnarray*}
for $\eps^2=0$. This recovers Deligne's original conception of schematic homotopy groups.
However, the  minimal models in $dg\hat{\cN}(R)$ seem to be more convenient in applications, since their generators correspond to homology rather than homotopy.

Although \cite{schematic}v1 Corollary 3.22 showed that schematic homotopy groups can be recovered from the minimal DGA, an explicit formula was not proved; the missing ingredient was a simple description of $\bL\cO_R(S^n\by BR)$ (their notation), which we have shown is just $D(k\oplus O(R)^{\vee}[-n]\eps)$.
\end{remark}

\begin{theorem}\label{defqs}
Given a cochain algebra $A\in DG\Alg(R)$, and $\g \in s\hat{\cN}(R)$, there is a canonical isomorphism
$$
N\co \Gg(DA, \g) \to \Gg(A, N(\g)),
$$
functorial in $A$ and $\g$, which is  compatible with $N\co \mc(DA,\g) \to \mc(A,N(\g))$, in the sense that
$$
N(g)(N(\omega))=N(g(\omega)),
$$
for
the respective gauge actions. Thus there is a canonical isomorphism $N\co \pi(DA,\g) \to \pi(A, N(\g))$ of torsor spaces.
\end{theorem}
\begin{proof}
Given  $g \in \Gg(A,\g)$, define $N \log g$ to be the projection of $\log g$ onto 
$$
 \prod_n D^nA\hat{\ten}^R\bar{N}_n\g.
$$
Since $\sigma^ig_n=\sigma_ig_{n+1}$, 
$$
N \log g \in  \prod_n N^n(DA)\hat{\ten}^R\bar{N}_n\g=\prod_n A^n\hat{\ten}^R\bar{N}_n\g,
$$
and both Lie algebra structures (using simplicial or cosimplicial shuffles) agree.
Then define $Ng=\exp(N\log g) \in \exp( \prod_n A^n\hat{\ten}^R\bar{N}_n\g)$. 
 That this is an isomorphism between $\Gg(DA, \g)$ and $\Gg(A, N(\g))$ follows from the proof of the dual Dold-Kan correspondence.

Now, take  $\omega \in \mc(DA,\g)$, giving
$$
N(g(\omega))_n \equiv (\pd_0g_{n+1}\nabla \pd^0g_n^{-1} -1) +  \pd_0g_{n+1}\nabla (\omega_n-1) \nabla\pd^0g_{n}^{-1}
$$

Since $\pd_0g\ten (\omega_n-1)\ten \pd^0g^{-1} \in N(DA^{\ten 3})\hat{\ten}^R \bar{N}(g)^{\hat{\ten} 3}$, we may compute the second term in terms of shuffles, giving
$$
\sum(-1)^{(\nu,\mu,\lambda)} \sigma^{\nu_1}\ldots \sigma^{\nu_r}\pd_0g_{n+1} \cdot \sigma^{\mu_1}\ldots \sigma^{\mu_q}(\omega_n-1)\cdot \sigma^{\lambda_1}\ldots \sigma^{\lambda_p}\pd^0g_n^{-1},
$$
with $p+q+r=2(n+1)$. Since we must have two of $\nu_1,\mu_1,\lambda_1$ being $0$, and since $\sigma^0\omega_n=1$, only shuffles with $\nu_1=\lambda_1=0$ contribute. The gauge and Maurer-Cartan identities make this equal
$$
\sum (-1)^{(\nu,\mu,\lambda)}\sigma_{\nu_r}\ldots \sigma_{\nu_1}g_{n-r} \cdot \sigma_{\mu_q}\ldots \sigma_{\mu_1}(\omega_{n-q}-1)\cdot \sigma_{\lambda_p}\ldots \sigma_{\lambda_1}g_{n-p}^{-1},
$$
summing over shuffles with $p+q+r=2n$. This is $N(g)N(\omega)N(g)^{-1}$.

For the first term, we have
\begin{eqnarray*}
N(g(1))&=& (\pd_0g_{n+1}\nabla \pd^0g_n^{-1} -1)\\
&=& \pd^0( g_{n}\nabla g_n^{-1} -1) + (\pd_0g_{n+1} -\pd^0g_n) \nabla \pd^0g_n^{-1}.
\end{eqnarray*}
The first term is $0$, and $(\pd_0g_{n+1} -\pd^0g_n) \nabla \pd^0g_n^{-1} \in N(DA^{\ten 2})\hat{\ten}^R \bar{N}(g)^{\hat{\ten} 2}$, so the second term may be computed in terms of shuffles, giving 
\begin{eqnarray*}
N(g(1)) &=& \sum_{p+q=n} (-1)^{(\nu,\mu)}\sigma^{\nu_1}\ldots \sigma^{\nu_q}(\pd_0g_{n+1} -\pd^0g_n)\cdot  \sigma^{\mu_1}\ldots \sigma^{\mu_p}\pd^0g_n^{-1}\\
&=& \sum_{p+q=n-1}(-1)^{(\nu,\mu)} \sigma_{\nu_q}\ldots \sigma_{\nu_1}(\pd_0g_{n+1-q} -\pd^0g_{n-q})\cdot  \sigma_{\mu_p}\ldots \sigma_{\mu_1}g_{n-q}^{-1},
\end{eqnarray*}
using the gauge identities. 

Now,
$$
(dg)_r = \sum_{i=0}^{r+1}(-1)^r\pd^ig_r -\sum_{i=0}^{r+1} (-1)^i\pd_ig_{r+1}= \pd^0g_r-\pd_0g_{r+1},
$$
all the other terms cancelling, so
$$
N(g(\omega))=N(g)(N(\omega)),
$$
as required.
\end{proof}

\subsection{Differential forms}

Take a topological space $X$ and a reductive pro-algebraic groupoid $R$, equipped with a morphism 
$$
\rho\co X \to |BR(k)|,
$$
such that $\pi_f\rho\co \pi_fX \to R(k)$ is essentially surjective on objects and Zariski-dense on morphisms. 

\begin{definition}
Recall that  $O(R)$ has the natural structure of an $R\by R$-representation. Since  every $R$-representation has an associated semisimple local system on $|BR(k)|$, we will also write $O(R)$ for the $R$-representation in semisimple local systems on $|BR(k)|$  corresponding to the $R\by R$-representation $O(R)$.
We then define the $R$-representation $\bO(R)$ in semisimple local systems on $X$ by $\bO(R):=\rho^{-1}O(R)$.
\end{definition}

\begin{definition}
For any sheaf $\sF$ on $X$, define 
$$
\CC^n(X,\sF):= \prod_{f\co \Delta^n \to X} \Gamma(\Delta^n, f^{-1}\sF).
$$
Together, these form a cosimplicial complex $\CC^{\bt}(X,\sF)$.
\end{definition}

\begin{lemma}
There is a canonical isomorphism
$$
\CC^{\bt}(X,\bO(R)) \cong \CC^{\bt}(\Sing X, \bO(R))
$$
in $c\Alg(R)$, where the latter is defined as in Definition \ref{bo}.
\end{lemma}

\begin{remark}\label{irreds}
$\FD\Rep(R)$ is equivalent to a tensor subcategory of the category of semisimple local systems on $X$. If $\{V_{\alpha}\,:\,\alpha \in I\}$ is s set of representatives of irreducible $R$-representations, let $\{\vv_{\alpha}\,:\,\alpha \in I\}$ be the corresponding semisimple local systems on $X$. Then Tannakian duality gives us an isomorphism
$$
O(R)\cong \bigoplus_{\alpha \in I} \Hom_{\End_R(V_{\alpha})}(V_{\alpha}\boxtimes k ,k\boxtimes V_{\alpha})
$$
as $R \by R$-representations. We then have an isomorphism
$$
\bO(R)\cong \bigoplus_{\alpha \in I} \Hom_{\End_R(V_{\alpha})}(V_{\alpha} ,\vv_{\alpha})
$$
as an $R$-representation in local systems on $X$. Note that $\End_R(V_{\alpha})$ is a division algebra of finite dimension over $k$, so must equal $k$ if $k$ is algebraically closed. 
\end{remark}

\begin{definition}
Given a manifold $X$, denote the sheaf  of real $n$-forms on $X$ by $\sA^n$. Given a real sheaf $\sF$ on $X$, write
$$
A^n(X,\sF):=\Gamma(X,\sF\ten_{\R} \sA^n).
$$ 
\end{definition}

\begin{proposition}\label{propforms}
The real Malcev  homotopy type of a manifold $X$ relative to $\rho\co \pi_fX \to R(\R)$ is given in $DG\Alg(R)$ by 
$$
A^{\bt}(X, \bO(R)).
$$
\end{proposition}
\begin{proof}
Consider the following morphisms in $c\Alg(R)$:
$$
\CC^{\bt}(X,\bO(R))\to \diag \CC^{\bt}(X,\bO(R)\ten_{\R}D\sA^{\bt}) \la \Gamma(X,\bO(R)\ten_{\R}D\sA^{\bt})=DA^{\bt}(X, \bO(R)).
$$
Taking cohomology, we have
$$
\H^*(X,\bO(R)) \to \bH^*(X, \bO(R)\ten_{\R}\sA^{\bt}) \la \H^*(A^{\bt}(\bO(R))).
$$
The first arrow is an isomorphism as $\R \to \sA^{\bt}$ is a resolution. The second is  an isomorphism as $\sA^n$ is flabby. Therefore $\CC^{\bt}(X,\bO(R))$ is weakly equivalent to $DA^{\bt}(X, \bO(R))$, as required.
\end{proof}

\begin{remark}
In the notation of Remark \ref{irreds}, 
$$
A^{\bt}(X, \bO(R)) \cong \bigoplus_{\alpha \in I} \Hom_{\End_R(V_{\alpha})}(V_{\alpha} ,A^{\bt}(X,\vv_{\alpha})).
$$
\end{remark}

\section{Automorphisms and formality}\label{hodge}

\begin{definition}
Given a reductive pro-algebraic groupoid $R$, we say that $A \in DG\Alg(R)$ is formal if it is weakly equivalent to its cohomology algebra $\H^*(A)$.
\end{definition}

\begin{definition}
We say that the Malcev homotopy type $(X,\rho)^{\mal}$ of a topological space $X$  relative to a  Zariski-dense homomorphism $\rho\co \pi_fX \to R$ is formal if $\bL D^*\CC^{\bt}(X, \bO(R))$ is formal. We say that the pro-algebraic homotopy type of a topological space $X$ is formal if $(X,\rho)^{\mal}$ is formal for $\rho\co \pi_fX \to (\pi_fX)^{\red}$.
\end{definition}

\begin{proposition}\label{formalpins}
If  $(X,\rho)^{\mal}$ is formal, then $N\Ru G(X,\rho)^{\mal}\in dg\hat{\cN}(R)$ is weakly equivalent to the free chain Lie algebra $\bar{G}(\H^*(X,\bO(R)))$ generated in degree $n$ by  $\H^{n+1}(X, \bO(R))^{\vee}$, with differential defined on generators by the coproduct
$$
\Delta\co  \H^{n+1}(X, \bO(R))^{\vee} \to \prod_{a+b=n-1} \H^{a+1}(X, \bO(R))^{\vee}\hat{\ten} \H^{b+1}(X, \bO(R))^{\vee},
$$
so
$$
\varpi_n(X^{\rho,\mal})= \H_{n-1}\bar{G}(\H^*(X,\bO(R)))
$$
for $n>1$, and
$$
\Ru\varpi_1(X^{\rho,\mal})=\H_{0}\bar{G}(\H^*(X,\bO(R))).
$$
\end{proposition}
\begin{proof}
Since $\bL D^*\CC^{\bt}(X, \bO(R))$ is formal, it is weakly equivalent to $\H^*(X, \bO(R))$. Now apply corollary \ref{qqs}, observing that the chain Lie algebra described above is just $\bar{G}(\H^*(X,\bO(R)))$.
\end{proof}

\begin{remark}
 If the pro-algebraic homotopy type is formal, we can therefore describe the pro-algebraic homotopy groups explicitly in terms of the cohomology ring.  This implies the  result in \cite{higgs} that  formality ensures the unipotent radical of the pro-algebraic fundamental group is quadratically presented. 
\end{remark}

\begin{theorem}\label{kformal}
The real pro-algebraic homotopy type of a compact K\"ahler manifold $X$ is formal.
\end{theorem}
\begin{proof}
By Proposition \ref{propforms}, the real homotopy type of $X$ relative to $R:=(\varpi_fX)^{\red}$ is given by the $R$-cochain algebra $A^{\bt}(X,\bO(R))$ of differential forms. As in \cite{higgs} Theorem \ref{higgs-formal} (or equivalently \cite{KTP} Corollary 2.1.3), consider the following morphisms of cochain algebras: 
$$
A^{\bt}(X,\bO(R)) \leftarrow (\z_{\dc}(A^*(X,\bO(R))),d) \to (\H_{\dc}^*(X,\bO(R)),d),
$$
where $\dc:=iD'_K-iD''_K$ is defined in terms of the harmonic metric $K$.  For any point $x \in X$, $\bO(R)(x)$ corresponds to a representation of $(\varpi_1(X,x))^{\red}$, so is a directed sum of finite-dimensional semisimple local systems on $X$. The principle of two types (\cite{Simpson} Lemmas 2.1 and 2.2, using the $d\dc$ lemma instead of the $\pd\bar{\pd}$ lemma) thus shows  that these morphisms are quasi-isomorphisms, and that $d=0$ on $\H_{\dc}(A(V))$.
\end{proof}

\subsection{Automorphisms and weights}

\begin{definition}
Given a pro-algebraic simplicial group $G$, define $c\Rep(G)$ to be the category of cosimplicial $O(G)$-comodules, so that dualisation gives a contravariant equivalence between $c\Rep(G)$ and $s\widehat{\FD\Rep}(G)$. Given $V \in c\Rep(G)$, define $\H^i(G,V):=\H^i(\R^i\Gamma V)$, where $\R^i\Gamma$ is  the right-derived invariant functor  from  $c\Rep(G)$ to $c\Vect$. The dualisation correspondence then gives an isomorphism
$$
\H^i(G,V)^{\vee} \cong H_i(G,V^{\vee}).
$$
In particular, if $V$ is a semisimple $\pi_fX$-representation with corresponding local system $\vv$, then the results of \S\ref{mc} imply that 
$$
\H^i(G(X)^{\alg},V) \cong \H^i(X,\vv).
$$
\end{definition}

\begin{definition}
Given $\g \in dg\hat{\cN}(R)$, and a chain complex $V \in dg\widehat{\FD\Rep}(R)$ equipped with the structure of a $\g$-module, we now define the chain complex $\Der_{\bt}(\g,V)$ of derivations. For each integer $n$, 
$$
\Der_n(\g,V)  
$$
 consists of all degree $n$ derivations, i.e. maps $\theta_i\co \g \to V_{n+i}$ such that  
$$
 \theta_{i+j}([a,b])= [\theta_i a, b] + (-1)^{in} [a,\theta_j b].
$$
 We then define the differential $D$ on $\Der_{\bt}(\g,V)$ by $D\theta=d\theta -(-1)^{\bar{\theta}}\theta d$. Note that if $V$ is an ideal in $\g$, then $\Der_{\bt}(\g,V)$ has the natural structure of a chain Lie algebra.
\end{definition}

\begin{definition}\label{cna}
Given a $k$-algebra $A$, we say that a Lie coalgebra $C$ over $A$ is ind-conilpotent if it is a union of conilpotent Lie coalgebras. Define $\hat{\cN}_A$ to be opposite to the category of ind-conilpotent Lie coalgebras over $A$.
For a reductive pro-algebraic groupoid $R$ over $k$, let $\hat{\cN}_A(R)$ be the category of $R$-representations in $\hat{\cN}_A$.

Similarly, let $dg\hat{\cN}_A(R)$ be the category of $R$-representations in the opposite category $dg\hat{\cN}_A$ to the category of ind-conilpotent $\N_0$-graded cochain Lie coalgebras over $A$.  Let $s\hat{\cN}_A(R)$ be the category of simplicial objects in $\hat{\cN}_A(R)$.
\end{definition}

Note that $\hat{\cN}_k(R)\cong \hat{\cN}(R)$, and that there is a continuous functor $\hat{\cN}(R)\to \hat{\cN}_A(R)$ given by $C^{\vee} \mapsto (C\ten_k A)^{\vee}$. We denote this by $\g \mapsto \g \hat{\ten} A$.

Recall that there are cofibrantly generated model structures on the categories $c\Mod_A(R)$ and $DG\Mod_A(R)$ of $R$-representations in cosimplicial $A$-modules and $\N_0$-graded cochain $A$-modules, in which fibrations are surjections, and weak equivalences are isomorphisms on cohomology. Dualisation allows us to  regard $\Mod_A(R)^{\opp}$ as the category of $A^{\vee}$-comodules over $\widehat{\FD\Vect}_k$.

\begin{lemma}\label{cnamod}
There is a closed model structure on $dg\hat{\cN}_A(R)$ (resp. $s\hat{\cN}_A(R)$)  in which a morphism $f\co \g \to \fh$ is a fibration or a weak equivalence whenever the underlying map 
$f^{\vee}: \fh^{\vee} \to \g^{\vee}$ in $DG\Mod_A(R)$ 
(resp. $c\Mod_A(R)$)  is a cofibration or a weak equivalence.
\end{lemma}
\begin{proof}[Proof (sketch).] The only non-trivial axioms are factorisation.
Let $F: (DG\Mod_A(R))^{\opp} \to dg\hat{\cN}_A(R)$ be left adjoint to the forgetful functor. Observe that whenever $f:M \to N$ is a fibration (resp. a trivial fibration) in $(DG\Mod_A(R))^{\opp}$, $F(f^{\vee})$ is a  cofibration (resp. a trivial cofibration). The idea behind the proof is to regard these as the generating cofibrations (resp. generating trivial cofibrations). Although the domains of these maps do not satisfy the small object argument in its usual form, nilpotence allows us to simplify the construction. Explicitly, we inductively construct  direct systems $\g^{(m)}$ with the property that $\g^{(n)}/\Gamma_n  \g^{(n)}\to \g^{(m)}/\Gamma_n  \g^{(m)}$ is an isomorphism for all $m \ge n$, and note that $F(M^{\vee})$ is small relative to systems of this form.
\end{proof}

\begin{definition}
Define the category $dg\cM_A(R)$ to have the same objects as $dg\hat{\cN}_A(R)$, with morphisms given by
$$
\Hom_{dg\cM_A(R)}(\g,\fh)=\Hom_{\Ho(dg\hat{\cN}_A(R))}(\g,\fh)/\exp(\fh_0^R),
$$
where $\fh_0^R$ is the sub-Lie algebra of $\fh_0(A)$ fixed by $R$, acting by conjugation on the set of homomorphisms. Define $s\cM_A(R)$ as a quotient of $\Ho(s\hat{\cN}_A(R))$ similarly.
\end{definition}

\begin{proposition}\label{anequiv}
The normalisation functor $N\co s\hat{\cN}_A(R) \to dg\hat{\cN}_A(R)$ is a right Quillen equivalence. Thus $\Ho(s\hat{\cN}_A(R))\simeq \Ho(dg\hat{\cN}_A(R))$, and $s\cM_A(R) \simeq dg\cM_A(R)$.
\end{proposition}
\begin{proof}
As for \cite{QRat} Theorem II.5.4.
\end{proof}

\begin{definition}\label{rout}
Given $\g \in s\hat{\cN}_R$, define a group-valued functor $\Out_R(\g)$  on the category of $k$-algebras by setting
$$
\Out_R(\g)(A):= \Aut_{s\cM_A(R)}(\g\hat{\ten} A).
$$

Given $G \in s\cE(R)$, define $\ROut(G):=\Out_R(\Ru(G))$, noting that
$$
\ROut(G)(k) \cong \Aut_{\Ho(s\cE(R))}(G).
$$
 For $G \in s\agpd$, set $\ROut(G):=\Out_{G^{\red}}(\Ru(G))$.
\end{definition}

\begin{theorem}\label{auto}
If  $G \in s\cE(R)$ is such that $\H^i(G,V)$ is finite-dimensional for all $i$ and all finite-dimensional irreducible $R$-representations $V$, then the group-valued functor
$$
\ROut(G)
$$
is represented by a  pro-algebraic group over $k$. The map
$$
\ROut(G) \to \prod_i  \Aut_R(\H^i(G,O(R)))
$$
of pro-algebraic groups has pro-unipotent kernel.
\end{theorem}
\begin{proof}
As in Proposition \ref{dgminimal}, take a minimal model $\m$ for $N\Ru(G)$. Since $\m$ is cofibrant, the map
$$
\Aut_{dg\hat{\cN}_A(R)}(\m\hat{\ten}A) \to \Aut_{dg\cM_A(R)}(\m\hat{\ten}A)\cong\ROut(G)(A)
$$
is surjective.

The group-valued functor $\Aut_{dg\hat{\cN}(R)}(\m)$, given by
$$
A \mapsto \Aut_{dg\hat{\cN}_A(R)}(\m\hat{\ten}A)
$$ 
is a closed subgroup of the functor $\Aut_{R}(\m)$ given by
$$
A \mapsto \Aut_{R,A}(\m\hat{\ten}A),
$$ 
which is a filtered direct limit of affine schemes, by the definition of a pro-category. In order to show that $\Aut_{dg\hat{\cN}(R)}(\m)$ is a pro-algebraic group, it will suffice to show that it is an affine scheme.

Now let $\m^{\ab}=\m/[\m,\m]$, so that
$$
\m^{ab}_i \cong \H_{i+1}(\Ru(G),k)\cong \H_{i+1}(G, O(R)^{\vee}), 
$$ 
since $R$ is reductive (so
$\H_i(G,V)\cong \H_{i}(\Ru(G),k)\hat{\ten}_R V$). 

Let $I$ be a set of representatives of irreducible $R$-representations. The canonical decomposition for  $\m^{ab}$ is then 
$$
\m^{\ab}_i\cong\prod_{\alpha\in I}(\m^{\ab}_i)^{\alpha},
$$
where each $(\m^{\ab}_i)^{\alpha}$ is a product of copies of the irreducible $V_{\alpha}$. From Remark \ref{irreds} and the finiteness hypotheses on $\H_*(G,V_{\alpha})$, it follows that $(\m^{\ab}_i)^{\alpha}$ is finite-dimensional.

Now, a choice of  Lie algebra generators for $\m$ (isomorphic to $\m^{\ab}$) gives an embedding
$$
\Aut_{dg\hat{\cN}(R)}(\m) \subset \prod_{\alpha, i} \Hom_R((\m^{\ab}_i)^{\alpha} , (\m_i)^{\alpha}),
$$
which is an affine scheme, since $(\m^{\ab}_i)^{\alpha}$ is finite-dimensional. This proves that the group-valued functor $\Aut_{s\hat{\cN}(R)}(\m)$ is  pro-algebraic.

The rest of the proof is an adapted from \cite{Sullivan} \S 6. The lower central series filtration $\Gamma_n \m$ on $\m$ gives a filtration on $\Aut_{dg\hat{\cN}(R)}(\m)$, namely
$$
F^n\Aut_{dg\hat{\cN}(R)}(\m)= \ker(\Aut_{dg\hat{\cN}(R)}(\m) \to \Aut_{dg\hat{\cN}(R)}(\m/\Gamma_{n+1}\m)),
$$
making $F^1\Aut_{dg\hat{\cN}(R)}(\m)$ into a pro-unipotent group. The associated Lie algebra is 
$$
\z_0\Der_{\bt}(\m,[\m,\m])=\ker(D\co \Der_0(\m,[\m,\m]) \to \Der_{-1}(\m,[\m,\m])).
$$

Note that the subspace
$$
 D\Der_1(\m,\m) \le \z_0\Der(\m,\m)
$$
 is an ideal in $\z_0\Der(\m,[\m,\m])$ (since minimality implies $d\m \le [\m,\m]$). Studying the path object in $dg\hat{\cN}_A(R)$ from Lemma \ref{dgpath}, we see that
$$
\ker(\Aut_{dg\hat{\cN}(R)}(\m) \to \Aut_{\Ho(dg\hat{\cN}(R))}(\m))=\exp(D\Der_1(\m,\m)).
$$

Inner derivations give us a map $\ad$ from $\m^R$ onto an ideal of $\z_0\Der_{\bt}(\m,[\m,\m])$, and so
$$
\ker(\Aut_{dg\hat{\cN}(R)}(\m) \to \Aut_{dg\cM(R)}(\m))=\exp(D\Der_1(\m,\m)+ \ad\m^R).
$$
Therefore $\Aut_{dg\cM(R)}(\m)$ is the pro-algebraic group
$$
\Aut_{dg\cM(R)}(\m)=\Aut_{dg\hat{\cN}(R)}(\m)/\exp(D\Der_1(\m,\m)+ \ad\m^R),
$$
with 
$$
\ker(\Aut_{dg\cM(R)}(\m) \to \Aut_{R}(\m^{\ab})) 
$$
a quotient of the pro-unipotent group $F^1\Aut_{dg\hat{\cN}(R)}(\m)$, as required.
\end{proof}

\begin{corollary}
If $X$ is a finite simplicial complex and $\rho\co \pi_fX \to R(k)$ essentially surjective on objects and Zariski-dense on morphisms, then $\Aut_{\Ho(s\cE(R))}(G(X,\rho)^{\mal})$ is  a pro-algebraic group, with 
$$
\ker(\Aut_{\Ho(s\cE(R))}(G(X,\rho)^{\mal}) \to \Aut_R(\H^*(X,\bO(R))))
$$
pro-unipotent.
\end{corollary}

\begin{definition}\label{wgtdef}
Given $G$ as above, define a weight decomposition on  $G \in s\cE(R)$ to be a morphism
$$
\bG_m \to \ROut(G)
$$
of pro-algebraic groups.
\end{definition}

\begin{lemma}\label{wgtpin}
A weight decomposition on $G \in s\cE(R)$ gives rise to weight decompositions on the homotopy groups $\pi_nG$, unique up to conjugation by $\Ru(\pi_0G)$.
\end{lemma}
\begin{proof}
Observe that there is a canonical map 
$$
\ROut(G) \to  \Aut_R(\pi_nG)/ \Ru(G_0) = \Aut_R(\pi_nG)/ \Ru(\pi_0G).
$$
Since $\Ru(\pi_0G)$ is pro-unipotent and $\bG_m$ reductive, the weight decomposition $\bG_m \to \ROut(G)$ lifts to a map $\bG_m \to \Aut_R(\pi_nG)$. The representation theory of $\bG_m$ then gives us a decomposition
$$
\pi_nG =\prod_{i\in \Z} (\pi_nG)^i,
$$
where $ (\pi_nG)^i$ is the subspace of $\pi_nG$ consisting of those $v$ for which $\lambda(v)=\lambda^iv$, for $\lambda \in \bG_m$.
\end{proof}

\begin{corollary}\label{formalpin}
If $(X,\rho)^{\mal}$ is formal, then there is a natural weight decomposition  on the homotopy groups $\varpi_n(X,\rho)^{\mal}$, unique up to conjugation by $\Ru(\varpi_1(X,\rho)^{\mal})$.
\end{corollary}
\begin{proof}
Since $(X,\rho)^{\mal}$ is formal, there is an isomorphism
$$
\Aut_{\Ho(s\cE(R)}(G) \cong \Aut_{\Ho(DG\Alg(R))}(\H^*(X,\bO(R))),
$$
and we define the weight decomposition 
$$
\bG_m \to \Aut_{\Ho(DG\Alg(R))}(\H^*(X,\bO(R)))
$$
by setting $\H^i(X,\bO(R))$ to be of weight $i$.
\end{proof}

\begin{remarks}
\begin{enumerate}
\item If $X$ is a compact K\"ahler manifold, this weight decomposition for homotopy groups extends the weight decompositions of \cite{DGMS} to the non-simply-connected case. Note that if $X$ also satisfies the conditions of Theorem \ref{classicalpimal}, then these pro-algebraic homotopy groups are just $\pi_n(X,x)\ten_{\Z}\R$.

\item In \cite{hodge}, we see that the weight decomposition on the real pro-algebraic homotopy type of a compact K\"ahler manifold $X$ can be extended to an analytic map
$$
\Cx^* \to \Hom_{\Ho(\bS)}(X, \bar{W}G(X)^{\alg}),
$$
giving a real Hodge structure in the sense of \cite{Hodge2}. In fact, it is shown that the pro-algebraic homotopy type of any (simplicial) proper variety has a real Hodge structure in this sense. The $\Cx^*$ action becomes algebraic on taking the Malcev completion relative to the groupoid $R$ defined by the property that  $\Rep(R) \subset \Rep((\varpi_fX)^{\red})$ is the full subcategory  whose objects underlie  variations of Hodge structure. 

\item Even if $X$ is formal, the map $\Aut_{\Ho(s\cE(R))}(G(X,\rho)^{\mal}) \to \Aut_R(\H^*(X,\bO(R))$ need not be injective, since $\H^*(X,\bO(R))$ is not in general cofibrant as a cochain algebra. However, abelian varieties give a  notable class of examples for which $\H^*(X,\bO(R))$ is cofibrant.
\end{enumerate}
\end{remarks}

\begin{proposition}\label{finiteextn}
If $k \subset K$ is a finite extension of fields and $X$ is any simplicial set, then $G(X)^{\alg, K}=G(X)^{\alg,k}\ten_k K$.
\end{proposition}
\begin{proof}
Since $K/k$ is a finite extension, for any pro-algebraic group $H$ over $K$, we may define, by Weil extension of scalars, a pro-algebraic group $H_{K/k}$ over $k$, with the property that $H_{K/k}(A)=H(R\ten_kA)$ for all $k$-algebras $A$. Then, for any $H \in s\agpd_K$, we have
\begin{eqnarray*}
\Hom_{s\agpd_K}(G(X)^{\alg, K}, H)&=& \Hom_{s\gpd}(G(X),H(K))\\
&=& \Hom_{s\gpd}(G(X),H_{K/k}(k))\\
&=& \Hom_{s\agpd_k}(G(X)^{\alg,k}, H_{K/k})\\
&=& \Hom_{s\agpd_K}( G(X)^{\alg,k}\ten_kK, H),
\end{eqnarray*}
as required.
\end{proof}

\begin{corollary}
The complex pro-algebraic homotopy type of a compact K\"ahler manifold is formal.
\end{corollary}

\subsection{Rational formality}

Let $k\subset K$ be an extension of fields in characteristic $0$ and $R$ a reductive pro-algebraic groupoid over $k$. Denote by $R\ten_kK$ the pro-algebraic groupoid over $K$ with the same objects as $R$, and $O(R\ten_kK)(x,y)=O(R)(x,y)\ten_kK$. Fix a representation $\rho\co \pi_fX \to R(k))$, essentially surjective on objects and Zariski-dense on morphisms. Write $\rho_K$ for the corresponding morphism $\rho_K\co \pi_fX \to (R\ten_kK)(K)=R(K)$.  Assume that $X$ satisfies the hypotheses of Theorem \ref{auto}.

\begin{theorem}
If $(X,\rho_K)^{\mal}$ has a weight decomposition $\cW$ such that the induced weight decomposition on $\H^*(X,\bO(R))\ten_kK$ is $k$-rational (i.e. comes from a $\bG_m$-action on  $\H^*(X,\bO(R))$ over $k$), then  there is a (non-canonical) weight decomposition $\cW'$ on    $(X,\rho)^{\mal}$ with the property that $\cW,\cW'$ agree on $\H^*(X,\bO(R))\ten_kK$.
\end{theorem}
\begin{proof}
We adapt the ideas of \cite{Morgan} \S 10.
We have a morphism 
$$
\cW\co \bG_m \to \Aut_{\Ho(s\cE(R)}(G(X,\rho_K)^{\mal} )\cong %\Aut_{\Ho(DG\Alg(R\ten_k K))_0}(\CC^{\bt}(X,\bO(R))\ten_kK)=
\Aut_{\Ho(DG\Alg(R))_0}(\CC^{\bt}(X,\bO(R)))\ten_kK.
$$

Let $H$ be the image of 
$$
\Aut_{\Ho(DG\Alg(R))_0}(\CC^{\bt}(X,\bO(R))) \to \Aut_R(\H^*(X,\bO(R))).
$$  
By hypothesis, the quotient map $\bar{\cW}\co \bG_m \to H\ten_kK$ comes from some map $\bar{\cW}\co \bG_m \to H$.
By Theorem \ref{auto}, we know that $\Aut_{\Ho(DG\Alg(R))_0}(\CC^{\bt}(X,\bO(R))) \to H$ is a unipotent extension. Since $\bG_m$ is reductive, we may therefore lift $\bar{\cW}$ to some map
$$
\cW'\co \bG_m \to \Aut_{\Ho(DG\Alg(R))_0}(\CC^{\bt}(X,\bO(R))),
$$
as required.
\end{proof}

The following is a partial generalisation of \cite{Sullivan} Theorem 12.1, which covers the case $R=1$: 
\begin{corollary}\label{ratformal}
If $(X,\rho_K)^{\mal}$ is formal and $X$ satisfies the hypotheses of Theorem \ref{auto}, then $(X,\rho)^{\mal}$ is formal. 
\end{corollary}
\begin{proof}
The weight decomposition on $(X,\rho)^{\mal}$ gives a weight decomposition on the minimal model $\m$, with $\m_i/[\m_i,\m_i]$ pure of weight $-i-1$. Since $[\m,\m]_i$ is of lower weights, we have a unique choice of homogeneous generators $V_i=\cW_{-i-1}\m_i \subset \m_i$. The differential must preserve the weight decomposition, so $d\co V_{i-1} \to \bigoplus_{a+b=i}[V_{a-1},V_{b-1}]$. Since $V_{i-1} \cong \H^{i}(X,\bO(R))^{\vee}$, Proposition \ref{spectralh} shows that $d$ must be dual to the cup product. Therefore
$$
\m \cong \bar{G}(\H^*(X,\bO(R)),
$$ 
so $(X,\rho)^{\mal}$ is formal. 
\end{proof}

\begin{corollary}
If $X$ is formal over $K$, then $X$ is formal over all subfields $k$ of $K$.
\end{corollary}
\begin{proof}
First observe that formality of $X$ over $(\varpi_fX)^{\red}$ implies formality of $X$ over all of its quotients $R$, since 
$$
\CC^{\bt}(X,\bO(R))=\CC^{\bt}(X, \bO(\varpi_fX^{\red}))\ten^{\varpi_fX^{\red}}O(R).
$$
Now $(\varpi_fX)_k^{\red}\ten_kK$ is a quotient of $(\varpi_fX)^{\red}_K$, since 
$$
\Rep((\varpi_fX)_k^{\red}\ten_kK) \subset \Rep((\varpi_fX)^{\red}_K).
$$
consists of those semisimple $K$-local systems $\vv$ on $X$ with $\vv=\ww\ten_k K$, for $\ww$ a $k$-local system. Taking $R=(\varpi_fX)_k^{\red}\ten_kK$ completes the proof.
\end{proof}

\begin{lemma}
If $K=\varinjlim K_{\alpha}$ is a filtered direct limit of fields, and $G$ is a finitely presented abstract groupoid, then
$$
G^{\alg,K} =\Lim G^{\alg,K_{\alpha}}\ten_{K_{\alpha}} K.
$$
\end{lemma}
\begin{proof}
For any $H \in \agpd_{\Q}$ we have
$$
\Hom_{\gpd}(G,H(K))=\varinjlim \Hom_{\gpd}(G,H(K_{\alpha})),
$$
since $G$ is finitely presented.

Consider the filtered direct system of categories 
$$
\FD\Rep(G^{\alg,K_{\alpha}}) \xra{\ten_{K_{\alpha}}K_{\beta}} \FD\Rep(G^{\alg,K_{\beta}}),
$$
multi-fibred over the filtered direct system
$$
s\Vect(K_{\alpha}) \xra{\ten_{K_{\alpha}}K_{\beta}} s\Vect(K_{\beta}).
$$

Since $\FD\Rep(G^{\alg,K_{\alpha}})$ is just the category of $G$-representations in $K_{\alpha}$-vector spaces, finite presentation then implies that
$$
\varinjlim s\FD\Rep(G(X)^{\alg,K_{\alpha}}) \to s\FD\Rep(G(X)^{\alg,K})
$$
is an equivalence. 

The result now follows by Tannakian duality.
\end{proof}

\begin{proposition}\label{filtered}
If $K=\varinjlim K_{\alpha}$ is a filtered direct limit of fields, and $X$  a finite simplicial set, then 
$$
G(X)^{\alg,K}= \Lim G(X)^{\alg,K_{\alpha}}\ten_{K_{\alpha}} K.
$$
In particular, if $k \subset K$ is a union of finite extensions, then 
$$
G(X)^{\alg,K}=G(X)^{\alg,k}\ten_k K,
$$
so homotopy types over $\Q$ and $\bar{\Q}$ are the same.
\end{proposition}
\begin{proof}
Since $G(X)$ is finitely presented, $G(X)_n$ is finitely presented for all $n$, and $\pi_fG(X)$ is also finitely presented. Therefore
$$
(G(X)_n)^{\alg,K}=\Lim (G(X)_n)^{\alg,K_{\alpha}}\ten_{K_{\alpha}} K, \quad  \pi_fG(X)^{\alg,K}=\Lim \pi_fG(X)^{\alg,K_{\alpha}}\ten_{K_{\alpha}} K.
$$

As $(G(X)^{\alg})_n$ is the pro-unipotent  completion of $ (G(X)_n)^{\alg}\to \pi_fG(X)^{\alg}$, the result follows. 
\end{proof}

\begin{corollary}
Assume that $K=\varinjlim K_{\alpha}$ and $X$ is a finite simplicial set, with weight decompositions on each $G(X)^{\alg,K_{\alpha}}$, such that the induced weight decompositions on cohomology are compatible with the maps
$$
\H^*(X, \bO(G(X)^{\red,K_{\beta}})) \to\H^*(X, \bO(G(X)^{\red,K_{\alpha}}))\ten_{K_{\alpha}}K_{\beta},
$$
for every field extension $K_{\alpha} \subset K_{\beta}$.  Then there is a weight decomposition on 
$$
G(X)^{\alg,K},
$$
compatible on cohomology with the weight decompositions on the $\H^*(X, \bO(G(X)^{\red,K_{\alpha}}))$.
\end{corollary}
\begin{proof}
Let $G=G(X)$. By Proposition \ref{filtered}, 
\begin{eqnarray*}
&\im(\Aut(G^{\alg,K})_{G^{\red,K}} \to \Aut_{G^{\red,K}}(\H^*(X,\bO( G^{\red,K}))))\\
= &\Lim \im(\Aut(G^{\alg,K_{\alpha}})_{G^{\red,K_{\alpha}}} \to \Aut_{G^{\red,K_{\alpha}}}(\H^*(X,\bO( G^{\red,K_{\alpha}}))))\ten_{K_{\alpha}}K.
\end{eqnarray*}
The hypotheses give us a morphism from $\bG_m$ to the latter group, and hence to the former group. Theorem \ref{auto} now completes the proof.
\end{proof}

\begin{corollary}\label{kformalall}
The pro-algebraic homotopy type of a compact K\"ahler manifold is formal over any field $K$ of characteristic zero.
\end{corollary}
\begin{proof}
We may write $K=\varinjlim K_{\alpha}$, with $K_{\alpha}$ of finite transcendence degree over $\Q$. Therefore $K_{\alpha}$ is isomorphic to a subfield of $\Cx$, so $G(X)^{\alg,K_{\alpha}}$ is formal, using Theorem \ref{kformal}. This gives us a weight decomposition on $G^{\alg,K}$, and the proof of Corollary \ref{ratformal} applies.
\end{proof}

\bibliographystyle{alphanum}
\addcontentsline{toc}{section}{Bibliography}
\bibliography{references.bib}
\end{document}